\documentclass[onecolumn,draftcls]{IEEEtran}  

\usepackage{amsmath}  % assumes amsmath package installed
\usepackage{amssymb}  % assumes amsmath package installed

\usepackage{url}
\usepackage{graphicx}
\usepackage{enumerate}
\usepackage{cite}

\newtheorem{thm}{Theorem}[section] 
\newtheorem{lem}[thm]{Lemma}
\newtheorem{defn}{Definition}[section] 
\newtheorem{prop}[thm]{Proposition}

\newcounter{taskcounter}
\newenvironment{hypot}{ \refstepcounter{taskcounter}\textbf{A\arabic{taskcounter}}:}
{}

\newenvironment{ohypot}{ \refstepcounter{taskcounter}\textbf{A\arabic{taskcounter}'}:}
{}

\newcommand {\ass}[1]{\textbf{A\ref{#1}}}
\newcommand {\oass}[1]{\textbf{A\ref{#1}'}}

\renewcommand \t{^\top}
\renewcommand \b{\mathbf}
\newcommand \tr{\mathrm{Tr}}
\newcommand \teq{\triangleq}

\newcommand \be{\begin{equation}}
\newcommand \ee{\end{equation}}
\newcommand \bml{\begin{multline}}
\newcommand \eml{\end{multline}}
\newcommand \ba{\begin{align}}
\newcommand \ea{\end{align}}

\newcommand \nqquad{\mspace{-36.0mu}}
\newcommand \qquadtwo{\qquad\qquad}

\newcommand \qquadeight{\qquad\qquad\qquad\qquad\qquad\qquad\qquad\qquad}

\newcommand \R{\mathbb{R}}			% The reals
\newcommand \N{\mathbb{N}}			% The naturals
\newcommand \Z{\mathbb{Z}}			% The integers
\newcommand \B{\mathcal{B}}			% Borel sets
\newcommand \Ne{\mathcal{N}}			% Neighborhood
\newcommand \F{\mathcal{F}}			% A sigma algebra F
\newcommand \probs{(\Omega,\mathcal{F},P)}	% probability space
\newcommand \U{\mathcal{U}}			% control set
\newcommand \W{\Omega}				% Capital omega
\newcommand \w{\omega}				% Lowercase omega
\renewcommand \iff{\Leftrightarrow}		% If and only if
\newcommand \conv{\rightarrow}			% Converges
\newcommand \e{\epsilon}			% Epsilon
\newcommand \p{\partial}			% Partial

\newcommand \E{\mathbb{E}}			% Expected Value
\newcommand \f{\infty}				% infty
\DeclareMathOperator*{\argmin}{arg\,min}

\newcommand \limt{\lim_{t \rightarrow \infty}}
\newcommand \limT{\lim_{T \rightarrow \infty}}

\newcommand \limN{\lim_{N \rightarrow \infty}}

\newcommand \lsT{\limsup_{T \rightarrow \infty}}

\newcommand \lsN{\limsup_{N \rightarrow \infty}}

\newcommand \xno{{(x^{N})}^{0}}

\renewcommand \th{\mathbf{\theta}}
\newcommand \Th{\mathbf{\Theta}}
\newcommand \hth{\hat{\mathbf{\theta}}}
\newcommand \z{\zeta}
\newcommand \hz{\hat{\zeta}}
\newcommand \hzn{\hat{\zeta}^{N_0}}
\newcommand \uo{u^{0}}
\newcommand \huo{\hat{u}^{0}}
\newcommand \uio{u_i^{0}}
\newcommand \huio{\hat{u}_i^{0}}
\newcommand \xo{x^{0}}
\newcommand \hxo{\hat{x}^{0}}
\newcommand \xio{x_i^{0}}
\newcommand \hxio{\hat{x}_i^{0}}
\newcommand \umio{u_{-i}^{0}}
\newcommand \humio{\hat{u}_{-i}^{0}}

\newcommand \mno{{(m^{N})}^{0}}
\newcommand \hmno{{(\hat{m}^{N})}^{0}}
\newcommand \mnu{m^N_{|_{u_i,u_{-i}^0}}}
\newcommand \hmnu{m^N_{|_{u_i,\hat{u}_{-i}^0}}}

\newcommand \ps{p^s}

\renewcommand{\IEEEQED}{\hfill \IEEEQEDclosed}

\title{Mean Field Stochastic Adaptive Control}

\author{Arman~C. Kizilkale~and~Peter~E. Caines% <-this % stops a space
\thanks{Arman~C. Kizilkale and Peter~E. Caines are with the Department of Electrical and Computer Engineering and the Centre for Intelligent Machines, McGill University, Montreal, Canada.
        {\tt\small \{arman;peterc\}@cim.mcgill.ca}}%
}

\begin{document}

\maketitle

%===============================================================================
\begin{abstract}
For noncooperative games the mean field (MF) methodology provides decentralized strategies which yield Nash equilibria for large population systems in the asymptotic limit of an infinite (mass) population. The MF control laws use only the local information of each agent on its own state and own dynamical parameters, while the mass effect is calculated offline using the distribution function of (i) the population's dynamical parameters, and (ii) the population's cost function parameters, for the infinite population case. These laws yield approximate equilibria when applied in the finite population. 

In this paper, these a priori information conditions are relaxed, and incrementally the cases are considered where, first, the agents estimate their own dynamical parameters, and, second, estimate the distribution parameter in (i) and (ii) above.

An MF stochastic adaptive control (SAC) law in which each agent observes a random subset of the population of agents is specified, where the ratio of the cardinality of the observed set to that of the number of agents decays to zero as the population size tends to infinity. Each agent estimates its own dynamical parameters via the recursive weighted least squares (RWLS) algorithm and the distribution of the population's dynamical parameters via maximum likelihood estimation (MLE). Under reasonable conditions on the population dynamical parameter distribution, the MF-SAC Law applied by each agent results in (i) the strong consistency of the self parameter estimates and the strong consistency of the population distribution function parameters; (ii) the long run average $L^2$ stability of all agent systems; (iii) a (strong) $\e$-Nash equilibrium for the population of agents for all $\e>0$; and (iv) the a.s.\ equality of the long run average cost and the non-adaptive cost in the population limit.

\end{abstract}


%\tableofcontents

%===============================================================================
\section{Introduction}\label{sec:intro}

\textit{Overview}

The control and optimization of large-scale stochastic systems is evidently of importance due to their ubiquitous appearance in engineering, industrial, social and economic settings. The complexity of these problems is amplified by the fact that for many such systems the agents involved have conflicting objectives; hence, it is appropriate to consider optimization methodologies based upon individual payoffs or costs. In particular, game theory has been formulated to capture such individual interest seeking behaviour of the agents in many social, economic and manmade systems. However, in a large-scale dynamic model, this approach results in an analytic complexity which is in general prohibitively high, and correspondingly leads to few substantive dynamic optimization results.

The optimization of large-scale linear control systems wherein (i) many agents are coupled with each other via their individual dynamics, and (ii) the costs are in an ``individual to the mass'' form was presented in \cite{2003HCM_CDC,2007HCM_TAC} where the theory of mean field (MF) control (previously termed Nash Certainty Equivalence) was introduced. It is to be noted that the dynamic large-scale cost coupled optimization structure of \cite{2007HCM_TAC} is motivated by a variety of scenarios, for instance, those analysed in \cite{2000HFV_N,2003TJP_CDC,2004LP_TAC,2000Lo_N}.

In the literature, studies of stochastic dynamic games and team problems may be traced to the 1960s (see e.g. \cite{1965HBB_TAC, 1967Va_SIAM, 1975Wi}) while within the optimal control context weakly interconnected systems were studied in \cite{1988Be}, and in a two player noncooperative nonlinear dynamic game setting Nash equilibria were analysed in \cite{1991SB}, where the coefficients for the coupling terms in the dynamics and costs are required to be small. In contrast to these studies, games with large populations are analyzed in \cite{1964Au,1988JR,2007HCM_TAC}. In \cite{2007HCM_TAC} the $\e$-Nash equilibrium properties are analysed for a system of competing agents where individual control laws use local information and the average effect of all agents taken together, henceforth referred to as \emph{the mass}. Overall, the MF methodology for noncooperative LQG games with mean field coupling has been developed in \cite{2003HCM_CDC,2007HCM_TAC,2009Ca} providing decentralized strategies which yield Nash equilibria. A nonlinear extension using McKean-Vlasov Markov process models is also presented in \cite{2006HMC_CIS}.

The central notion of MF theory is that for general classes of large population stochastic dynamic games there exist game theoretic Nash equilibria for the individual agents when each applies certain competitive strategies (i.e. control laws) with respect to the mass effect resulting from all the agents' strategies. Here each agent is modelled by an individually controlled stochastic system and the systems interact through their individual cost functions and possibly via weak dynamical interaction. The key feedback nature of the mean field solutions is that the individual competitive actions against the mass, plus local feedback control, act so as to collectively reproduce that mass behaviour. The mass effect and associated feedback control laws are calculated offline for the infinite population case and yield approximate equilibria when applied in the finite population case.
% (i.e. the overall effect of the resulting population behaviour)

For this class of game problems, a related approach has been independently developed in \cite{2008WBR,2010AJWG_CDC}, where the notion of oblivious equilibrium by use of a mean field approximation for models of many firm industry dynamics is proposed. The asymptotic equilibrium properties of a market with a large population of agents is studied in \cite{2011Bo}. Another related work is \cite{2007LL} where a mean field Nash equilibrium is studied subject to the assumed existence of a factorizing mean field distribution corresponding to the propagation of chaos for the infinite population system. The work in \cite{2009TLEA} presents mean field control results for a Markov Decision Problem (MDP) formulation of evolutionary games and teams where the basic system hypothesis is the exchangeability of the underlying random processes.

\textit{Stochastic Adaptive Control}

For discrete time dynamics the long run average (LRA) asymptotically optimal adaptive tracking problem was solved in \cite{1981GRC_SIAM}; subsequently, it was shown in \cite{1984CL_TAC} that strongly consistent parameter estimates may be obtained by the use of persistently excited controls. %v3 Results for continuous time systems have been obtained in \cite{1991GGW_SIAM,1991GGMW,1990DD,1996CDD}. 
The LRA stochastic (sample path) mean square stability for continuous time linear stochastic adaptive systems was established in \cite{1992Ca_SCL}. The weighted least squares (WLS) scheme introduced in \cite{1995Be_SIAM} was shown in \cite{1996Gu_TAC} to be convergent without stability and excitation assumption, and a LRA asymptotically optimal solution to the continuous time adaptive LQG control problem under controllability and observability assumptions using the WLS scheme for identification was subsequently obtained in \cite{1999DGD_TAC} following \cite{1986CG_IJC,1991CG,1996Ga} and \cite{1997GD_SCL}.
%v3
%In order to relax the persistent excitation hypothesis,
%under a persistent excitation hypothesis \cite{1982LW}

% ACK-LV
% In the standard setting of single agent adaptation, the long run average (LRA) linear stochastic (sample path) mean square stability for continuous time linear stochastic adaptive systems was established in \cite{1992Ca_SCL} by use of the recursive least squares (RLS) algorithm in the adaptive control algorithm under a persistent excitation hypothesis \cite{1982LW} using the differentiation technique of \cite{1987CM_TAC} and \cite{1990CG_SIAM}.
% ACK-LV
%Results for continuous time parameter partially observed systems have been obtained in \cite{1991GGW,1991GGMW}.

\textit{MF Stochastic Adaptive Control}

It is important to note that in the non-adaptive MF theory \cite{2003HCM_CDC, 2007HCM_TAC} each agent uses its \emph{self state} and \emph{self dynamical parameters} (i.e. its own state and its own dynamical parameters) and statistical information on the dynamical parameters of the population in order to generate the control action. %v3 The relaxation of the requirements of a priori known information above naturally leads to the use of the methods of stochastic adaptive control in the MF stochastic control context. %while the behaviour of the mass is precomputable from knowledge of the distribution of the dynamical parameters throughout the mass population. All this information is assumed known to each agent in the basic non-adaptive MF theory.
The natural initial problem in the development of adaptive MF stochastic system theory is that where each agent needs to estimate its own dynamical parameters, while its control actions are permitted to be explicit functions of the parameter distribution of the entire population of competing agents \cite{2010KC_MTNS}. Subsequent problem generalizations are such that (i) each agent also needs to estimate the distribution parameter of the population's dynamical parameters \cite{2010KC_CDC}, and (ii) cost function parameters also vary over the population and this distribution parameter is unknown to each agent and hence needs to be estimated \cite{2010KC_ACCA}. In this paper we provide a solution to the most general problem in this sequence.

The inclusion of learning procedures for the identification by a given agent of the  dynamical and cost function parameters of other competing agents in a stochastic dynamic system, or of the statistical distribution of these parameters in a mass of competing agents, introduces new features into the system theoretic MF setup. In this connection we note that in the economics literature the so-called ``privacy of information'' on dynamical parameters and cost function parameters is an important issue \cite{2002Kr,1989MS_ET,1999MM}. %v3 In many micro-economic domains, for instance in the theory of market mechanisms and of auctions, it is a significant aspect of the dynamical models that the valuation function or utility function of an agent should not initially be known to other agents in the system; correspondingly there is a vast literature concerned with the identification of agents' valuation functions and hidden dynamics \cite{1989MS_ET,1999MM}.
%  \cite{2001NR} \cite{2002Kr}

This paper presents an MF stochastic adaptive control (SAC) law in which each agent observes a random subset of the population of agents. The MF-SAC Law specifies that the ratio of the cardinality of the observed set of agents to that of the population of agents is chosen so that it decays to zero as the population size tends to infinity. When the MF-SAC Law is applied by each member of the population, each agent estimates its self dynamical parameters via the recursive weighted least squares (RWLS) algorithm and the distribution of the population's dynamical parameters via maximum likelihood estimation (MLE). 

Under reasonable conditions on the population dynamical parameter distribution, the MF-SAC Law results in (i) the strong consistency of the self parameter estimates and the strong consistency of the population distribution function parameters; (ii) the long run average $L^2$ stability of all agent systems; (iii) a (strong) $\e$-Nash equilibrium for the population of agents for all $\e>0$; and (iv) the a.s.\ equality of the long run average cost and the non-adaptive cost in the population limit.

%v3 Finally, we present the simulation results for a population of 400 agents, where each agent estimates its self dynamical parameters and the population dynamical and cost function distribution parameter based on observations on the trajectories and control actions of self and 20 randomly chosen competing agents.

%%%%%%%%%%%%%%%%%%%%%%%%%%%%%%%%%%%%%%%%%%%%%%%%%%%%%%%%%%%%%
\textit{Notation}

%v3Henceforth, we adopt the following notation. 
We denote the set of nonnegative real numbers by $\R_+$, the set of nonnegative integers by $\Z_+$, and the set of strictly positive integers by $\Z_1$. The norm $\lVert \cdot \rVert$ denotes the 2-norm of vectors and matrices, and $\lVert x  \rVert_Q^2 \triangleq x\t\b{Q}x$. $\b{C}_b=\{ x:x\in \b{C},\sup_{t\geq 0} \lVert x(t) \rVert < \f\}$ denotes the family of all bounded continuous functions, and for any $x\in\b{C}_b$, $\lVert \cdot \rVert_\f$ denotes the supremum norm: $\lVert x \rVert_\f\teq \sup_{t\geq 0} \lVert x(t) \rVert$. $\tr(\b{X})$ denotes the trace, and $\b{X}\t$ denotes the transpose of a matrix $\b{X}$. %Moreover, $\lambda_{max}$ and $\lambda_{min}$ denote the maximum and minimum eigenvalues of symmetric matrices respectively.

\section{Problem Formulation and MF-SAC Law Specification}\label{sec:pro_form}

\subsection{Review of Non-Adaptive MF Stochastic Control}\label{ssec:NCESC}

We consider a large population of $N$ stochastic dynamic agents which (subject to independent controls) are stochastically independent, but which shall be cost coupled, where the individual dynamics are defined by
\begin{equation}\label{eqn:dynamics}
dx_i=[\b{A}_i x_i + \b{B}_i u_i]dt + \b{D} dw_i,\quad t\geq 0,\quad 1\leq i \leq N,
\end{equation} 
where for agent $A_i$, $x_i \in \R^n$ is the state, $u_i \in \R^m$ is the control input, $w_i \in \R^r$ is a standard Wiener process on a sufficiently large underlying probability space $\probs$ such that $w_i$ is progressively measurable with respect to $\F^{w_i}\triangleq \{\F_t^{w_i};t\geq 0\}$. We denote the state configuration by $x=(x_1, \cdots, x_N)\t$, and (with an abuse of notation) the population average state by $x^N=(1/N)\sum_{i=1}^N x_i$. 

The long run average (LRA) cost function for the agent $A_i,1\leq i\leq N$, is given by
\begin{equation}\label{eqn:cost}
J_i^N(u_i,u_{-i}) = \lsT\frac{1}{T}\int_0^T
\{ \lVert x_i - m^N \rVert_{Q_i}^2 + \lVert u_i \rVert_{R}^2 \}dt,
\end{equation}
w.p.1, where we assume the cost-coupling to be of the form $m^N(t) \triangleq m( x^N (t) + \eta),\, \eta \in \R^n$. The coefficients $\th_i\t \triangleq [\b{A}_i,\b{B}_i, \b{Q}_i]\in \Th \subset \R^{n(n+m+(n+1)/2)},$ will be called the \emph{dynamical} and \emph{cost function} parameters. The disturbance weight matrix $\b{D}$ and the control action penalizing matrix $\b{R}$ are constant matrices, which are assumed to be known by all agents, and assumed to be the same for all agents in the population. The choice of homogeneous parameters for $\b{D}$ and $\b{R}$ is only for notational brevity; the analysis is similar for varying $\b{D}$ and $\b{R}$. The function $u_i(\cdot)$ is the control input of the agent $A_i$ and $u_{-i}$ denotes the control inputs of the complementary set of agents $A_{-i} = \{ A_j, j \neq i,1\leq j\leq N\}$. 

For the basic MF control problem, the following assumptions are adopted.

\begin{hypot}\label{ass:ind}
The disturbance processes $ w_i,1\leq i \leq N $, are mutually independent and independent of the initial conditions, and $\sup_{i\geq 1} [\tr{\Sigma_i} + \E\lVert x_i(0) \rVert^2 ] < \infty$, where $\E w_iw_i\t = \Sigma_i,\, 1\leq i \leq N$.
\IEEEQED
\end{hypot}

\begin{hypot}\label{ass:conobs}
$\tilde{\Th}$ is an open set such that for each $\th\t= [\b{A}_{\tilde{\theta}},\b{B}_{\tilde{\theta}},\b{Q}_{\tilde{\theta}}] \in \b{\tilde{\Theta}}$, $[\b{A}_{\tilde{\theta}},\b{B}_{\tilde{\theta}}]$ is controllable and $[\b{Q}_{\tilde{\theta}}^{1/2},\b{A}_{\tilde{\theta}}]$ is observable.%, $\b{Q}_{\tilde\theta}$ is symmetric and positive definite, and $\b{B}_{\tilde{\theta}}$ is invertible (and hence, necessarily, $[\b{A}_{\tilde{\theta}},\b{B}_{\tilde{\theta}}]$ is controllable).
\IEEEQED
\end{hypot}

\begin{hypot}\label{ass:compact}
Let the parameter set $\Th$ be a compact set such that $\Th \subset \tilde{\Th} \subset \R^{n(n+m+(n+1)/2)}$, and\\ $ \lVert \b{R}^{-1} \rVert \lvert \gamma \rvert \int_{\th \in \Th} \lVert \b{Q}(\th) \rVert \lVert \b{B}(\th) \rVert^2    ( \int_0^\infty \lVert e^{\b{A}_*(\th)\tau} \rVert d\tau )^2 dF_\z(\th)<1$, where $\b{A}_* = \b{A}-\b{B}\b{R}^{-1}\b{B}\t\b{\Pi}$, $\z$ is the distribution parameter and $\gamma$ is defined in the next hypothesis.
% ACK, 29 Mar 2011
% Let the parameter set $\Th$ be a compact set such that $\Th \subset \tilde{\Th} \subset \R^{n(n+m+(n+1)/2)}$, and $ \lVert \b{R}^{-1} \rVert  \gamma M_Q M_B^2 M_{A_*}^2 < 1$, where $M_Q = \sup_{\th\in\Th}\lVert \b{Q}_{\theta} \rVert,\, M_B = \sup_{\th\in\Th}\lVert \b{B}_{\theta} \rVert,\, M_{A_*} = \sup_{\th\in\Th}\left\lVert  \int_0^\infty \lVert e^{\b{A}_*(\th)\tau} \rVert d\tau  \right\rVert$. $\gamma$ is defined in the next hypothesis.
\IEEEQED
\end{hypot}

\begin{hypot}\label{ass:cost_coup}
The cost-coupling is of the form: $m^N(\cdot) \triangleq m((1/N)\sum_{k=1}^N x_k + \eta),\eta\in\R^n$, where the function $m(\cdot)$ is Lipschitz continuous on $\R^n$ with a Lipschitz constant $\gamma>0$, i.e. $\lVert m(x)-m(y) \rVert \leq \gamma \lVert x - y \rVert$ for all $x, y \in \R^n.$ 
\IEEEQED
\end{hypot}

% The cost coupling function $m^N(\cdot)$ is estimated by a deterministic function $x^*(t), t\geq 0$, and the problem is first solved in \cite{2007HCM_TAC} for a discounted expected cost function when $m^N(\cdot)$ is substituted for $x^*(t),\,t\geq 0$. Following \cite{2007HCM_TAC}, the problem is solved in \cite{2008LZ_TAC} for the cost function \eqref{eqn:cost}. The positive solution is obtained for the following algebraic Riccati equation

For dynamics \eqref{eqn:dynamics} and cost function \eqref{eqn:cost}, a production output planning example is provided in \cite{2007HCM_TAC} that satisfies the assumptions given above. Each agent's production level $x_i$ is modeled by \eqref{eqn:dynamics}, and each agent's cost function is of tracking type \eqref{eqn:cost}, where the tracked signal is a function of price, which is an averaging function of production levels: $m^N(t) \triangleq m( x^N (t) + \eta),\, \eta \in \R^n$.

Following \cite{2007HCM_TAC}, the long run average (LRA) mean field (MF) problem is formulated in \cite{2008LZ_TAC}. Each agent $A_i,\, 1\leq i \leq N$, obtains the positive definite solution to the algebraic Riccati equation 
\begin{equation} \label{eqn:DE_Pi} 
\b{A}_i\t\b{\Pi}_i + \b{\Pi}_i \b{A}_i-\b{\Pi}_i
\b{B}_i\b{R}^{-1}\b{B}_i\t \b{\Pi}_i + \b{Q}_i = 0.
\end{equation}
Moreover, for a given \emph{mass tracking signal} $x^*\in \b{C}_b[0,\infty)$ the mass offset function $s_i(t)$ is generated by the differential equation
\begin{equation}\label{eqn:DE_s}
-\frac{ds_i(t)}{dt} = \b{A}_i\t s_i(t) - \b{\Pi}_i \b{B}_i\b{R}^{-1}\b{B}_i\t s_i(t) - \b{Q}_i x^*(t), \quad t\geq 0.
\end{equation}
Then, the optimal tracking control law \cite{1992Be} is given by
\begin{equation}\label{eqn:opt_tracking} 
u_i(t)=-\b{R}^{-1}\b{B}_i\t(\b{\Pi}_i x_i(t) + s_i(t)), \quad t\geq 0,
\end{equation}
where $u_i(\cdot)$ solves $\inf J_i(u_i,x^*)$, which is defined below by an abuse of notation: 
\[
J_i(u_i,x^*) \triangleq \lsT\frac{1}{T}\int_0^T
\{ \lVert x_i - x^* \rVert_{Q_i}^2 + \lVert u_i \rVert_{R}^2 \}dt \quad \text{w.p.1}.
\]
Note that the procedure above assumes a given mass tracking signal $x^*$. The equation system to calculate $x^*$ will be given subsequently.

We first define the empirical distribution associated with the first $N$ agents:\\ $F_\z^N(\th)=\frac{1}{N}\sum_{i=1}^N \mathbb{I}_{(\th_i < \th)}$, $\th\in\R^{n(n+m+(n+1)/2)}$, where $\{\th_i,\, 1\leq i \leq N\}$ is a set of random matrices on $\probs$ with the probability distribution $F_\z(\th)$, parameterized by $\z\in P \subset \tilde P \subset \R^p$, the \emph{population dynamical and cost function distribution parameter} such that $P$ is compact and $\tilde P$ is an open set. Then we employ the following assumption.

\begin{hypot}\label{ass:weakdist}
There exists a family of distribution functions $\{F_\z(\th);\, \th\in\Th\},\z\in \tilde P$, such that $F_\z^N(\cdot)\conv F_\z(\cdot)$ w.p.1 weakly on $\th \in \Th$ and uniformly over $\z \in P$ as $N\rightarrow \infty$.%, i.e., $\lim_{N \rightarrow \infty } F_\z^N (\cdot)=F_\z(\cdot)$.
% first submission 
% if $F_\z$ is continuous at $\th\in\Th \subset \R^{n(n+m+(n+1)/2)}$.
% first submission
\IEEEQED
\end{hypot}

% \ass{ass:weakdist} is satisfied by the Glivenko-Cantelli theorem \cite{1997CT}. 
Each agent solves the equation system below to calculate the mass tracking signal $x^*(\tau,\z),\, t_0\leq\tau<\infty$, offline, for an infinite population of agents.
\begin{defn}\label{defn:MF_Equations}
\emph{Mean Field (MF) Equation System on $[t_0,\infty)$}:
\begin{equation}\label{eqn:ncee}
\begin{aligned}
-\frac{ds_\theta}{d\tau} & = (\b{A}_\theta\t - \b{\Pi}_\theta \b{B}_\theta\b{R}^{-1}\b{B}_\theta\t ) s_\theta - \b{Q}_\theta x^*(\tau,\z),\\
\frac{d\bar x_\theta}{d\tau} & = (\b{A}_\theta-\b{B}_\theta\b{R}^{-1}\b{B}_\theta\t \b{\Pi}_\theta) \bar{x}_\theta -\b{B}_\theta\b{R}^{-1}\b{B}_\theta\t s_\theta ,\\
\bar x(\tau,\z) & = \int_\Theta \bar x_\theta d F_\z(\th),\\
x^*(\tau,\z) & = m(\bar x(\tau,\z)+\eta),\quad t_0 \leq \tau < \infty.
\end{aligned}
\end{equation}
\IEEEQED
\end{defn}

Under \ass{ass:ind}-\ass{ass:cost_coup}, the MF Equation System admits a unique bounded solution \cite{2007HCM_TAC}.
% under the following technical assumption:
% \begin{hypot}\label{ass:fpt}
% \begin{equation}
% \lVert \b{R}^{-1} \rVert  \gamma\int_{\th \in \Th} \lVert \b{Q}(\th) \rVert \lVert \b{B}(\th) \rVert^2    ( \int_0^\infty \lVert e^{\b{A}_*(\th)\tau} \rVert d\tau )^2 dF_\z(\th)<1. 
% \end{equation}
% \end{hypot}

\textit{The Global Observation Control Set $\U_g^N$:\,} For the optimality analysis, we first introduce the global observation control set. The set of control inputs $\U_g^N$ consists of all feedback controls adapted to $\{ \th_j, 1\leq j \leq N ; \, F_\z(\th);\, \F_t^N, t\geq 0 \}$, where $\F_t^N$ is the $\sigma$-field generated by the set $\{ x_j(\tau);\, 0\leq \tau \leq t, 1\leq j\leq N\}$.

\textit{The Local Observation Control Set $\U_{l,i}^N$:\,} The local observation control set of agent $A_i$ is the set of control inputs $\U_{l,i}^N$ which consists of the feedback controls adapted to the set $\{ \th_i;\, F_\z(\th);\, \F_{i,t}, t\geq 0 \}$. The $\sigma$-field $\F_{i,t}$ is generated by $( x_i(\tau);\, 0\leq \tau \leq t)$,  and  $\F_t^N$ is the $\sigma$-field generated by the set $\{ x_j(\tau);\, 0\leq \tau \leq t, 1\leq j\leq N\}$.

%ACK-LV
% \begin{defn}\label{defn:enash}
% Given $\e>0$, the set of controls $\U^0 = \{ {u}_i^0;1\leq i \leq N \}$ generates an \emph{$\e$-Nash Equilibrium} w.r.t.\ the costs $\{ J_i; 1\leq i\leq N \}$ if for each $i,\, 1\leq i \leq N,$
% \begin{align*}
% J_i^N(u_i^0, u_{-i}^0) - \e \leq \inf_{u_i \in \U_g^N} J_i^N(u_i,u_{-i}^0)\leq J_i^N(u_i^0,u_{-i}^0). 
% \end{align*}
% \IEEEQED
% \end{defn}
%ACK-LV

\begin{thm}\label{thm:Main_MF_thm} \emph{Non-Adaptive MF Stochastic Control (SC) Theorem} \cite[following \cite{2007HCM_TAC}]{2008LZ_TAC}

Let \ass{ass:ind}-\ass{ass:weakdist} hold. The MF Stochastic Control Law \eqref{eqn:opt_tracking} generates a set of controls  $\U_{MF}^N \triangleq \{ {u}_i^0;1\leq i \leq N \},\, 1\leq N < \infty,$ with
\begin{equation}\label{eqn:MF_control} 
u_i^0(t)=-\b{R}^{-1}\b{B}_i\t(\b{\Pi}_i x_i(t) + s_i(t)),\quad t\geq 0,
\end{equation}
such that
\begin{enumerate}[(i)]
\item the MF equations \eqref{eqn:ncee} have a unique solution;
\item all agent system trajectories $x_i,\, 1\leq i \leq N,$ are $LRA-L^2$ stable w.p.1; 
\item  $\{ \U_{MF}^N;1\leq N < \infty\}$ yields an $\e$-Nash equilibrium for all $\e>0$, i.e., for all $\e>0$, there exists $N(\e)$ such that for all $N \geq N(\e)$
\begin{equation*} 
J_i^N({u}_i^0, {u}_{-i}^0)-\e \leq\inf_{u_i \in\U_g^N } J_i^N(u_i, {u}_{-i}^0) \leq J_i^N({u}_i^0, {u}_{-i}^0).
\end{equation*}

\end{enumerate}

\IEEEQED
\end{thm}

Conceptually, Theorem \ref{thm:Main_MF_thm} may be paraphrased to say that individual competitive actions against the mass effect collectively produce the mass behaviour, and hence the $\e$-Nash equilibrium is obtained. In the proof of Theorem \ref{thm:Main_MF_thm}, the results are first established for an infinite population and then are shown to be approximated by a large finite population with the approximation error decaying to zero as the population size goes to infinity; it is this which gives the $\e$-Nash property.

% second submission version
%===============================================================================
\subsection{MF Stochastic Adaptive Control (SAC)}\label{ssec:ncesac}

In this section we first present the identification schemes to be used by each agent under the MF Stochastic Adaptive Control (SAC) Law to estimate both the self dynamical parameters and the population dynamical and cost function distribution parameter. In other words, the analysis concerns a family of agents $A_i,\, 1\leq i \leq N,$ whose control action at any instant is not permitted to be an explicit function of the self dynamical parameters $[\b{A}_i,\b{B}_i]$ and the dynamical and cost function distribution parameter $\z$. At time $t\geq 0$, the self dynamical parameters are estimated from the input-output sample path $\{ x_i(\tau), u_i(\tau); 0\leq \tau \leq t \}$ of $A_i$; in other words, each agent $A_i$ performs the identification based upon observations of its own trajectory. The distribution parameter $\z$ is estimated from observations $\{ x_j(\tau), u_j(\tau); 0\leq \tau \leq t, j\in Obs_i(N) \}$ on a random subset of agents $Obs_i(N)$ where $\lvert Obs_i(N) \rvert \conv \infty$, and $\lvert Obs_i(N)\rvert/N \conv 0$ as $N \conv \infty$. 

\textit{The Adaptive Agent Control Set $\U_{a,i}^N$:\,} We next define the set of control inputs $\U_{a,i}^N$, the admissible control set of an adaptive agent $A_i$, which consists of all feedback controls adapted to the set $\{ \F_{i,t} , \, \F_{i,t}^{obs} ,\,  t\geq 0;\, \b{Q}_i\}$. The $\sigma$-field $\F_{i,t}$ is generated by the agent's own trajectory and control input, $\{ x_i(\tau),u_i(\tau);\, 0\leq \tau \leq t\}$, and $\F_{i,t}^{obs},\, t\geq 0$, is the observation $\sigma$-field generated by the trajectories and control inputs in the set $Obs_i(N)$, $\{ x_j(\tau), u_j(\tau);\, 0\leq \tau \leq t,\, j\in Obs_i(N)\}$. For definiteness in this paper, the identification algorithms employed are recursive weighted least squares (RWLS) for the self dynamical parameter identification and maximum likelihood estimation (MLE) for the distribution parameter identification. However, any identification scheme which generates consistent estimates w.p.1 (subject to the given hypotheses) will also yield the system asymptotic equilibrium properties to be established.

%===============================================================================
\subsubsection{Self Dynamical Parameter Identification (SDPI)}\label{sssec:sel_dyn_par_est}

We denote the self estimate of the matrix $\th_i$ by $\hth_{i,t}\t= [\hat{\b{A}}_{i,t},\hat{\b{B}}_{i,t},\b{Q}_i]$, $t\geq 0$, and the estimate of $\z$ by $\hzn_{i,t},\, t\geq 0$, where $N_0 \teq \lvert Obs_i(N) \rvert$, and assume $\hth_{i,t}$ and $\hzn_{i,t}$ are generated at $t\geq 0$ by the identification algorithm. Note that the self cost function parameter $\b{Q}_i$ is in the information set of agent $A_i$, and is therefore not to be estimated. We adopt the notation $\z^0\triangleq \z,\, \th^0 \triangleq \th$ for the true parameters in the system. At time $t\geq 0$, agent $A_i$ solves the RWLS equations with the measurement variable set as $dx_t$ with the regression vector $[x_t\t, u_t\t]$ in order to obtain the estimates $[\hat{\b{A}}_{i,t},\hat{\b{B}}_{i,t}]$. To ensure controllability and observability of the estimates, a projection method is used; the estimates are projected onto the compact set $\Th_{|Q_i}\subset \tilde{\Th}_{|Q_i}$, where given $\b{Q}_i$, $[\b{A}_\theta,\b{B}_\theta]$ is controllable and $[\b{Q}_i^{1/2},\b{A}_\theta]$ is observable. Note that $\Th$ is known to all agents in the system. 

\subsubsection{Population Dynamical and Cost Function Distribution Parameter Identification} \label{sssec:pop_par_est}

\paragraph{Population Dynamical Parameter Identification (PDPI)} 

At $t\geq 0$, agent $A_i$ estimates dynamical parameters $\{[\hat{\b{A}}_{j,t},\hat{\b{B}}_{j,t} ];\, j \in Obs_i(N)\}$ of the agents in its observation set, $Obs_i(N)$ . The admissible control set of agent $A_i$ is $\U_{a,i}^N$, consisting of observations of the trajectories and control inputs of all the agents in the set $Obs_i(N)$. Based upon this observation set, agent $A_i$ obtains estimates $\{[\hat{\b{A}}_{j,t},\hat{\b{B}}_{j,t}];\, j\in Obs_i(N)\}$ solving the RWLS equations using $\{dx_{j,t};\, j\in Obs_i\}$ as the measurement variable with the regression vector $\{[x_{j,t}\t,u_{j,t}\t];\, j\in Obs_i\}$.%

\paragraph{Population Cost Function Parameter Identification (PCPI)}

The solution to the RWLS equations with the inputs described above generates the estimates $\{ [\hat{\b{A}}_{j,t},\hat{\b{B}}_{j,t} ];\, j \in Obs_i(N)\}$. The objective at this point for each agent is to obtain the estimates $\{ \hat{\b{Q}}_{j,t};\, j \in Obs_i(N)\}$. The RWLS equations are then solved employing the observed control inputs $\{u_{j}(t);j \in Obs_i(N)\}$ such that agent $A_i$ calculates $\{-(\hat{\b{B}}_{j,t}\t)^{-1}\b{R}u_{j}(t);\, j\in Obs_i(N) \}$ and sets as the measurement vector. Note that one needs the following additional assumption.

\begin{ohypot}\label{ass:invB}
$\b{B}_\theta$ is invertible (and hence, necessarily, $[\b{A}_\theta,\b{B}_\theta]$ is controllable) for all $\th\in\Th$.
\IEEEQED
\end{ohypot}

This rather restrictive assumption is only needed for the cost function parameter identification; therefore, PCPI will be given as an optional procedure in the MF-SAC Law. The observed control action is in the form \eqref{eqn:MF_control}; therefore arranging the variables in a certain way to be specified later, agent $A_i$ obtains the estimates $\{\hat{\b{\Pi}}_{j,t},\hat{s}_{j}(t); \, j\in Obs_i(N)\}$. Solving the algebraic Riccati equation for $\hat{\b{Q}}_{j,t}$ agent $A_i$ obtains its estimates $\{\hat{\b{Q}}_{j,t},\, j\in Obs_i(N)\}$. The symmetry of $\{ \b{Q}_{j,t};\, j\in Obs_i(N)\}$ is guaranteed. To ensure the positive definiteness of the obtained estimates $\{\hat{\b{Q}}_{j,t},\, j\in Obs_i(N)\}$, $[\b{A},\b{B}]$ controllability, $[\b{Q}^{1/2},\b{A}]$ observability, and that the requirement in \ass{ass:compact} holds, the set $\{\hat{\b{A}}_{j,t},\hat{\b{B}}_{j,t},\hat{\b{Q}}_{j,t},\, j\in Obs_i(N)\}$ is projected onto $\Th$.

\paragraph{Distribution Parameter Identification (DPI)}

Once the projected estimates $\hth_{i,t}^{[1:N_0]}\triangleq \{\hat{\b{A}}_{j,t},\hat{\b{B}}_{j,t},\hat{\b{Q}}_{j,t},\, j\in Obs_i(N)\},$ $N_0=\lvert Obs_i(N) \rvert$, are obtained, agent $A_i$ forms the scaled log-likelihood-type function
\begin{equation*}
L(\hth_{i,t}^{[1:N_0]};\z) \triangleq - \frac{1}{N_0}\log\left ( \prod_{j\in Obs_i(N)} f_\z(\hth_{j,t}) \right),
\end{equation*}
calculates $\hzn_{i,t}$, the estimate of the distribution parameter, solving $\argmin_{\z \in P}L(\hth_{i,t}^{[1:N_0]};\z)$. Note that $P$ is known to all agents in the system.

Overall using the identification procedures explained above, agent $A_i$ obtains estimates $[\hat{\b{A}}_{i,t},\hat{\b{B}}_{i,t}]$ and $\hzn_{i,t}$ and forms the self estimated dynamical parameter vector $\hth_{i,t}\t = [\hat{\b{A}}_{i,t},\hat{\b{B}}_{i,t},\b{Q}_i]$. 

\subsubsection{Certainty Equivalence Adaptive Control}

At time $t$, employing $\hzn_{i,t}$ agent $A_i$ solves the MF Equation System \eqref{eqn:ncee} to obtain $x^*(\tau,\hzn_{i,t}),\, t\leq \tau < \infty$. Then using $\hth_{i,t}$ agent $A_i$ solves the Riccati equation \eqref{eqn:DE_Pi}, obtains $\hat{\b{\Pi}}_{i,t} \triangleq \b{\Pi}(\hth_{i,t})$ and solves the mass offset differential equation \eqref{eqn:DE_s} to obtain $\hat s_i(t) \triangleq s(t;\hth_{i,t},\hzn_{i,t})$. The certainty equivalence adaptive control for the admissible control set $\U_{a,i}^N$ is then given by $\huio(t) \triangleq \uio(t;\hth_{i,t},\hzn_{i,t}) = -\b{R}^{-1}\hat{\b{B}}_{i,t}\t(\hat{\b{\Pi}}_{i,t} x_i(t) + \hat s_i(t)),\, t\geq 0$.%, where
%\begin{equation}\label{eqn:opt_control}
%\huio(t)=-\b{R}^{-1}\hat{\b{B}}_{i,t}\t(\hat{\b{\Pi}}_{i,t} x_i(t) + \hat s_i(t)), \quad t\geq 0.
%\end{equation}

% first submission
% We observe that the control law is based on estimates of local parameters obtained from the agent's own trajectory and the population dynamical and cost function distribution parameter estimate obtained from observation on a subset of all agents in the system. 
% first submission
To obtain the main MF-SAC result stated in Theorem \ref{thm:CESACMF}, we first establish the strong consistency for the family of estimates $\{ \hth_{i,t};\, t\geq 0 ,1\leq i \leq N \}$ and $\{ \hzn_{i,t};\, t\geq 0, 1 \leq i \leq N \}$.

\subsubsection{Control Excitation for Consistent Identification}

In order to generate a consistent sequence of estimates $(\hth_{i,t};\, t\geq 0)$ w.p.1, a diminishing excitation is added to the adaptive control in \eqref{eqn:opt_tracking} to give
\begin{equation}\label{eqn:control_with_dit} \huio(t)=-\b{R}^{-1}\hat{\b{B}}_i\t(\hat{\b{\Pi}}_i
x_i(t) + \hat s_i(t)) + \xi_k\left [\e_i(t)-\e_i(k) \right ],\quad t\in(k,k+1],\quad k\in \N,\quad 1\leq i \leq N,
\end{equation}
% or
% \begin{multline} 
% d\hat u_i(t) = -\b{R}^{-1}\hat{\b{B}}_i\t(\hat{\b{\Pi}}_i dx_i(t) + d\hat s_i(t)) + \xi_k d\e_i(t),\\ t\in(k,k+1], \quad k\in \N,\quad 1 \leq i \leq N,
% \end{multline}
where $\huio(0) = 0,\, ( \xi_k^2=\log k/\sqrt{k};\, k\in \Z_1)$, and the process $(\e(t),t\geq 0)$ is an $\R^m$-valued standard Wiener process that is independent of $(w_i(t);\,t\geq 0)$. The sequence of random processes $( \e(t)-\e(k);\,t\in(k,k+1],\,k\in \N)$ is assumed to be mutually independent and all members of the set have the same probability law on $(0,1]$. Since the sequence $(\xi_k;k \in \N)$ converges to zero at a suitable rate, it will be established following \cite{1999DGD_TAC} that the diminishing control excitation $ ( \xi_k [\e(t)-\e(k)];\,t\in[0,1),\,k\in \N )$ provides sufficient excitation for almost sure consistent identification and decreases sufficiently rapidly enough not to affect the limiting performance of the system with respect to $\hth_{i,t}=\th_i^o,t\geq 0$, i.e.\ the non-adaptive case. In other words, the asymptotic performance achieved is equal to the one obtained in the non-adaptive case almost surely. The diminishing control excitation \eqref{eqn:control_with_dit} was introduced in \cite{1986CG_IJC,1991CG}, and it was shown in \cite{1999DGD_TAC} to generate strongly consistent parameter estimates via RWLS for dynamical parameters of the system \eqref{eqn:dynamics} under certainty equivalence adaptive control.%

%===============================================================================
\subsection{The MF Stochastic Adaptive Control (SAC) Law}\label{ssec:ncesacl}

We observe that the control law \eqref{eqn:control_with_dit} has three terms computed from the local state information, the self dynamical parameter estimates and the population distribution parameter estimate. It can be written for each agent $A_i,\, 1\leq i \leq N$, in the form of $\uio(t;\hth_{i,t},\hzn_{i,t})=u_i^{loc}(t;\hth_{i,t})+u_i^{pop}(t;\hth_{i,t},\hzn_{i,t})+u_i^{dit}(t)$, $t\geq 0$, where $u_i^{loc}(\cdot)$ is the LQG feedback for the system of agent $A_i$ based on \emph{local} information; $u_i^{pop}(\cdot)$ is the mass offset term based on local information and \emph{population} information received from the observed set; and $u_i^{dit}(\cdot)$ is the locally generated \emph{dither} input. In this section we present the MF-SAC Law which generates the feedback control law $\huio(t)\triangleq \uio(t;\hth_t,\hzn_t),\, t\geq 0$, that leads to the $\e$-Nash equilibrium. The continuous time MF-SAC Law for agent $A_i,\, 1\leq i \leq N,$ with parameter $\th_i\in \Th,\, 1\leq i \leq N$, is summarized in three major steps in Table \ref{tab:MFSAC}.

%===============================================================================
\begin{table}

{\noindent} 
% \fbox{\begin{minipage}{3.2in}
\fbox{\begin{minipage}{0.965\hsize}

\begin{center}
\textbf{Specification of the MF-SAC Law}
\end{center}

For agent $A_i,\, t\geq 0$:
\begin{enumerate}[(i)]

\item \label{nslspe} Self parameter $\hth_i$ identification:

Solve the RWLS equations \eqref{eqn:nslspe} for the \emph{dynamical parameters}: 

\hfill (subscript $i$ suppressed for clarity)
\begin{equation}\label{eqn:nslspe}
\begin{aligned}
 \upsilon_t\t & =  [\hat{\b{A}}_t,\hat{\b{B}}_t ],\qquad \psi_t\t=[x_t\t,u_t\t],\\
 d\upsilon_t & =  a(t)\Psi_t\psi_t[dx_t\t -\psi_t\t\upsilon_t dt],\\
 d\Psi_t  & =  -a(t) \Psi_t \psi_t \psi_t\t \Psi_t dt,
\end{aligned}
\end{equation}
% and calculate
\be \label{eqn:nslspep}
\text{and calculate }\upsilon_t^{pr} = \argmin_{\psi \in \Th _{|\b{Q}_i}}\lVert \upsilon_t - \psi \rVert, \, \hth_{i,t}^{pr} = [\upsilon_t\t,\b{Q}_i].
\ee

\item Population-parameter $\hzn_i$ identification:
\begin{enumerate}[(a)] 
\item \label{nslppe1} Solve the RWLS equations \eqref{eqn:nslspe} for the \emph{dynamical parameters} $\{\hat{\b{A}}_{j,t},\, \hat{\b{B}}_{j,t},\, j\in Obs_i(N)\}$.

\item \label{nslppe2} Solve the RWLS equations \eqref{eqn:nslppe} for the \emph{cost function parameters} $\{ \hat{\b{\Pi}}_{j,t}, \hat{s}_{j,t},\, j\in Obs_i(N)\}$: (subscript $j$ suppressed for clarity)
\begin{equation}\label{eqn:nslppe}
\begin{aligned}
 \upsilon_t\t  &= [\hat{\b{\Pi}}_t,\hat{s}(t) ],\qquad \psi_t\t=[x_t\t,1],\\
 d\upsilon_t  &=  a(t)\Psi_t\psi_t[(-(\hat{\b{B}}_t\t)^{-1}\b{R}u_t)\t -\psi_t\t\upsilon_t ],\\
 d\Psi_t  &=  -a(t) \Psi_t \psi_t \psi_t\t \Psi_t dt,
\end{aligned}
\end{equation}
solve the algebraic Riccati Equation \eqref{eqn:nslq} for $\{\hat{\b{Q}}_{j,t},\, j \in Obs_i(N)\}$,
\begin{equation}\label{eqn:nslq}
\hat{\b{Q}}_{j,t} = -\hat{\b{A}}_{j,t}\t \hat{\b{\Pi}}_{j,t}\t - \hat{\b{\Pi}}_{j,t} \hat{\b{A}}_{j,t} +
\hat{\b{\Pi}}_{j,t}\t \hat{\b{B}}_{j,t} {\b{R}}^{-1} \hat{\b{B}}_{j,t}\t \hat{\b{\Pi}}_{j,t},
\end{equation}
% set $\hth_{j,t}\t = [\hat{\b{A}}_{j,t},\hat{\b{B}}_{j,t},\hat{\b{Q}}_{j,t}]$, and calculate 
set $\hth_{j,t}\t = [\hat{\b{A}}_{j,t},\hat{\b{B}}_{j,t},\hat{\b{Q}}_{j,t}]$, and calculate 
\begin{equation}\label{eqn:nslqp} 
\hth_{j,t}^{pr} = \argmin_{\psi \in \Th} \lVert \hth_{j,t}\t - \psi \rVert.
\end{equation}

\item \label{nslmle} Solve the MLE equation \eqref{eqn:nslmle} at $\hth_{i,t}^{[1:N_0]} = [\hat{\b{A}}_{j,t}^{pr},\hat{\b{B}}_{j,t}^{pr},\hat{\b{Q}}_{j,t}^{pr}],\, j \in Obs_i(N)$, to estimate $\z^0$ via:
\begin{equation}\label{eqn:nslmle}
\begin{aligned}
L(\hth_{i,t}^{[1:N_0]};\z) = -\frac{1}{N_0}\log\left ( \prod_{j\in Obs_i(N)} f_\z(\hth_{j,t}) \right),\\
\hzn_{i,t} = \argmin_{\z \in P} L(\hth_{i,t}^{[1:N_0]} ;\z),\quad N_0 = \lvert Obs_i(N) \rvert,
\end{aligned}
\end{equation}
and solve the set of MF Equations \eqref{eqn:ncee} for all $\th\in\Th$ generating $x^*\left (\tau,\hzn_{i,t}\right ),\, t\leq \tau < \infty$. 
\end{enumerate}

\item \label{nslcle} Solve the MF Control Law Equation at $\hth_{i,t}^{pr}$ and $\hzn_{i,t}$:

\begin{enumerate}[(a)]

\item \label{nslr} $\hat{\b{\Pi}}_{i,t}$: Solve the Riccati Equation \eqref{eqn:nslr} at $\hth_{i,t}^{pr}$:
\begin{equation}\label{eqn:nslr}
\hat{\b{A}}_{i,t}\t \hat{\b{\Pi}}_{i,t} + \hat{\b{\Pi}}_{i,t} \hat{\b{A}}_{i,t} -\hat{\b{\Pi}}_{i,t} \hat{\b{B}}_{i,t} \b{R}^{-1}\hat{\b{B}}_{i,t}\t \hat{\b{\Pi}}_{i,t} + \b{Q}_i = 0.
\end{equation}

\item \label{nsls} $\hat{s}_i(t) \triangleq s(t;\hth_{i,t}^{pr},\hzn_{i,t})$:  Solve the mass offset differential equation \eqref{eqn:nsls} at $\hth_{i,t}^{pr}$ and $\hzn_{i,t}$:
\begin{equation} \label{eqn:nsls} 
-\frac{d\hat{s}_i(\tau)}{d\tau} = (\hat{\b{A}}_{i,t}\t - \hat{\b{\Pi}}_{i,t} \hat{\b{B}}_{i,t} \b{R}^{-1}\hat{\b{B}}_{i,t}\t)  \hat{s}_i(\tau)  - \b{Q}_i x^*(\tau,\hzn_{i,t} ),\quad t\leq \tau <\infty.
\end{equation}

\item \label{nslu} Obtain the Certainty Equivalence Adaptive Control at $\hth_{i,t}^{pr}$ and $\hzn_{i,t}$:
\begin{equation}\label{eqn:nslu}
\huio(t) = -\b{R}^{-1} \hat{\b{B}}_{i,t}\t \Big ( \hat{\b{\Pi}}_{i,t}  x_i (t) + \hat{s}_i(t) \Big )  + \xi_k\left [\e_i(t)-\e_i(k) \right ],\quad t\in (k,k+1],\quad k\in \N.
\end{equation}

\end{enumerate}
\end{enumerate}
\end{minipage}}

\caption{MF-SAC LAW}
\label{tab:MFSAC}
\end{table}

%===============================================================================

The function $a(t),t\geq 0$, in \eqref{eqn:nslspe} is in the form of $a(t)=1/f(r(t))$, where $r(t)=\lVert \Psi_0^{-1} \rVert + \int_0^t \lvert \psi(s) \rvert^2 ds$, and $f\in \{ f: \R_+ \conv \R_+, \, f \text{ is slowly increasing and} \, \int_c^{\infty} 1/(xf(x)) dx < \infty;\, c\geq 0 \}$. The function $f(.)$ is slowly increasing if it is increasing and satisfies $f(.)\geq 1$ and $f(x^2)=O(f(x))$ \cite{1999DGD_TAC}. 

Note that a positive definite solution to the Riccati equation \eqref{eqn:nslr} exists as the projected estimate is in the set of controllable and observable dynamical parameters: $\hth_t^{pr}\in \Th \subset \tilde{\Th} \subset \R^{n(n+m+(n+1)/2)}$.%

\subsection{Asymptotic Properties of the MF-SAC Law}

% first submission
% Note that each agent starts with no prior information on its self parameter, nor on the population dynamical and cost function distribution parameter. 
% first submission
A key feature of the work in this paper is that the state aggregation integration in \eqref{eqn:ncee} is performed by use of the estimated distribution $F_{\hzn}(\cdot)$ in place of the true distribution $F_{\z^0}(\cdot)$ (see \eqref{eqn:nceea} below). Then the central results of this paper are the following: under the MF-SAC Law, asymptotically as the population tends to infinity, the competitive best response actions of the adaptive agents with no prior information on self dynamical parameters and no prior statistical information on dynamical and cost function parameters of the mass give rise to a unique Nash equilibrium. Moreover, the resulting cost for each agent from the MF-SAC Law is asymptotically almost surely equal to the cost resulting from the non-adaptive MF Stochastic Control Law. 
% The theorem that gives this result is presented in Section \ref{sec:mai_the}. 

%%%%%%%%%%%%%%%%%%%%%%%%%%%%%%%%%%%%%%%%%%%%%%%%
%%%%%%%%%%%%%%%%%%%%%%%%%%%%%%%%%%%%%%%%%%%%%%%%
\begin{thm}\label{thm:CESACMF} \emph{MF-SAC Theorem} 

Let \ass{ass:ind}-\ass{ass:weakdist}, \ass{ass:pdf}, \ass{ass:MLE_ident} hold. Then, assume each agent $A_i,\, 1\leq i \leq N$, is such that it:
\begin{enumerate}[(i)]
\item observes a random subset $Obs_i(N)$ of the total population $N$ such that $\lvert Obs_i(N) \rvert \conv \infty,\, \lvert Obs_i(N) \rvert/N \rightarrow 0,$ as $ N\rightarrow \infty;$
 
\item estimates its own parameter $\hth_{i,t}$ via the RWLS \eqref{eqn:nslspe}; 

\item estimates the population dynamical and cost function distribution parameter $\hzn_{i,t}$ via MLE \eqref{eqn:nslmle}; and 

\item computes $\uio(t;\hth_{i,t},\hzn_{i,t})$ via the extended MF equations plus dither.% where $ \bar x_\theta  = \int_{\Th}\bar x_\tau$ $d F_{\hzn_t}(\th)$. 
\end{enumerate}

Then,

\begin{enumerate}[(a)]
\item  $\hth_{i,t} \conv \th_i^0$ w.p.1 as $t\conv \infty,\, 1 \leq i \leq N$ (strong consistency);

\item $\hzn_{i,t}\conv \z^0$ w.p.1 as $t\conv\infty$, and $N\conv\infty,\, 1 \leq i \leq N$. 

The MF-SAC Law generates a set of controls $\hat{\U}_{MF}^N = \{ \huio;\, 1\leq i \leq N \},\, 1\leq N < \infty,$ such that:

\item all agent system trajectories $x_i,\, 1\leq i \leq N,$ are $LRA-L^2$ stable w.p.1; 

\item \textit{$\e-$Nash Property:} $\{ \hat{\U}_{MF}^N;\, 1\leq N < \infty\}$ yields an $\e$-Nash Equilibrium for all $\e$, i.e., for all $ \e > 0$, there exists $N(\e)$ such that for all $N\geq N(\e)$
\begin{equation*}
J_i^N(\huio,\humio)-\e \leq \inf_{u_i\in \U_g^N} J_i^N (u_i,\humio) \leq  J_i^N(\huio,\humio) \quad \text{w.p.1,} \quad 1\leq i \leq N; 
\end{equation*}

\item \textit{Equal Adaptive and Non-adaptive $(\th^0,\z^0)$ MF Equilibrium Performance:}
\begin{equation*}
\limN J_i^N ( \huio,\humio) = \limN J_i^N(u_i^0,u_{-i}^0 )\quad \text{w.p.1},\quad 1\leq i \leq N;
\end{equation*}
\item \textit{Adaptive Control Performance Equals Complete Information Performance:}
\begin{equation*}
\limN J_i^N(\huio,\humio) = \limN \inf_{u_i \in \U_g^N}J_i^N (u_i,\humio) \quad \text{w.p.1}, \quad 1\leq i \leq N.
\end{equation*}

\end{enumerate}\IEEEQED
\end{thm}

% The proof consists of the unification of the principal theorems proved earlier, the Propositions \ref{prop:eq_costs} and \ref{prop:eq_costs_g} above and the key lemmas established in Appendix \ref{sec:app_MT}.

%v3
The proof consists of the unification of the principal Theorems \ref{thm:str_con_est}, \ref{thm:scmle}, \ref{thm:L2_stable} and the Propositions \ref{prop:eq_costs} and \ref{prop:eq_costs_g} that are presented in the remaining sections. The outline of the proof is given in Appendix \ref{thm:CESACMF_pro}.

The technical plan of the paper is presented in three layers. The main theorem of the paper is Theorem \ref{thm:CESACMF}. In the first layer, Propositions \ref{prop:eq_costs}, \ref{prop:eq_costs_g} and Theorems \ref{thm:str_con_est}, \ref{thm:scmle} and \ref{thm:L2_stable} support Theorem \ref{thm:CESACMF}. In the second layer, Lemmas \ref{lem:stability}, \ref{lem:con_x_m}, \ref{lem:lim_Ji_opt_N}, \ref{lem:lim_Ji_N}, \ref{lem:mass_conv}, Theorem \ref{thm:pop_L2_conv} and Proposition \ref{prop:pop_u_conv} support Proposition \ref{prop:eq_costs} whereas Lemma \ref{lem:uni_con} supports Theorem \ref{thm:str_con_est}. In the third layer, Lemmas \ref{lem:state_transition}, \ref{lem:state_transition_conv}, \ref{lem:kronecker} and Proposition \ref{prop:s_conv} support Theorem \ref{thm:pop_L2_conv}.

%===============================================================================
\section{Convergence Properties of the MF-SAC Parameter Estimates} \label{sec:cppe}

% Finally, the distribution parameter for the dynamical and cost function parameters is estimated by the MLE equation.  

We show that for self dynamical parameter identification, the RWLS equations for \emph{dynamical parameters} \eqref{eqn:nslspe} with the projection method \eqref{eqn:nslspep} provide strongly consistent, uniformly controllable and observable estimates. The \emph{population dynamical and cost function distribution parameter} identification is handled in three steps. First, each agent obtains the dynamical parameter estimates for the agents in its observation set solving the RWLS equations \eqref{eqn:nslspe}. It is shown that the RWLS equations \eqref{eqn:nslspe} with the projection method \eqref{eqn:nslspep} applied on the observed agents' controlled trajectories also provide strongly consistent, uniformly controllable and observable estimates. Secondly, another set of RWLS equations \eqref{eqn:nslppe} are solved using the previously obtained dynamical parameter estimates as inputs; and finally cost function parameter estimates are obtained for the agents in the observation set \eqref{eqn:nslq}. We show that the estimates obtained are positive definite and uniformly bounded by use of a projection method \eqref{eqn:nslqp}. Finally, we show that the MLE scheme \eqref{eqn:nslmle} employed using these estimates provides strongly consistent population distribution parameter estimates.

% for population dynamical and cost function distribution parameter each agent first obtains the dynamical parameter estimates for the agents in its observation set. Then solves the RWLS equations \eqref{eqn:nslppe}, which with the projection method provide positive definite cost function parameters. The MLE equation \eqref{eqn:nslmle} employed on these estimates provides strongly consistent population distribution parameter estimates.

%===============================================================================
\subsection{Asymptotic Convergence of the Dynamical Parameter Estimates} \label{ssec:acrwlspe}

The RWLS algorithm is self-convergent \cite{1996Gu_TAC}, i.e., it converges to a certain random vector almost surely irrespective of the control law design, but there is no guarantee that the estimated dynamical parameters will be controllable and observable, or the cost function estimates will be positive definite. To ensure that the sequence of estimated dynamical parameters are controllable, observable, uniformly bounded and the sequence of estimated cost parameters are positive definite and uniformly bounded we use the \emph{projection method} \cite{1992Ca_SCL}.

For self dynamical parameter identification, the self dynamical parameter estimates with the cost function parameter $\b{Q}_i\in \R^{n(n+1)/2}, \, \hth_{i,t}\t = [\hat{\b{A}}_{i,t},\hat{\b{B}}_{i,t},\b{Q}_i] \in \R^{n\left (n+m+(n+1)/2 \right)},\, t\geq 0$, ($\b{Q}_i$ known by agent $A_i$) is projected (denoted by $\hth_{i,t}^{pr}$ in \eqref{eqn:nslspep}) onto the compact set $\Th_{|Q_i}\subset \tilde{\b{\Th}}_{|Q_i}$, where for the given $\b{Q}_i$, $[\b{A}_\theta,\b{B}_\theta]$ is controllable and $[\b{Q}_i^{1/2},\b{A}_\theta]$ is observable.

For the distribution parameter identification, the population dynamical parameter estimates together with the cost function parameter estimates are projected onto the compact subset $\Th$ of the set of controllable and observable dynamical parameters $\tilde{\Th}$ where, in addition, $\b{Q}_{\tilde{\th}},\, \tilde\th\in\tilde\Th$, is positive definite (for which the control law generated by \eqref{eqn:nslr} necessarily exists and is asymptotically stabilizing). 

\begin{lem}\label{lem:uni_con}
Let $\Th$ be a compact set such that $\th_i^0 \in \Th\subset \tilde{\Th} \subset \R^{n\left(n+m+(n+1)/2\right)},\, 1\leq i \leq N$. Set $\hth_{i,t}\t = [\hat{\b{A}}_{i,t},\hat{\b{B}}_{i,t},\hat{\b{Q}}_{i,t}],\, t\geq 0$. Let $[\hat{\b{A}}_{i,t},\hat{\b{B}}_{i,t}]\t$ be the estimate of $[\b{A}_i^0,\b{B}_i^0]$ obtained by the RWLS equations \eqref{eqn:nslspe}, and let $\hat{\b{Q}}_{i,t}$ be the estimate of $\b{Q}_i^0 $ obtained by the RWLS equations \eqref{eqn:nslppe} and \eqref{eqn:nslq}. Assume $\hth_{i,t}\conv\th_i^0$ w.p.1 as $t\conv\infty,\, 1\leq i \leq N$. Then, $\hth_{i,t}^{pr}  \teq [\hat{\b{A}}_{i,t}^{pr},\hat{\b{B}}_{i,t}^{pr},\hat{\b{Q}}_{i,t}^{pr}] \triangleq \argmin_{\psi\in \Th} \lVert \hth_{i,t} - \psi\rVert$ (together with a co-ordinate ordering measurable tie breaking rule), satisfies $\hth_{i,t}^{pr} \in \Th$ w.p.1 for all $t\geq 0$, and $\hth_{i,t}^{pr}\conv\th_i^0$ w.p.1 as $t\conv\infty$. In the SDPI case the corresponding result is achieved by setting $\hat{\b{Q}}_{i,t} = \b{Q}_i^0$ for all $t\geq 0$.\IEEEQED
\end{lem}
The Lemma is proved in Appendix \ref{lem:uni_con_pro}.
%%%%%%%%%%%%%%%%%%%%%%%%%%%%
% \begin{lem}\label{lem:uni_con}\textit{Self Dynamical Parameters Identification Case} 

% Let $\Th$ be a compact set such that $\th^0 \in \Th\subset \tilde{\Th} \subset \R^{n(n+m+(n+1)/2)}.$ Let $\hth_t$ be the estimate of $\th^0 \in \Th$ obtained by the RWLS equations \eqref{eqn:nslspe}. Then, $\hth_t^{pr} \triangleq \argmin_{\psi\in \Th} \lVert \hth_t - \psi\rVert$ (together with a co-ordinate ordering measurable tie breaking rule), satisfies $\hth_t^{pr} \conv \th^0$ w.p.1 as $t\conv \infty$.\IEEEQED
% \end{lem}
%The result is proved in Appendix \ref{lem:uni_con_pro}.
%===============================================================================
% \subsection{Asymptotic Convergence of the RWLS Cost Function Parameter Identification} \label{ssec:acrwlspe}
%%%%%%%%%%%%%%%%%%%%%%%%%%%%
% \begin{lem}\label{lem:uni_con_q}\textit{Population Cost Function Parameter Identification Case}

% Let $\Th$ be a compact set such that $\b{Q}^0 \in \Th\subset \tilde{\Th} \subset \R^{n(n+m+(n+1)/2)}.$ Let $\hat{\b{Q}}_t$ be the estimate of $\b{Q}^0 \in \Th$ obtained by the RWLS equations \eqref{eqn:nslppe}. Then, $\hat{\b{Q}}_t^{pr} \triangleq \argmin_{\psi\in \Th} \lVert \hat{\b{Q}}_t - \psi\rVert$ (together with a co-ordinate ordering measurable tie breaking rule), satisfies $\hat{\b{Q}}_t^{pr} \conv \b{Q}^0$ w.p.1 as $t\conv \infty$.\IEEEQED
% \end{lem}

Now, given the projection method lemma, we show that the RWLS equations for dynamical parameters \eqref{eqn:nslspe} and the RWLS equations for cost function parameters \eqref{eqn:nslppe} generate strongly consistent estimates.
%%%%%%%%%%%%%%%%%%%%%%%%%%%%%%
\begin{thm}\label{thm:str_con_est}
Let hypotheses \ass{ass:ind}-\ass{ass:compact} hold, $x^* \in \b{C}_b[0,\infty)$, and let $( [\hat{\b{A}}_{i,t},\hat{\b{B}}_{i,t}] ;\, t\geq 0),\,1\leq i \leq N$, be the process of estimates obtained by the RWLS equations \eqref{eqn:nslspe}, and $(\hat{\b{Q}}_{i,t};\, t\geq 0)$ be the process of estimates obtained by \eqref{eqn:nslq} along the controlled trajectory $((x_{i,t},\hat{u}_{i,t}^0);\, t\geq 0)$, generated by the control $(\hat{u}_{i,t}^0;\,t\geq 0)$ according to the MF-SAC Law \eqref{eqn:nslu}. Furthermore, let $( \hth_{i,t}^{pr} =[\hat{\b{A}}_{i,t}^{pr},\hat{\b{B}}_{i,t}^{pr},\hat{\b{Q}}_{i,t}^{pr}];\, t\geq 0 )$, be the projected estimates according to Lemma \ref{lem:uni_con}. Then, % given $\hzn_i\in\b{C}_b[0,\infty)$,
\begin{enumerate}[(i)] 
\item the input process given in \eqref{eqn:nslu} is well defined,% and is given at $\hth_t^{pr}$ by \eqref{eqn:nslu},
% ACK-LV
%\begin{equation}\label{eqn:th_control} 
%\uo (t;\hth_t^{pr},\hzn_t) = -\b{R}^{-1}\hat{\b{B}}_{t}\t (\hat{\b{\Pi}}_{t} x_t + s(t;\hth_t^{pr},\hzn_t) ) + \xi_k\left [\e(t)-\e(k) \right ],
%\end{equation}
%ACK-LV
\item $[\hat{\b{A}}_{i,t},\hat{\b{B}}_{i,t}] \conv [\b{A}_i^0,\b{B}_i^0] \quad \text{w.p.1 as } t\conv \infty$, $1\leq i \leq N$,
\item with the optional assumption \oass{ass:invB}, $\hat{\b{Q}}_{i,t} \conv \b{Q}_i^0$ w.p.1 as $t\conv\infty$, $1\leq i \leq N$.
\end{enumerate}\IEEEQED
\end{thm}

The theorem is proved in Appendix \ref{thm:str_con_est_pro} using the methodology of \cite{1999DGD_TAC}, which establishes the convergence of the RWLS estimates \eqref{eqn:nslspe} with diminishing excitation in the controls \eqref{eqn:nslu}. The required uniform controllability and observability of the estimates is a consequence of Lemma \ref{lem:uni_con} since $\hth_t^{pr} \in \Th,\, t\geq 0$.

%===============================================================================
\subsection{Asymptotic Convergence of the Population Distribution Parameter Estimates} \label{ssec:acppes}

The MF-SAC Law specifies that the distribution parameter identification is such that each agent $A_i,\, 1\leq i \leq N,$ observes the control and state trajectories of a random subset of agents $Obs_i(N),\, 1\leq i \leq N$, and at each time iteration applies \eqref{eqn:nslspe} to obtain the dynamical parameter estimates of each agent in its set. The MLE scheme \eqref{eqn:nslmle} is then applied to these estimated parameters of the agents $Obs_i(N),\,1\leq i \leq N$, for $t\geq 0$, to obtain an estimate of the distribution parameter. To obtain the strong consistency of the distribution parameter estimates we adopt the hypotheses \ass{ass:pdf} and \ass{ass:MLE_ident} below.

\begin{hypot}\label{ass:pdf}
There exists a bounded continuous (on $\Th\times \tilde P$) family of densities $f_\z \triangleq \{ f_\z(\th);\th \in \Th, \z\in \tilde P\}$ for the family of dynamical and cost function parameter distributions $\{ F_\z(.);\, \z \in P\}$. Further, the distribution function $f_\z(\th)$ is bounded away from 0 uniformly over $\Th\times \tilde P$, i.e., $f_\z(\th)\geq \delta$ for some $\delta>0$ for all $\th \in \Th$ and $\z \in \tilde P$. Moreover, for each $j,\, 1\leq j \leq p$, $(\p f_\z/\p \z_j)(\th)$ exists for all $\z \in \tilde P$, and is uniformly bounded on $\Th\times \tilde P$, except possibly on a Lebesgue null set independent of $\z\in \tilde P$.
%, where $P$ is a compact set such that $P\subset \tilde P \subset \R^p$, and $\tilde P$ is an open set
\IEEEQED
\end{hypot}

For \eqref{eqn:nslmle}, let $f(\th^{[1:N_0]};\, \z)$ be the likelihood function of $f_\z$ at $\th^{[1:N_0]}\triangleq \{\b{A}_j,\b{B}_j,\b{Q}_j,\, j\in Obs_i(N),\, N_0 = \lvert Obs_i(N) \rvert \},$ and let $L(\th^{[1:N_0]};\, \z)$ be the continuously differentiable monotonically decreasing function of $f(\th^{[1:N_0]} ;\, \z)$ given by the scaled log-likelihood function $-(1/N)$ $\log f(\th^{[1:N_0]};\, \z).$

\begin{hypot}\label{ass:MLE_ident}
$\{f_\z(\cdot);\,\z\in P\}$ satisfies: \[ \E_{\z^0} [\log f_{\z}(\th)] = \E_{\z^0} [\log f_{\z'}(\th)]    \iff  \z = \z',\] for all $\z,\z',\z^0 \in P$, where $\z^0$ is the true parameter.
\IEEEQED
\end{hypot}

\begin{thm}\label{thm:scmle}
Let \ass{ass:ind}-\ass{ass:compact}, \ass{ass:pdf}, \ass{ass:MLE_ident} hold; let $\lvert Obs_i(N) \rvert\conv\infty$ and $\lvert Obs_i(N) \rvert/N\conv 0$ as $N\conv\infty$, $1\leq i \leq N$, and let $( \hzn_{i,t}(\hth_{i,t}^{[1:N_0]}),\, t\geq 0),\, N_0 = \lvert Obs_i(N) \rvert$, be the MLE process given by \eqref{eqn:nslmle} along the controlled trajectories of the observed set of agents $( (x_{j,t}, \hat{u}_{j,t}^0 ); \, t\geq 0, \, j \in Obs_i(N) )$ generated by the controls $(\huo_{j,t};\, t\geq 0,j\in Obs_i(N))$ \eqref{eqn:nslu}. Then,
%\begin{enumerate}[(i)] 
$\hzn_{i,t}$ is strongly consistent at $\z^0$, that is, $\limN\lim_{t\conv\infty} \hzn_{i,t}( \hth_{i,t}^{[1:N_0]} ) = \z^0$ w.p.1, $1\leq i \leq N$.
% \item $\limN \lim_{t\conv\infty} f_{\hzn_{i,t}}(\cdot) = f_{\z^0}(\cdot)$ uniformly on $P$ \ack{???}, w.p.1, $1\leq i \leq N$. \ack{uniform convergence proof for the distribution is missing}
%\end{enumerate}
\IEEEQED
\end{thm}
The proof is given in Appendix \ref{thm:scmle_pro}.
%%%%%%%%%%%%%%%%%%%%%%%%%%%%%%%%%%%%%%%%%%%
% Dynamical Properties of Parameter Convergent Systems
%%%%%%%%%%%%%%%%%%%%%%%%%%%%%%%%%%%%%%%%%%%
% \section{Dynamical Properties of Parameter Convergent Systems}\label{sec:dyn_pro_par_con_sys}

%===============================================================================
\section{The Principal Asymptotic Results}\label{sec:pr}

%===============================================================================
\subsection{Asymptotic Behaviour of the MF Equations}\label{ssec:abncee}

% In this section we present two results under the hypotheses of convergent individual dynamical parameters and the population dynamical and cost function distribution parameters. The convergence results for the mass offset function and the control function are shown below. 

The MF Equations \eqref{eqn:ncee} that permit the calculation of the mass tracking signal $x^*(\tau,z),\, t_0\leq \tau < \infty$, are dependent on the population distribution parameter $\z$. Correspondingly, the MF Equations of the MF-SAC Law on $[t,\infty),\, t\geq 0,$ with the strongly consistent distribution parameter estimate $\hzn_{i,t},\, t\geq 0, \, 1\leq i \leq N$, are given below.
\begin{defn}\label{MFSAC_Equations}
\emph{MF-SAC Equation System on $[t,\infty)$}:
\begin{equation}\label{eqn:nceea}
\begin{aligned}
-\frac{ds_\theta}{d\tau} & = (\b{A}_{\th}\t - \b{\Pi}_{\th} \b{B}_{\th}\b{R}^{-1}\b{B}_{\th}\t ) s_\theta- \b{Q}_\theta x^*(\tau,\hzn_{i,t}),\\
\frac{d\bar x_\theta}{d\tau} & = (\b{A}_{\th}-\b{B}_{\th}\b{R}^{-1}\b{B}_{\th}\t \b{\Pi}_{\th})\bar{x}_\theta -\b{B}_{\th}\b{R}^{-1}\b{B}_{\th}\t s_\theta,\\
\bar x(\tau,\hzn_{i,t}) & = \int_{\Th}\bar x_\theta d F_{\hzn_{i,t}}(\th),\\
x^*(\tau,\hzn_{i,t}) & = m(\bar{x} (\tau,\hzn_{i,t})+\eta),\qquad t \leq \tau < \infty.
\end{aligned}
\end{equation}\IEEEQED
\end{defn}

%%%%%%%%%%%%%%%%%%%%%%%%%%%%%%%%%%%%%%%%%%%%
% Proposition
%%%%%%%%%%%%%%%%%%%%%%%%%%%%%%%%%%%%%%%%%%%%
% S convergence
\begin{prop}\label{prop:s_conv} 
For the system \eqref{eqn:dynamics} let \ass{ass:ind}-\ass{ass:cost_coup}, \ass{ass:pdf}, \ass{ass:MLE_ident} hold. For agent $A_i,\, 1\leq i \leq N$, let: (i) $\hth_{i,t}^{pr}$ be the solution to \eqref{eqn:nslspep}, $\hzn_{i,t}$ be the solution to \eqref{eqn:nslmle} in the MF-SAC Law; (ii) $x^*(\tau,\hzn_{i,t}),\, t\leq \tau < \infty$, be the solution to the MF-SAC Equation System \eqref{eqn:nceea}; $x^*(\tau,\z^0),\, t\leq \tau < \infty$, be the solution to the MF Equation System \eqref{eqn:ncee}; (iii) $s(t;\hth_{i,t}^{pr},\hzn_{i,t})$ be the solution to \eqref{eqn:nsls} in the MF-SAC Law; and $s(t;\th_i^0,\z^0)$ be the solution to the mass offset function differential equation \eqref{eqn:DE_s}. Then,
%$\hth_t \conv \th^0$ w.p.1 as $t\conv\infty$; $\hzn_t\conv \z^0$ w.p.1 as $t\conv \infty,\, N\conv \infty$.
\begin{enumerate}[(i)]
\item $\limN\lim_{t\conv\infty} x^*(\tau, \hzn_{i,t}) = x^* (\tau, \z^0)$ w.p.1, $t\leq \tau <\infty,\,1\leq i \leq N$, 
\item $\limN\lim_{t\conv\infty} s(t;\hth_{i,t}^{pr},\hzn_{i,t}) = s(t;\th_i^0,\z^0)$ w.p.1, $1\leq i \leq N$,
\item The input process given in \eqref{eqn:nslu} is well defined and is given at $\hth_{i,t}^{pr}$ and $\hzn_{i,t}$ by
\begin{equation*} 
\uio (t;\hth_{i,t}^{pr},\hzn_{i,t}) = -\b{R}^{-1}\hat{\b{B}}_{i,t}\t (\hat{\b{\Pi}}_{i,t} x_{i,t} + s(t;\hth_{i,t}^{pr},\hzn_{i,t}) ) + \xi_k\left [\e_i(t)-\e_i(k) \right ].
\end{equation*}

\end{enumerate}

\IEEEQED

\end{prop}

The result is proved in Appendix \ref{prop:s_conv_pro}.

%===============================================================================
\subsection{Asymptotic Behaviour of System Trajectories}\label{ssec:abst}

We show that under the hypotheses that the \emph{self dynamical parameter estimates} and the \emph{population distribution parameter estimates} converge to their true values, the trajectories of adaptive individual agents are stable in the $L^2-LRA$ sense. Moreover, these trajectories and the corresponding control actions converge to the non-adaptive values obtained with the true parameters. 

Recall that $\U_{MF}^N = \{ u_i^0;1\leq i \leq N\}$ is the set of controls generated by the non-adaptive MF Stochastic Control Law, while $\hat{\U}_{MF}^N = \{ \huio;1\leq i \leq N\}$ is the set of controls generated by the MF-SAC Law.

Using the notation $\hth_i^{0,t}\triangleq (\hth_{i,\tau},\, 0\leq \tau \leq t)$, and $\hz_i^{0,t}(N_0)\triangleq (\hzn_{i,\tau},\, 0\leq \tau \leq t)$, let $\hat{x}_i^0\triangleq  \xio(t;\hth_i^{0,t},\hz_i^{0,t}(N)   )  $ be the state trajectory of agent $A_i,\, 1\leq i \leq N$, under the control law $\uio(t;\hth_{i,t},\hzn_{i,t}) \in \hat{\U}_{MF}^N$, and $x_i^0 \triangleq x_i(t;\th_i^o,\z_i^0)$ be the state trajectory of agent $A_i$ under the control law $u_i^0 \triangleq u_i^0(t;\th_i^0,\z_i^0) \in \U_{MF}^N $, where $\hth_{i,t}$ is the solution to \eqref{eqn:nslspep}, and $\hzn_{i,t}$ is the solution to \eqref{eqn:nslmle}.

%%%%%%%%%%%%%%%%%%%%%%%%%%%%%%%%%%%%%%%%%%%%%%%%%
% Theorem
% Stability of adaptive systems
\begin{thm}\label{thm:L2_stable}
Let \ass{ass:ind}-\ass{ass:cost_coup} hold; then, the process $(\hxio(t);\, t\geq 0),\, 1\leq i \leq N$, is stable in the sense that 
\begin{equation*}
\sup_{N\geq 1}\max_{1\leq i \leq N} \limT\frac{1}{T}\int_0^T \lVert \hxio(t) \rVert^2 dt < \infty \quad \text{w.p.1}.
\end{equation*}
\IEEEQED
\end{thm}
%===============================================================================
%%%%%%%%%%%%%%%%%%%%%%%%%%%%%%%%%%%%%%%%%%%%%%%%%%%%%%%%%%%%%%%%%%%%%%%%%%%%%%%%
% Proof of Theorem thm:L2_stable
%%%%%%%%%%%%%%%%%%%%%%%%%%%%%%%%%%%%%%%%%%%%%%%%%%%%%%%%%%%%%%%%%%%%%%%%%%%%%%%%
% \subsection{Proof of Theorem \ref{thm:L2_stable}}\label{thm:L2_stable_pro}
\begin{IEEEproof}
It has been shown in Theorem \ref{thm:str_con_est} that $\hth_{i,t}\conv\th_i^0$ w.p.1, and in Theorem \ref{thm:scmle} that $\hzn_{i,t}\conv\z^0$ w.p.1 as $t\conv\infty$ and $N\conv\infty,\, 1\leq i \leq N$. Moreover, it has already been shown in Proposition \ref{prop:s_conv} that the tracking signal $x^*(\tau,\hzn_t)\in \b{C}_b[0,\infty)$, and the input process is well defined. All the hypotheses in \cite{1999DGD_TAC} are satisfied, and Theorem 1 in \cite{1999DGD_TAC} proves the claim.
\end{IEEEproof}
%%%%%%%%%%%%%%%%%%%%%%%%%%%%%%%%%%%%%%%%%%%%%%%%%
% Theorem
% MF convergence proof with consistent estimates
\begin{thm}\label{thm:pop_L2_conv}
For the system \eqref{eqn:dynamics}, under \ass{ass:ind}-\ass{ass:cost_coup}, \ass{ass:pdf}, \ass{ass:MLE_ident}
\[
\lsN\lsT \frac{1}{T}\int_0^T\lVert  \hat{x}_i^0 - x_i^0 \rVert^2 dt = 0\, \text{w.p.1},\, 1\leq i \leq N.
\]\IEEEQED
\end{thm}
% let $\hth_{i,t} \conv \th_i^0$ w.p.1 as $t\conv \infty$ and $\hzn_{i,t} \conv \z_i^0$ w.p.1 as $t\conv \infty,\, N\conv\infty,\, 1\leq i \leq N$. 
The result is proved in Appendix \ref{thm:pop_L2_conv_pro}.

%===============================================================================
\subsection{Asymptotic Behaviour of Cost Functions}\label{ssec:abcf}

In the population limit, the asymptotic cost of an agent performing the MF-SAC Law in a system within which all of the agents are adaptive is almost surely equal to the cost of an agent in a system of agents all of which are performing the non-adaptive MF-SC Law. This is shown in Proposition \ref{prop:eq_costs} whose proof is given in Appendix \ref{prop:eq_costs_pro}. Moreover, Proposition \ref{prop:eq_costs_g} shows that in the population limit, the best response of an agent in a population of agents performing the MF-SAC Law is almost surely equal to the best response of an agent in a population of agents performing the non-adaptive MF-SC Law. The proof is given in Appendix \ref{prop:eq_costs_g_pro}. 
%v3
%Finally, the main result of the paper is given in Theorem \ref{thm:CESACMF}, and it is shown that for each agent the self dynamical parameter estimates and the population distribution parameter estimates are strongly consistent; all agent systems are long run average $L^2$ stable; the set of controls yields a (strong) $\e$-Nash equilibrium for all $\e$; and in the population limit for each agent the long run average cost obtained is almost surely equal to the cost obtained in a system of all non-adaptive agents. 
% Then we show in Lemma \ref{prop:eq_costs_g} that the MF-SAC Law gives the asymptotic optimal control solution when all other agents apply the MF-SAC Law, and we present the proof in Section \ref{prop:eq_costs_g_pro}.
%We set 
%\begin{equation*}
%J_i^\infty(u_i,u_{-i}) = \lsT \frac{1}{T} \int_0^T \left \{ \lVert x_i - x^* \rVert_Q^2 + \lVert u_i \rVert_R^2 \right \} \quad \text{w.p.1},\quad 1\leq i \leq N, 
%\end{equation*}
%where $x^*$ is given by \eqref{eqn:nceea} in the SAC $(\huio)$ case, and by \eqref{eqn:ncee} in the non-SAC $(u_i^0)$ case.

%%%%%%%%%%%%%%%%%%%%%%%%%%%%%%%%%%%%%%%%%%%%%%%%%%%%%%%%%%%%%%
% Proposition
%%%%%%%%%%%%%%%%%%%%%%%%%%%%%%%%%%%%%%%%%%%%%%%%%%%%%%%%%%%%%%
\begin{prop}\label{prop:eq_costs}
For the system \eqref{eqn:dynamics}, let \ass{ass:ind}-\ass{ass:cost_coup}, \ass{ass:pdf}, \ass{ass:MLE_ident} hold, let $ u_i^0 \in \U_{MF}^N,\, 1\leq i \leq N,$ be the set of controls generated by the non-adaptive $(\th^0,\z^0)$ MF Stochastic Control Law, and let $\huio \in \hat{\U}_{MF}^N,\, 1\leq i \leq N,$ be the set of controls generated by the MF-SAC Law. Then,
\begin{equation}\label{eqn:eq_costs}
\limN J_i^N ( \huio,\humio) = \limN J_i^N(\uio,\umio) \, \text{w.p.1},\, 1\leq i \leq N.
\end{equation}\IEEEQED
\end{prop}

%%%%%%%%%%%%%%%%%%%%%%%%%%%%%%%%%%%%%%%%%%%%%%%%%%%%%%%%%%%%%%
% Proposition
%%%%%%%%%%%%%%%%%%%%%%%%%%%%%%%%%%%%%%%%%%%%%%%%%%%%%%%%%%%%%%
\begin{prop}\label{prop:eq_costs_g} 
For the system \eqref{eqn:dynamics}, under \ass{ass:ind}-\ass{ass:cost_coup}, \ass{ass:pdf}, \ass{ass:MLE_ident}, $u_i \in \U_g^N,\, u_i^0 \in \U_{MF}^N,\, \huio \in \hat{\U}_{MF}^N,\, 1\leq i \leq N$, the following holds:
\begin{equation}\label{eqn:eq_costs_g}
\limN\inf_{u_i \in \U_g^N} J_i^N(u_i,\humio) = \limN\inf_{u_i\in\U_g^N}J_i^N(u_i,\umio ) \quad
\text{w.p.1},\quad 1\leq i \leq N.
\end{equation}\IEEEQED
\end{prop}

% The result is proved in Appendix \ref{thm:CESACMF_pro}.

% Therefore we have shown that the cost obtained when each agent applies the MF-SAC Law in a system of adaptive agents all of whom are applying the same law, the cost obtained is almost surely equal to the cost obtained in a system of non-adaptive agents where all agents apply the non-adaptive MF Law. In the next lemma we show that the MF-SAC Law is indeed the optimal solution with respect to $\U_g^N$ when all other agents apply the same law. 

%===============================================================================
\section{Simulations}\label{sec:sims}

\begin{figure}[ht]
\centering
\includegraphics[width=\hsize]{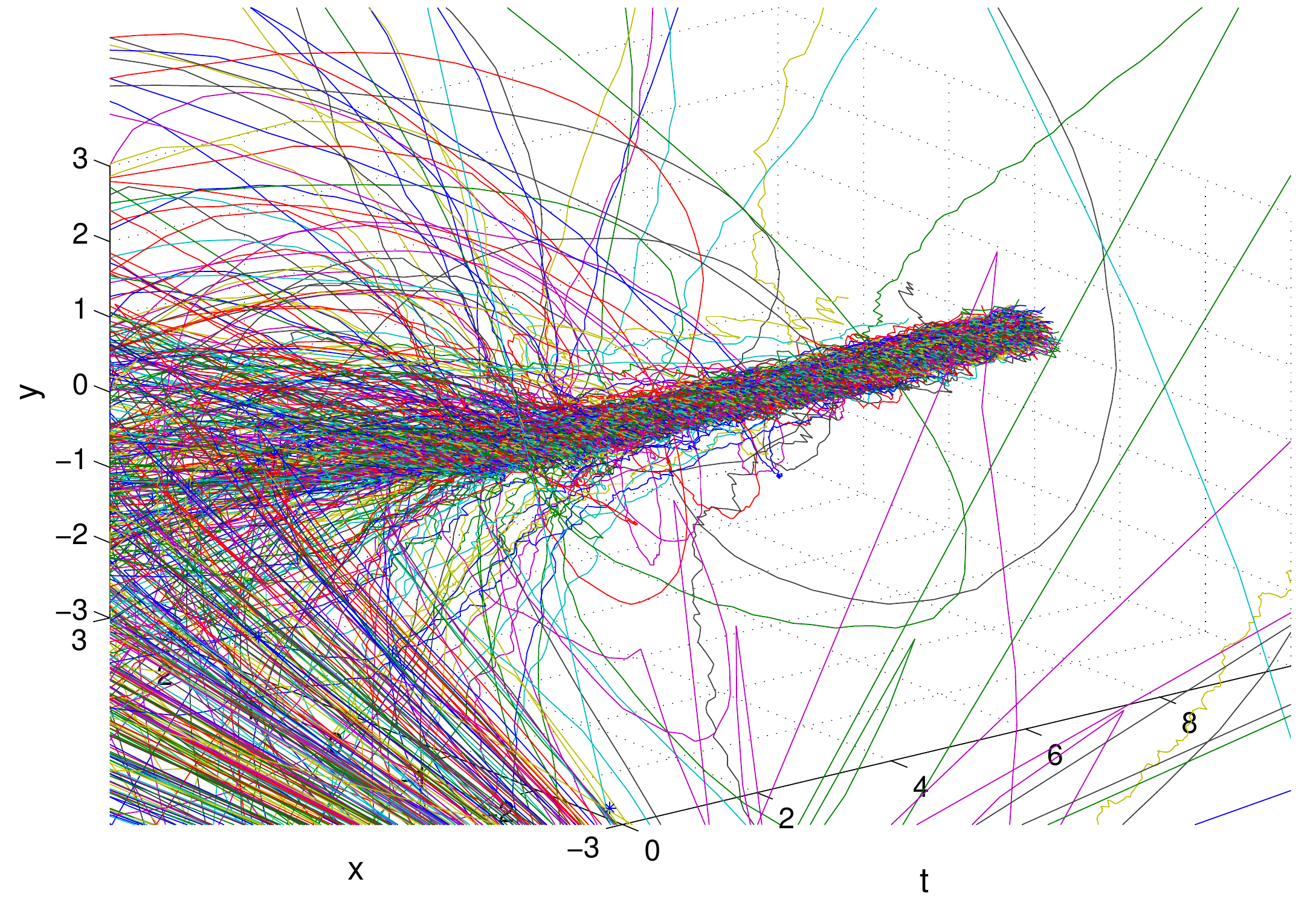}
\caption{State Trajectories}
\label{fig:traj}
\end{figure}

\begin{figure}[ht]
\begin{minipage}[b]{0.45\hsize}
\centering
\includegraphics[width=\hsize,height=3cm]{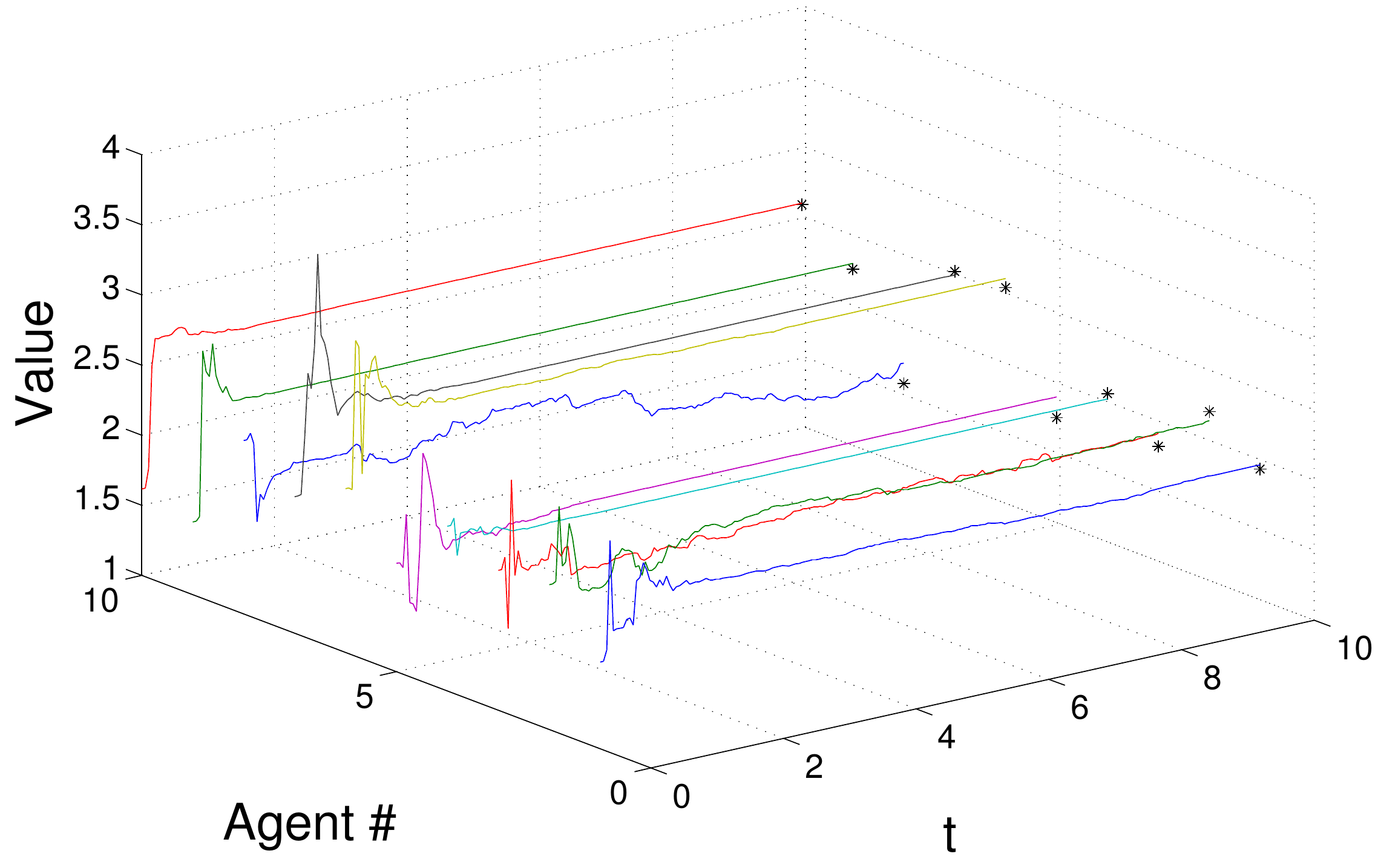}
\caption{Self Dynamical Parameter Identification}
\label{fig:A}
\end{minipage} 
\begin{minipage}[b]{0.45\hsize}
\centering
\includegraphics[width=\hsize,height=3cm]{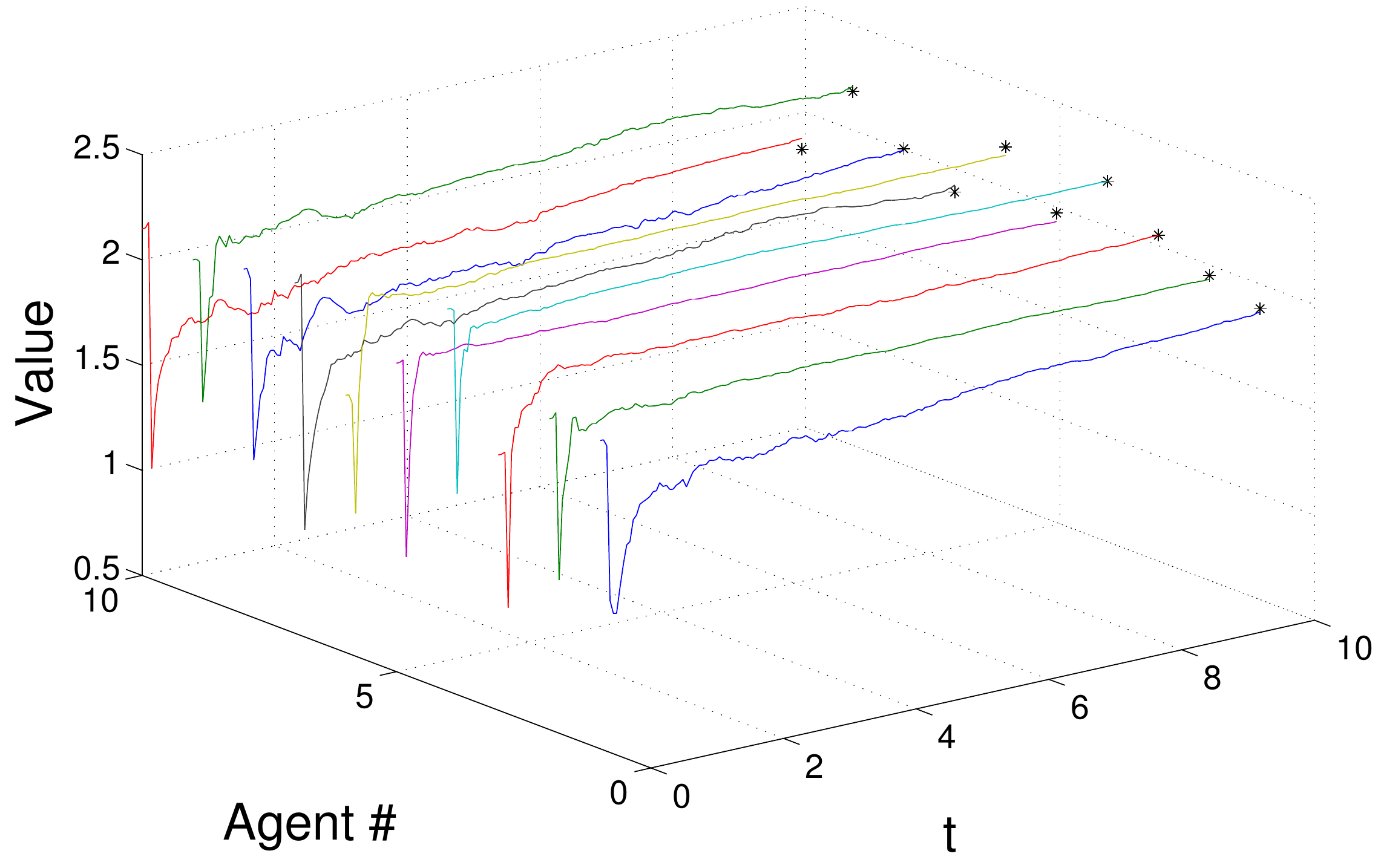}
\caption{Population Parameter Identification}
\label{fig:popn}
\end{minipage}
\end{figure}

Consider a system of 400 agents where each agent is modeled by a 2 dimensional system. 
% v3 The system matrices $\{\b{A}_k \}, \{ \b{B}_k \}, \{\b{Q}_k\},\, 1\leq k \leq 400$, are defined as
% \[ \b{A} \triangleq \left[ \begin{array}{cc}
%         -0.2+a_{11} & -2+a_{12} \\
%         1+a_{21} & 0.3 + a_{22} \end{array} \right], \quad \b{B} \triangleq \left[ \begin{array}{c}
%         1+b_{1}\\
%         1+b_{2} \end{array} \right],\quad \b{Q} \triangleq \left[ \begin{array}{cc}
%         2+q_{11} & 0.5+q_{12} \\
%         0.5+q_{12} & 2 + q_{22} \end{array} \right].  \]
% The population dynamical and cost function distribution parameters $a_{ij}$'s, $b_i$'s and $q_{ij}$'s in the simulation are mutually independent and distributed according to $\{a_{ij} \sim N(0,0.2),\, 1\leq i,j \leq 2\}$, $\{b_{i} \sim N(0,0.2),\, 1\leq i \leq 2\}$, and $\{q_{11},q_{12},q_{22} \sim N(0,0.2) \}$. 
All agents apply the MF-SAC Law; each of 400 agents observes its own 20 randomly chosen agents' outputs and control inputs, as well as its own trajectory. Rapid convergence of the state trajectories of all agents to the steady state values can be seen Fig.\ \ref{fig:traj} where `x' and `y' represent the two dimensions of each agent's state and `t' denotes time. In order to plot the convergence of the self identification of dynamical parameters $\b{A}_i,\, 1\leq i \leq N$, we plot the norm trajectories of the estimates in Fig.\ \ref{fig:A}. The symbol `*' denotes the true value of the parameter for each agent. Only 10 randomly chosen agents are shown in Fig.\ \ref{fig:A} for clarity of presentation. In Fig.\ \ref{fig:popn}, we depict each agent's estimate of the mean of the dynamical parameter $\b{A}$ (i.e., the mean of the random variable $\b{A}$), and we display 10 randomly chosen agents' estimate trajectories for clarity. Again, the norm of the estimates and the true values are displayed in this diagram. The resulting parameter estimate is different for each agent due to the fact that each agent only observes 20 randomly chosen agents out of a system of 400 agents.
%v3
% For the cost function parameter distribution identification, agent $A_i$ solves the RWLS equations \eqref{eqn:nslppe} for $\b{Q}_j, \, j\in Obs_i(N)$. The convergence of the norm of cost function parameter estimates of 10 randomly chosen agents in $Obs_i(N)$ is presented in Fig.\ \ref{fig:Q}.
% In Fig.\ \ref{fig:hist_A}, the histogram of the norm of the true values of $\b{A}_i,\, 1\leq i \leq N$, and the histogram of the norm of the self estimated values at the final instant of the simulation are shown.  

%and the cost function parameter $\b{Q}_i,\, 1\leq i \leq N$


%v3
%\begin{figure}[ht]
%\begin{minipage}[b]{0.49\hsize}
%\centering
%\includegraphics[width=\hsize]{TAC_Q_estimation}
%\caption{Cost Function Parameter Identification}
%\label{fig:Q}
%\end{minipage}
%\begin{minipage}[b]{0.49\hsize}
%\centering
%\includegraphics[width=\hsize]{TAC_A_mean_estimation_10}
%\caption{Population Parameter Identification}
%\label{fig:popn}
%\end{minipage} 
%\end{figure}
%\ack{one is value, the other one is values}

%===============================================================================
\section{Conclusion}\label{sec:con}

This paper presents a study of the mean field stochastic adaptive control problem where the cost functions of the agents in a population are coupled, and each agent estimates its own dynamical parameters based upon observations of its own trajectory, and furthermore estimates the distribution parameter of the population's dynamical and cost function parameters by observing a randomly chosen fraction of the population. This work makes a contribution to the mean field literature by extending the established $\e$-Nash equilibrium results of a large population of egoistic agents to a large population of adaptive egoistic agents. The information requirement for each agent is kept limited in the sense that the distribution parameter is estimated only through an observed set of agents, where the ratio of the cardinality of the observed set to the number of agents in the population becomes negligible as the population size grows to infinity. The strong consistency of the self parameter estimates and the distribution parameter estimates, the stability of the all agent systems, and an $\e$-Nash Equilibrium property are all established in the paper.

Future research directions include: (i) investigation of the influence of various rates of observed population fraction decay and rates of convergence on the results in this paper, together with (ii) the extension of adaptive MF theory to (a) the currently developing areas of distance dependent cost function influence among agents \cite{2010HCM_TAC}, (b) altruist and egoist MF theory \cite{2010HCM_CDC} and (c) problems involving partially observed systems.
% Another challenging line of future research is to study full Bayesian game problems within the MF framework.

%===============================================================================
% Appendix
\appendices

%===============================================================================
\section{}%{Dynamical Properties of Parameter Convergent Systems} 
%===============================================================================
% Unless otherwise stated $\b{A},\b{A}_t,\hat{\b{A}}_t$, etc, denote bounded real $n\times n$ matrix functions of $t\geq 0$, and analogously for the real $n\times m$ matrices $\b{B},\b{B}_t,\hat{\b{B}}_t$, etc. In case $\b{A},\hat{\b{A}}_t,\b{B},\hat{\b{B}}_t$ are random processes, they are bounded on almost all sample paths.
% Here we present the lemmas that are needed to prove the main result of the paper.
\subsection*{Preparatory Lemmas on Asymptotic Dynamics and Dither Inputs}\label{ssec:pro_exc_pro}
Four basic properties to be used in the sequel are given in the following lemmas.
\begin{lem}\label{lem:state_transition}
Let $\b{A}(\w)$ be an asymptotically stable random matrix on $\ps$ for all $\w \in \W$ except on a $P-$null set $\mathcal{N}$, and $\b{A}_t(w),t\geq 0$, be a bounded random matrix function of $t\geq 0$. If for all $\w \in \W \backslash \mathcal{N}$ and all $\e>0$, there exists $T_\w = T_\w(\e)$ such that $t>T_\w$ implies $\lVert \b{A}_t - \b{A} \rVert < \e$, i.e. $\b{A}_t \conv \b{A}$ w.p.1 as $t\conv \infty$, then $( \b{A}_t,t\geq 0 )$ is an exponentially stable time varying matrix w.p.1, in the sense that its fundamental matrix satisfies the estimate $\lVert \b{\Phi}_{t,s} \rVert \leq \beta e^{-\rho (t-s)}$ for $t_0 \leq s \leq t$, where $\rho(\w)>0$ and $0<\beta(\w)<\infty$. 
\end{lem}
%===============================================================================
%%%%%%%%%%%%%%%%%%%%%%%%%%%%%%%%%%%%%%%%%%%%%%%%%%%%%%%%%%%%%%%%%%%%%%%%%%%%%%%%
% Proof of Lemma lem:state_transition
%%%%%%%%%%%%%%%%%%%%%%%%%%%%%%%%%%%%%%%%%%%%%%%%%%%%%%%%%%%%%%%%%%%%%%%%%%%%%%%%
%\subsection{Proof of Lemma \ref{lem:state_transition}}\label{lem:state_transition_pro}
 \begin{IEEEproof}
  
 Suppressing mention of $\w\in \W \backslash \mathcal{N}$, whenever possible for simplicity of notation we consider,  
 \begin{align}
 \dot x(t) & =\mathbf{A} x(t), \quad t\geq 0; \quad x(0)=x_0 \in \R^n,\text{ and}\\
 \dot x_a(t) & =\mathbf{A}_t x_a(t), \quad t\geq 0; \quad x_a(0)=x_a \in \R^n.
 \end{align}
 Since $\mathbf{A}$ is asymptotically stable, we may form the Lyapunov function $V(x)=x\t \mathbf{\Pi} x$, where $\mathbf{\Pi}>0$ satisfies $\mathbf{\Pi} \mathbf{A} + \mathbf{A}\t \mathbf{\Pi} = -\mathbf{Q}$ for some $\mathbf{Q}>0$.
  
 Now, 
 \begin{align} 
 \dot V\left(x(t)\right) = & \dot x\t(t) \mathbf{\Pi} x(t) + x\t(t)\mathbf{\Pi} \dot x(t)\\
 = & x\t(t)[\mathbf{\Pi} \mathbf{A} + \mathbf{A}\t \mathbf{\Pi}]x(t)\\ 
 = & -x\t(t)\mathbf{Q}x(t), \quad t \geq 0.
 \end{align}
 Then writing
 \be \mathbf{A}_t\t \mathbf{\Pi} + \mathbf{\Pi} \mathbf{A}_t = (\mathbf{A}_t-\mathbf{A})\t\mathbf{\Pi} + \mathbf{A}\t \mathbf{\Pi} + \mathbf{\Pi}(\mathbf{A}_t-\mathbf{A})+\mathbf{\Pi} \mathbf{A},\quad t \geq 0, \ee
 we see that for all $\w \in \W \backslash \mathcal{N}$ there exists sufficiently large $T_\w$ such that for all $t>T_\w$,
  
 \be\mathbf{A}_t\t \mathbf{\Pi} + \mathbf{\Pi} \mathbf{A}_t < -\mathbf{Q}+\frac{\mathbf{Q} } {2} = \frac{-\mathbf{Q}} {2}. \ee
  
 Therefore, \be \dot V_a(x_a(t)) = \frac{d}{dt} \left [ x_a\t(t)\mathbf{\Pi} x_a(t)\right ] < -x_a\t(t)\frac{\mathbf{Q}}{2}x_a(t) < -\left (\frac{\lambda_{max} (\mathbf{Q}) } {2\lambda_{min}(\mathbf{\Pi})}\right ) ( x_a\t(t) \mathbf{\Pi} x_a(t)) < 0, \ee which implies $ \dot V_a(x_a(t)) < -a V_a[x_a(t)],$ where $a \triangleq \left ( \frac {\lambda_{max} (\mathbf{Q})} {2\lambda_{min}(\mathbf{\Pi})} \right )$, which gives
  \be V_a(x_a(t)) \leq V_a(x_a(t_0)) e^{-a(t-t_0)}. \ee
  
 Now, for the fundamental matrix, we have
 \be \lVert \mathbf{\mathbf{\Phi}}(t,t_0) \rVert_0 = \sup_{x_{t_0} \neq 0} \frac{\lVert \mathbf{\Phi}(t,t_0)x_{t_0} \rVert }{ \lVert x_{t_0} \rVert}= \sup_{x_{t_0} \neq 0} \frac{\sqrt{\lVert x_t \rVert^2}}{\lVert x_{t_0} \rVert}.\ee
  
  Without loss of generality, take $\lVert x_{t_0} \rVert=1$. Then,  
 \begin{align} 
 \lVert \mathbf{\Phi}(t,t_0) \rVert_0 = & \sup_{x_{t_0}} \sqrt{\lVert x_t \rVert^2}\\
 & \leq \sup_{x_{t_0} } \left ( \frac{V(x_t)}{\lambda(\mathbf{\Pi})_{min}} \right ) ^{\frac{1}{2}} \\
 & \leq \frac{1}{\sqrt{\lambda(\mathbf{\Pi})_{min}}}\sup_{x_{t_0}} (V(x_t))^{\frac{1}{2}}\\
 & \leq \frac{1}{\sqrt{\lambda(\mathbf{\Pi})_{min}}}\sup_{x_{t_0} } (e^{-a (t-t_0)}V(x_{t_0}))^{\frac{1}{2}}\\
 & \leq \frac{1}{\sqrt{\lambda(\mathbf{\Pi})_{min}}} e^{-(a/2) (t-t_0)} \sup_{x_{t_0} } (x_{t_0}\t \mathbf{\Pi} x_{t_0})^{\frac{1}{2}}\\
 & \leq e^{-\rho (t-t_0)} \left (  \frac{\lambda(\mathbf{\Pi})_{max}}{\lambda(\mathbf{\Pi})_{min}}   \right )^{\frac{1}{2}}\\
 & \leq \beta e^{-\rho (t-t_0)}, \enspace t\geq t_0; \text{ when }\rho=\frac{a}{2} \text{ and }\beta=\sqrt{\frac{\lambda(\mathbf{\Pi})_{max}}{\lambda(\mathbf{\Pi})_{min}}}.
 \end{align}
 \end{IEEEproof}
%%%%%%%%%%%%%%%%%%%%%%%%%%%%%%%%%%%%%%
\begin{lem}\label{lem:state_transition_conv}
Let $( \b{A}_t,t\geq 0 )$ be a random bounded matrix sequence on $(\W,\F,P)$, which converges almost surely to the asymptotically stable matrix $\b{A}^o$ as $t \conv \infty$; let $\b{\Psi}_{t,t_0}$ be defined by $ \frac{d}{dt}\b{\Psi}_{t,t_0}=\b{A}^o\b{\Psi}_{t,t_0}$, i.e.~ $\b{\Psi}_{t,t_0}=e^{\b{A}^o(t-t_0)}$, and let $\frac{d}{dt}\b{\Phi}_{t,t_0}=\b{A}_t\b{\Phi}_{t,t_0}$ with $\b{\Psi}_{t_0,t_0}= \b{\Phi}_{t_0,t_0} = \b{I}$. Then the following limit holds:
\begin{equation*}%\label{eqn:state_tran_int}
\lim_{T\conv \infty}\frac{1}{T} \int_{t_0}^T \left \lVert    \left ( \b{\Phi}_{t,t_0} - \b{\Psi}_{t,t_0} \right ) \right\rVert^2 dt = 0 \quad \text{w.p.1}.
\end{equation*}
\end{lem}
% The result is proved in Appendix \ref{lem:state_transition_conv_pro}.
%===============================================================================
%%%%%%%%%%%%%%%%%%%%%%%%%%%%%%%%%%%%%%%%%%%%%%%%%%%%%%%%%%%%%%%%%%%%%%%%%%%%%%%%
% Proof of Lemma lem:state_transition_conv
%%%%%%%%%%%%%%%%%%%%%%%%%%%%%%%%%%%%%%%%%%%%%%%%%%%%%%%%%%%%%%%%%%%%%%%%%%%%%%%%
%\subsection{Proof of Lemma \ref{lem:state_transition_conv}}\label{lem:state_transition_conv_pro}

 \begin{IEEEproof}
  
The proof is given in four steps below.

 %%%%%%%%%% (i) %%%%%%%%%
 (i) \textit{Integral Representation $I^T$:}
  
 For almost all $\w \in \W$, we have $\mathbf{A}_t(\w) \conv \mathbf{A}^0(\w)$ as $t\conv \infty$, restricting attention to the probability 1 subset of $\W_0 \subset \W$ on which a unique solution exists. Since
 \begin{equation}\label{eqn:difeqns} 
 \frac{d}{dt}\mathbf{\Psi}_{t,t_0}=\mathbf{A}^0\mathbf{\Psi}_{t,t_0} \text{ with } \mathbf{\Psi}_{t_0,t_0}=\b{I}, \quad \text{ and } \quad \frac{d}{dt}\mathbf{\Phi}_{t,t_0}=\mathbf{A}_t\mathbf{\Phi}_{t,t_0} \text{ with }\b{\Phi}_{t,t_0}=\b{I},
 \end{equation}
  
 %\mathbf{\Psi}_{t_0,t_0}=\mathbf{\Psi}_{t_0,t_0}=I\\
  
 we have
 \begin{align}
 \label{eqn:est_difference} \frac{d}{dt}(\mathbf{\Phi}_{t,t_0} - \mathbf{\Psi}_{t,t_0}) = & \mathbf{A}_t\mathbf{\Phi}_{t,t_0} - \mathbf{A}^0 \mathbf{\Psi}_{t,t_0}\\
 = & \mathbf{A}_t\mathbf{\Phi}_{t,t_0} - (\mathbf{A}^0-\mathbf{A}_t)\mathbf{\Psi}_{t,t_0} - \mathbf{A}_t\mathbf{\Psi}_{t,t_0}\\
 = & \mathbf{A}_t(\mathbf{\Phi}_{t,t_0} -\mathbf{\Psi}_{t,t_0})-(\mathbf{A}^0- \mathbf{A}_t)\mathbf{\Psi}_{t,t_0}.
 \end{align}
  
 Integrating we obtain,
  
 \begin{equation} \mathbf{\Phi}_{t,t_0} - \mathbf{\Psi}_{t,t_0}=\mathbf{\Phi}_{t,t_0}(\mathbf{\Phi}_{t_0,t_0}-\mathbf{\Psi}_{t_0,t_0})-\int_{t_0}^t \mathbf{\Phi}_{t,s}(\mathbf{A}^0-\mathbf{A}_s)e^{\mathbf{A}^0(s-t_0)}ds, \end{equation}
  
 with the initial condition $ \mathbf{\Phi}_{t_0}-\mathbf{\Psi}_{t_0} = \mathbf{I} - \mathbf{I} =0 $. Therefore, 
 \begin{equation}\label{eqn:wtic}
 \mathbf{\Phi}_{t,t_0} - \mathbf{\Psi}_{t,t_0} = -\int_{t_0}^t \mathbf{\Phi}_{t,s}(\mathbf{A}^0-\mathbf{A}_s)e^{\mathbf{A}^0(s-t_0)}ds, \enspace t\geq t_0,
 \ee
 and so,
 \begin{equation}\label{eqn:state_tran_int2_detail}
 I^T \triangleq \lim_{T\conv \infty}\frac{1}{T} \int_{t_0}^T \left \lVert    \left ( \mathbf{\Phi}_{t,t_0} - \mathbf{\Psi}_{t,t_0} \right ) \right\rVert^2 dt = \frac{1}{T} \int_{t_0}^T \left \lVert  \int_{t_0}^t \mathbf{\Phi}_{t,s}(\mathbf{A}^0-\mathbf{A}_s)e^{\mathbf{A}^0(s-t_0)}ds \right \rVert^2 dt.
 \end{equation}
  
 %%%%%%%%%% (ii) %%%%%%%%%
 (ii) \textit{$I^T = I_1^T + I_2^T;$ Convergence of $I_1^T$:}
  
 Let us split the integrals in \eqref{eqn:state_tran_int2_detail} as follows:
 \begin{equation}\label{eqn:state_tran_sep_int2}
 I^T = \frac{1}{T} \left ( \int_{t_0}^{T_w} \left \lVert . \right \rVert^2 dt +  \int_{T_w}^T \left \lVert . \right \rVert^2 dt \right ) =: I_1^T + I_2^T,
 \end{equation}
 where the inner integrals are defined by `$\cdot$' for brevity in this definition and $T_\w>t_0$ is a random instant whose value is to be determined later.
  
 We take the norm inside the integral in $I_1^T$; then by use of the Cauchy Schwarz Inequality (henceforth termed CS) we may bound $I_1^T$ above as in
 \be 
 I_1^T \leq \frac{1}{T}  \int_{t_0}^{T_\w} \left ( (t-t_0) \int_{t_0}^t \lVert  \mathbf{\Phi}_{t,s} e^{\mathbf{A}^0(s-t_0)} \rVert^2 \lVert \tilde{\mathbf{A}}_s  \rVert^2 ds \right ) dt, \quad t_0\leq t \leq T_\w, 
 \ee
 where $ \tilde{\mathbf{A}}_s = (\mathbf{A}^0 - \mathbf{A}_s),\,s\geq 0.$
  
 Next, we may bound $\lVert \tilde{\mathbf{A}}_s \rVert$ above by some $M_{\tilde A}^{T_\w}$ for $t_0\leq s \leq T_w$, and we may bound $\lVert e^{\mathbf{A}^0(s-t_0)} \rVert$ by $\beta_0 e^{-\rho_{A^0}(s-t_0)}$, for some $\beta_0 >0$. Moreover $\lVert \mathbf{\Phi}_{t,s} \rVert \leq M_\Phi^{T_\w}$, for all $t_0\leq s \leq t \leq T_\w$, for some $M_\Phi^{T_\w}<\infty$, by the continuity of solutions to \eqref{eqn:difeqns}.
  
 Then, 
 \begin{equation}\label{eqn:rnd1}
 \limsup _{T\conv \infty} I_1^T \leq \limsup_{T\conv \infty}\frac{\beta_0^2}{T} \left ( M_{\tilde{A}}^{{T_\w}}\right)^2 \left ( M_{\Phi}^{{T_\w}} \right)^2   \int_{t_0}^{T_\w} (t-t_0) \left (   \int_{t_0}^t   e^{-2\rho_{A^0}(s-t_0)}  ds \right ) dt =: \lim_{T\conv \infty} \frac{1}{T} \kappa g(t_0,T_\w), 
 \end{equation}
 where $\kappa=\beta_0^2 M_{\tilde{A}}^{{T_\w}^2} M_{\Phi}^{{T_\w}^2}<\infty$ and $g(\cdot)$ is a bounded continuous function of $t_0$ and $T_\w$. Hence for a fixed $T_\w$, $\limsup_{T\conv \infty} \frac{1}{T} \kappa g(t_0,T_\w) = 0 $. Therefore $I_1^T$ tends to 0 as $T$ tends to $\infty$.
  
 %%%%%%%%%% (iii) %%%%%%%%%
  
 (iii) \textit{$I^T = I_1^T + I_2^T;$ Convergence of $I_2^T:$}
  
 For the second integral $I_2^T$ in \eqref{eqn:state_tran_sep_int2}, we have
  
 \begin{align}
 I_2^T = & \frac{1}{T} \int_{T_\w}^T \left \lVert  \int_{t_0}^t \mathbf{\Phi}_{t,s}(\mathbf{A}^0-\mathbf{A}_s)e^{\mathbf{A}^0(s-t_0)}ds \right \rVert^2 dt\\
 & \leq \frac{1}{T}  \int_{T_w}^T  \left ( (t-t_0)   \int_{t_0}^t \lVert  \mathbf{\Phi}_{t,s} e^{\mathbf{A}^0(s-t_0)} \rVert^2 \lVert \tilde{\mathbf{A}}_s \rVert^2 ds \right ) dt , \quad t_0 \leq T_\w \leq t \leq T,\\
 & = \frac{1}{T} \int_{T_w}^T     \left ( (t-t_0)\int_{t_0}^{T_w} \lVert  . \rVert^2 ds + (t-t_0) \int_{T_w}^t \lVert  . \rVert^2 ds \right ) dt =: I_{21}^T + I_{22}^T,\end{align}
 where we split the inner integral and use the $(\cdot)$ notation for brevity.
  
 Using the semi-group property of the state transition matrix, we may write $\mathbf{\Phi}_{t,s}=\mathbf{\Phi}_{t,T_w}\mathbf{\Phi}_{T_w,s}$ for all $t_0\leq s \leq T_\w < T$. But we have $\sup_{t_0\leq s \leq T_\w}\lVert \mathbf{\Phi}_{T_\w,s} \rVert =: M_{\Phi}^{T_\w} < \infty$, and we have $M_{\tilde{A}}^{T_\w} := \sup_{t_0\leq s \leq T_\w}\lVert \tilde{\b{A}}_s \rVert$. Therefore,
  
 \be 
 I_{21}^T \leq \frac{1}{T}   \left ( M_{\tilde{\b{A} } }^{T_\w}\right)^2 \left ( M_{\mathbf{\Phi}}^ {T_w} \right )^2 \int_{T_w}^T     \left ( (t-t_0)\int_{t_0}^{T_w} \lVert  \mathbf{\Phi}_{t,T_w} e^{\mathbf{A}^0(s-t_0)} \rVert^2 ds \right ) dt, \quad t_0\leq s \leq T_\w \leq t \leq T.
 \ee
  
 Concerning $I_{22}^T$, the random time $T_\w$ is chosen so that for $s\geq T_w$ (increasing the value over that used in \eqref{eqn:rnd1} without affecting that argument), $\lVert \tilde{\mathbf{A}}_s \rVert < \e$. Hence,
 \be  
 I_{22}^T \leq \frac{1}{T}  \e^2 \int_{T_w}^T  \left ( (t-t_0) \int_{T_w}^t \lVert  \mathbf{\Phi}_{t,s}  e^{\mathbf{A}^0(s-t_0)} \rVert^2 ds \right ) dt,\quad t_0 \leq T_\w\leq s \leq t \leq T.  
 \ee
  
 From Lemma \ref{lem:state_transition}, $\mathbf{\Phi}_{t,t_0}$ satisfies the bound  $\lVert \mathbf{\Phi}_{t,t_0} \rVert \leq \beta_1 e^{-\rho_{\Phi} (t-t_0)},t \geq t_0$, where $\beta_1=\beta_1(\w),\rho = \rho(\w)$. Finally, bounding $\lVert e^{\mathbf{A}^0(s-t_0)} \rVert$ by $\beta_0 e^{-\rho_{A^0}(s-t_0)}$ yields
 \begin{multline} 
 I_{21}^T \leq \frac{\beta_0^2}{T}  \left ( M_{\tilde{ A } }^{T_\w}\right)^2 \left(M_{\Phi}^{T_w}\right)^2 \beta_1^2  \int_{T_w}^T     \left ( (t-t_0)\int_{t_0}^{T_w}   e^{-2\rho_\Phi (t-T_\w)} e^{-2\rho_{A^0} (s-t_0)}  ds \right ) dt,\\ \quad t_0 \leq s \leq T_\w \leq t \leq T. \end{multline}

 For simplicity, in the following we use $\rho = \min \left [ \rho_\Phi,\rho_{A^0}\right ]$ and $\beta = \max \left [ \beta_0,\beta_1\right ]$; then,
 \begin{align}
 I_{21}^T & \leq \frac{1}{T}  \kappa' \int_{T_w}^T     \left ( (t-t_0)\int_{t_0}^{T_w}   e^{-2\rho (t-T_\w)} e^{-2\rho (s-t_0)}  ds \right ) dt, \quad t_0 \leq s \leq T_\w \leq t \leq T,\\
 & \qquadeight \text{where } \kappa' = \beta^4 \left ( M_{\tilde{ A } }^{T_\w}\right)^2 \left(M_{\Phi}^{T_w}\right)^2,\\
 I_{21}^T & \leq \limsup_{T\conv \infty}\frac{1}{T} \kappa'\kappa'' e^{2\rho T_\w}\left( \frac{T}{e^{2\rho T}} + \frac{t_0}{e^{2\rho T}}+ \frac{t_0}{e^{2\rho T_\w}} + \frac{T_\w}{e^{2\rho T_\w}}  \right),
 \end{align}
 for a suitable constant $\kappa''$ independent of $T_\w$, which tends to 0 as $T \conv \infty$.
  
 %%%%%%%%%% (iv) %%%%%%%%%
 (iv) \textit{$I_2^T = I_{21}^T + I_{22}^T;$ Convergence of $I_2^T$:} Employing the hypothesis $\b{A}_t \conv \b{A}^0$ w.p.1 as $t\conv\infty$, we shall fix $T_\w$ such that $\lVert \tilde{\b{A}}_s \rVert<\e$
  
 For $I_{22}^T$, applying Lemma \ref{lem:state_transition} for $T_\w\leq s \leq t \leq T$, we obtain
 \begin{align} 
 I_{22}^T & \leq \frac{1}{T}  \e^2 \beta^4 \int_{T_w}^T   \left ( (t-t_0)\int_{T_w}^t   e^{-2\rho_{\mathbf{\Phi}} (t-s)}  e^{-2\rho_{\mathbf{A}^0}(s-t_0)}  ds \right ) dt,\quad t_0 \leq T_\w \leq s \leq t \leq T,\\  
 & \leq \frac{1}{T} \e^2 \beta^4 \int_{T_w}^T   \left ( (t-t_0)\int_{T_w}^t   e^{-2\rho (t-t_0)}  ds \right ) dt,\qquad \rho:=\min[\rho_{\mathbf{\Phi}},\rho_{\mathbf{A}^0}],\\
 & \leq \frac{g(T_\w)}{T} + \e^2\beta^2\exp(-2\rho(T-T_\w))\frac{2\rho^2T^2+2\rho T+1}{4\rho^3 T},
 \end{align}
 where $g(\cdot)$ is a bounded continuous function. Then, $\lsT I_{22}^T = 0$. Therefore $\limsup_{T\conv \infty}I_2^T \leq \limsup_{T\conv \infty} I_{21}^T + \lsT I_{22}^T =0$. 
 
 Since we have established that $I_1^T \conv 0, I_2^T \conv 0$ w.p.1 as $T\conv \infty$, we obtain $\limsup_{T\conv \infty} I^T \leq \lsT I_1^T + \lsT I_2^T=0$.

 Hence, we have proved that
 \be
 \lim_{T\conv \infty}\frac{1}{T} \int_{t_0}^T \left \lVert   \mathbf{\Phi}_{t,t_0} - \mathbf{\psi}_{t,t_0}  \right\rVert^2 dt = 0,\enspace w.p.1.
 \ee
 %\be
 %\lim_{T\conv \infty}\frac{1}{T} \int_{t_0}^T \left \lVert  f(t)  \left  ( \mathbf{\Phi}_{t,t_0} - \mathbf{\psi}_{t,t_0} \right ) \right\rVert^2 dt = 0,
 %\ee
 %which by the Cauchy-Schwarz inequality implies
 %\be
 %\lim_{T\conv \infty}\frac{1}{T} \int_{t_0}^T \left \lVert f(t) \left  (\mathbf{\Phi}_{t,t_0} - \mathbf{\psi}_{t,t_0} \right ) \right  \rVert dt = 0.
 %\ee
  \end{IEEEproof}
%%%%%%%%%%%%%%%%%%%%%%%%%%%%%%%%%%%%%%%%%%%%
% Lemma: A Property of the Excitation Process
%%%%%%%%%%%%%%%%%%%%%%%%%%%%%%%%%%%%%%%%%%%%
\begin{lem}\label{lem:kronecker} \cite[Duncan et al.\ (1999)]{1999DGD_TAC} 
Assume the process $(\e(t),t\geq 0)$ is an $\R^m$-valued standard Wiener process that is independent of $(w(t),t\geq 0)$, and assume the countable set of random processes $\{ (\e(t+k)-\e (k),t\in(0,1]);k\in \N\}$ to be mutually independent and all members of the set have the same probability law on $(0,1]$. Then, for all $f(\cdot)\in L^\infty[0,\infty)$:
\begin{equation*}%\label{eqn:kronecker}
\lim_{T\conv \infty}\frac{1}{T} \int_0^T  \bigg\lVert  \sum _{k=0}^{\lfloor t \rfloor}  \Big( \int_k^{\min[t,k+1]}\\ f(\tau)\xi_k\left [\e(\tau)-\e(k) \right ] d\tau  \Big)  \bigg\rVert^2 dt = 0. \quad \text{w.p.1}.
\end{equation*}
\IEEEQED
\end{lem}
The proof is given in \cite[Lemma 5]{1999DGD_TAC}.
% ACK-LV
 %%%%%%%%%%%%%%%%%%%%%%%%%%%%%%%%%%%%%%%%%%%%
 % Lemma: Ergodicity
 %%%%%%%%%%%%%%%%%%%%%%%%%%%%%%%%%%%%%%%%%%%%
 \begin{lem}\label{lem:ergodicity}\cite[Chen and Guo (1991)]{1991CG} 
 Let $\b{A} \in \R^{n^2}$ be an asymptotically stable matrix, and let $\b{D}\in \R^{nr}$. Then
 \begin{equation*}%\label{eqn:ergodicity}
 \lsT\frac{1}{T}\int_0^T \left\lVert \int_0^t e^{\b{A} (t-\tau)}\b{D}dw(\tau)\right\rVert^2dt = \int_0^\infty \tr (e^{\b{A} t}\b{D}\b{D}\t e^{\b{A}\t t})dt.
 \end{equation*}
 % Therefore $\left\lVert \int_0^t e^{\b{A} (t-\tau)}\b{D} dw(\tau)\right\rVert^2$ is ergodic.
 \end{lem}
 %The result is proved in Appendix \ref{lem:ergodicity_pro} following \cite[Lemma 12.4]{1991CG}.
% ACK-LV
%===============================================================================
%%%%%%%%%%%%%%%%%%%%%%%%%%%%%%%%%%%%%%%%%%%%%%%%%%%%%%%%%%%%%%%%%%%%%%%%%%%%%%%%
% Proof of Lemma lem:ergodicity
%%%%%%%%%%%%%%%%%%%%%%%%%%%%%%%%%%%%%%%%%%%%%%%%%%%%%%%%%%%%%%%%%%%%%%%%%%%%%%%%

 \begin{IEEEproof}[Proof (after \cite{1991CG})]\label{lem:ergodicity_pro}

 Consider the stochastic differential equation 
 \be dx_t = \b{A} x_t d_t + \b{D}dw_t,\enspace t\geq 0.\ee
 Since $\b{A}$ is asymptotically stable, there exists a positive definite matrix $\b{\Pi} >0$ such that
 \begin{equation}
 \b{\Pi} \b{A} + \b{A}\t \b{\Pi} = -\b{I} .
 \end{equation}
 
 Following \cite{1991CG}, applying the It\^{o} formula to the Lyapunov function $x_t\t \b{\Pi} x_t,\, t\geq 0$,
 \begin{align}
 \label{eqn:dif_Lyap1}d[x_t\t \b{\Pi} x_t] & = x_t\t(\b{\Pi} \b{A} + \b{A}\t \b{\Pi})x_tdt +  \tr ( \b{\Pi} \b{D}\b{D}\t ) dt + 2 x_t\t \b{\Pi} \b{D} dw_t\\
 \label{eqn:dif_Lyap2} & = -\lVert x_t \rVert^2 dt + \tr (\b{\Pi} \b{D}\b{D}\t) dt + 2 x_t\t \b{\Pi} \b{D} dw_t.
 \end{align}

  Integrating \eqref{eqn:dif_Lyap2} and using the result in Lemma 4 of Christopeit \cite{1986Ch} to estimate the third term on the RHS of \eqref{eqn:dif_Lyap1}, we obtain
  \begin{equation}
  x_t\t  \b{\Pi} x_t \leq - \int_0^t \lVert x_s \rVert^2ds + \tr ( \b{\Pi} \b{D} \b{D}\t ) t + o \left (\left \{  \int_0^t \lVert x_s \rVert^2ds  \right \}^{1/2+\eta} \right ) + O(1),\quad \text{where } 0 <  \eta < \frac{1}{2} ,
  \end{equation}
  and hence, $\int_0^t \lVert x_s \rVert^2 ds = O(t)\quad \text{w.p.1}$.
  
 We apply the It\^{o} formula to the outer product $x_t x_t\t,\, t\geq 0$,
 \begin{equation}
 d[x_t x_t\t] =  x_t x_t\t \b{A}\t dt +\b{A} x_t x_t\t dt + \b{D}\b{D}\t dt + \b{D}dw_t x_t\t  + x_t dw_t\t\b{D}\t.
 \end{equation}
 
 Integrating the outer product $x_t x_t\t$ from $t=0$ yields
  \be x_t x_t\t = \left ( \int_0^t x_s x_s\t ds\right ) \b{A}\t + \b{A} \left ( \int_0^t x_s x_s\t ds\right ) +  (\b{D}\b{D}\t)t + \int_0^t  (\b{D}dw_s x_s\t)  + \int_0^t  (x_s dw_s\t\b{D}\t).  \ee
  
 % ??? of Lyapunav ODE integration type
 A ``Lyapunav integral move" yields
 \begin{equation}\label{eqn:dif_Lyap3}
  \begin{aligned}
  & \int_0^t x_s x_s\t ds =\\ 
  & \qquad \int_0^t e^{\b{A} (t-s)} (\b{D}\b{D}\t) s e^{\b{A}\t (t-s)} ds + \int_0^t e^{\b{A} (t-s)} \left ( \int_0^s \left \{ (\b{D}dw_\tau x_\tau\t) +  (x_\tau dw_\tau\t\b{D}\t) \right \}d\tau \right ) e^{\b{A}\t (t-s)} ds.
  \end{aligned}
  \end{equation}
  
  We deal with the second term of RHS of \eqref{eqn:dif_Lyap3}. Using Christopeit's \cite{1986Ch} estimate again we write,
 \begin{align} 
 & \left \lVert \int_0^t  e^{\b{A} (t-s)} \left ( \int_0^s\left \{ (\b{D}dw_\tau x_\tau\t) +  (x_\tau dw_\tau\t\b{D}\t) \right \}d\tau \right ) e^{\b{A} (t-s)} ds \right \rVert\\
 & \qquad \leq \left ( \int_0^t e^{-2\rho (t-s)} \left ( o\left \{ \int_0^s \lVert x_\tau \rVert^2 d\tau \right \}^{1/2+\eta} + O(1) \right ) ds \right ),\quad \text{where}\enspace 0<\eta<\frac{1}{2},\\
 & \qquad =  \int_0^t e^{-2 \rho (t-s)} o \left (  s^{1/2+\eta} \right) ds = o(t^{1/2 +\eta} ),\quad \eta>0. 
 \end{align}
  
  As $\lim_{t\conv \infty} \frac{1}{t} \int_0^t e^{-2\rho(t-s)}s ds = \int_0^\infty e^{-2\rho s} ds$, we take the time average limit of \eqref{eqn:dif_Lyap3} and get
  \[  \lim_{t\conv \infty } \frac{1}{t} \int_0^t x_s x_s\t ds  = \int_0^\infty e^{\b{A} s}  ( \b{D}\b{D}\t ) e^{\b{A}\t s} ds, \]
  
  which implies
  \begin{equation}
  \lsT\frac{1}{T}\int_0^T \left\lVert \int_0^t e^{\b{A} (t-\tau)}\b{D}dw(\tau)\right\rVert^2dt = \int_0^\infty \tr (e^{\b{A} t}\b{D}\b{D}\t e^{\b{A}\t t})dt.
  \end{equation}
  
  Thus we obtain the desired result.\end{IEEEproof}

%===============================================================================
\section{}\label{sec:app_cppe}%{Convergence Properties of the Parameter Estimates}
%===============================================================================

%===============================================================================
%%%%%%%%%%%%%%%%%%%%%%%%%%%%%%%%%%%%%%%%%%%%%%%%%%%%%%%%%%%%%%%%%%%%%%%%%%%%%%%%
% Proof of Lemma lem:uni_con
%%%%%%%%%%%%%%%%%%%%%%%%%%%%%%%%%%%%%%%%%%%%%%%%%%%%%%%%%%%%%%%%%%%%%%%%%%%%%%%%
\subsection*{Proof of Lemma \ref{lem:uni_con}}\label{lem:uni_con_pro}
We drop the subscript $i$ for clarity. By definition, when $\hth_t$ is the solution to RWLS equations \eqref{eqn:nslspe}, $\hth_t^{pr}$ satisfies $\hth_t^{pr} \triangleq \argmin_{\psi\in \Th} \lVert \hth_t - \psi \rVert$, employing a co-ordinate ordering measurable tie breaking rule, if necessary. Since $\hth_t \in \R^{n\left (n+m+(n+1)/2 \right)},\, \hth_t^{pr} \in \Th$ and $\th^0 \in \Th$, the definition of $\hth_t^{pr}$ gives $\lVert \hth_t - \hth_t^{pr} \rVert \leq \lVert \hth_t - \th^0\rVert$. But by hypothesis, $\lVert \hth_t - \th^0\rVert \conv 0\enspace \text{w.p.1 as}\enspace t\conv \infty$; therefore, $\hth_t^{pr}  \conv \th^0$ w.p.1 as $t\conv \infty$.
\IEEEQED

%%%%%%%%%%%%%%%%%%%%%%%%%%%%%%%%%%%%%%%%%%%%%%%%%%%%%%%%%%%%%%%%%%%%%%%%%%%%%%%%
% Proof of Theorem thm:str_con_est
%%%%%%%%%%%%%%%%%%%%%%%%%%%%%%%%%%%%%%%%%%%%%%%%%%%%%%%%%%%%%%%%%%%%%%%%%%%%%%%%
\subsection*{Proof of Theorem \ref{thm:str_con_est}}\label{thm:str_con_est_pro}
(i) Since the solution $\b{\Pi}_\theta \in \R^{n^2},\, \th\in\R^{n(n+m+(n+1)/2)}$, to the algebraic Riccati equation \eqref{eqn:DE_Pi} parametrized by $\th\in \Th$ is a smooth function of $\th$ (see \cite{1984De_TAC}), and since $\hth_{i,t}^{pr}\in\Th,\, t\geq 0,\, 1\leq i \leq N$, $\b{\Pi}(\hth_{i,t}^{pr}),\,t\geq 0$, satisfies $\b{\Pi}(\hth_{i,t}^{pr})<\infty$ w.p.1 for all $t\geq 0$. It is given that $x^* \in \b{C}_b[0,\infty)$; therefore, $s(t;\hth_{i,t}^{pr},\hzn_{i,t})<\infty$ w.p.1 for $t\geq 0$ evaluated along $\hth_{i,t}^{pr},\, t\geq 0$. Hence, $\huio (t;\hth_{i,t}^{pr},\hzn_{i,t})$ given in \eqref{eqn:nslu} is well defined. 

(ii) The strong consistency of the dynamical parameter estimates $[\hat{\b{A}}_{i,t},\hat{\b{B}}_{i,t}],\, t\geq 0,\, 1\leq i \leq N,$ is shown in \cite{1999DGD_TAC} under the controllability and observability of the true parameters (A1 and A2 in \cite{1999DGD_TAC}) and the uniform controllability and observability of the estimates (Definition 1 in \cite{1999DGD_TAC}). In our work, the controllability and observability assumptions are satisfied since $[\b{A}_{\th},\b{B}_{\th}]$ is controllable and $[\b{Q}_{\th}^{1/2},\b{A}_{\th}]$ is observable for all $\th \in \Th$ by \ass{ass:conobs}, and moreover, the uniform controllability and observability of the estimates are satisfied due to Lemma \ref{lem:uni_con}. 

(iii) Dropping the subscript $i$ for clarity we set the estimation vector $\upsilon_t\t  =  [\hat{\b{\Pi}}_t,\hat{s}(t) ]$ and the regression vector as $\psi_t\t=[x_t\t,1]$. The persistence of excitation is satisfied since \[ \liminf_{T\conv\infty} \frac{1}{T} \lambda_{min} \left ( \int_0^T \psi_t \psi_t\t dt \right ) >0.\] Setting the measurement vector to be $ (-(\hat{\b{B}}_t\t)^{-1}\b{R}u_t)\t$ and employing \oass{ass:invB} we get $[\hat{\b{\Pi}}_t,\hat{s}(t) ]\conv [\b{\Pi}^0,s^0(t)]$ w.p.1 as $t\conv\infty$. Also, as shown in Part (i) $[\hat{\b{A}}_t,\hat{\b{B}}_t]\conv [\b{A}^0,\b{B}^0]$ w.p.1 as $t\conv\infty$. Each estimated parameter in $\hat{\b{Q}}_{t} = -\hat{\b{A}}_{t}\t \hat{\b{\Pi}}_{t}\t - \hat{\b{\Pi}}_{t} \hat{\b{A}}_{t} + \hat{\b{\Pi}}_{t}\t \hat{\b{B}}_{t} {\b{R}}^{-1} \hat{\b{B}}_{t}\t \hat{\b{\Pi}}_{t}$ converges to its true value as $t\conv\infty$. Hence, $\hat{\b{Q}}_t\conv \b{Q}^0$ w.p.1 as $t\conv\infty$.

We observe that instead of the random regularization method used in \cite{1996Gu_TAC} and \cite{1999DGD_TAC}, we employ here the projection method (Lemma \ref{lem:uni_con}), which guarantees the uniform controllability and observability of the estimates. 

\IEEEQED
%===============================================================================
%%%%%%%%%%%%%%%%%%%%%%%%%%%%%%%%%%%%%%%%%%%%%%%%%%%%%%%%%%%%%%%%%%%%%%%%%%%%%%%%
% Proof of Theorem thm:scmle
%%%%%%%%%%%%%%%%%%%%%%%%%%%%%%%%%%%%%%%%%%%%%%%%%%%%%%%%%%%%%%%%%%%%%%%%%%%%%%%%
\subsection*{Proof of Theorem \ref{thm:scmle}}\label{thm:scmle_pro}
Recall that $\th^{[1:N_0]}\teq \{\th_j;\,j\in Obs(N)\} $ is an independently selected subset of $\th^{[1:N]}$ of cardinality $N_0(N)$, and $\th_i,1\leq i \leq N$, is an independently distributed sequence with each $\th_i$ having the density $f_\z(\th)$, and hence $\th^{[1:N_0]}$ possesses a density in product form. Consequently, the scaled log-likelihood function of $\z$ at $\th^{[1:N_0]}$ is given by $L_{N_0}(\z)\teq L\left(\th^{[1:N_0]};\z\right) = -\frac{1}{N}\log\left(\prod_{j\in Obs(N)}f_\z(\th_{j})\right)$, $N_0\triangleq \lvert Obs(N) \rvert$. Note that the subscript $i$ is suppressed for clarity. The maximum likelihood estimate of $\z$ given $\th^{[1:N_0]}$ is then given by $\hzn = \argmin_{\z \in P} L(\th^{[1:N_0]};\z)$.

Now, it has been established in Theorem \ref{thm:str_con_est} that $(\hth_t^{[1:N_0]};\, t\geq 0)$ for each $N_0(N),N\in\Z_1$, constitutes a strongly consistent estimate of $\th^{[1:N_0]}$, i.e.,\ $\limt \hth_t^{[1:N_0]}= \th^{[1:N_0]}$ w.p.1. %In the probability space $(\Th,\B(\Th),P_\z)$, where $P_\z$ is the probability measure associated with the density function $f_\z,\, \z\in P$, the independence of the selection of $\th^{[1:N_0]}\in \Th$ results in the product form of the joint density of $\th^{[1:N_0]}$. 
Based upon this, the proof of the theorem consists of an analysis of the convergence (as $N\conv\f$ and hence $N_0(N)\conv\f$, and $t\conv\f$) of the likelihood function $ L(\th^{[1:N_0]};\z)$ with $\hth_t^{[1:N_0]}$ substituted for $ \th^{[1:N_0]}$ and hence of the associated sequence of estimators $(\hz^{N_0}_t;\, N_0\in\Z_1)$ to $\z^0$. 

First we present two lemmas that will be used in the sequel for the proof of Theorem \ref{thm:scmle}. 

\paragraph*{Convergence of the Likelihood Functions $L\left(\th^{[1:N_0]};\z\right)$}

\begin{lem}\label{lem:scmle1}
Subject to \ass{ass:pdf}, \ass{ass:MLE_ident} we have \[L_{N_0}(\z) \teq L\left(\th^{[1:N_0]};\z\right) \conv L_{\z^0}(\z) \teq -\E_{\z^0}\log f_\z(\th) \, \text{w.p.1},\] as $N_0\conv\infty$ uniformly over $\z\in P$.
\IEEEQED
\end{lem}
%v3 The proof of Lemma \ref{lem:scmle1} is given after the proof of Theorem \ref{thm:scmle}.
The proof of Lemma \ref{lem:scmle1} is given later in Appendix \ref{lem:scmle1_pro}.

\begin{lem}\label{lem:scmle2}
$L_{\z^0}(\z^0) \leq L_{\z^0}(\z)$ for all $\z\in P$, with equality holding if and only if $\z=\z^0$.
\IEEEQED
\end{lem}
The proof of Lemma \ref{lem:scmle2} follows a standard argument. A typical treatment can be found in \cite{1988Ca}.

\paragraph*{Convergence of the Functions $L\left(\hth_t^{[1:N_0]};\z\right)$}

Now $P$ is a compact set, so it is sequentially compact \cite{1955Ke}, and the sequence $(\hzn_t;\, N_0\in\Z_1)$ has a convergent subsequence $(\hz_t^{N_M};\, M\in\Z_1)$ for all $t\geq 0$, for which $\lim_{M\conv\infty}\lim_{t\conv\infty} \hz_t^{N_M}= \z^*\in P$, in the topology of $P$. Further, we observe that $\z^*$ is a $\B$-measurable $P$-valued random variable. We will adopt the notation $(\hzn;\, N_0\in\Z_1 ) \teq (\limt \hzn_t;\, N_0\in \Z_1)$ in order to denote the sequence of MLE estimates indexed by the size of the population. 

%We showed in Lemma \ref{lem:scmle2} that $L_{\z^0}(\z^0) \leq L_{\z^0}(\z)$ for all $\z\in P$, with equality holding if and only if $\z=\z^0$. Now, 
We will show that $L_{\z^0}(\z^*)\leq L_{\z^0}(\z^0)$ for any $\z^0\in P$. This, together with Lemma \ref{lem:scmle2} with $\z$ set equal to $\z^*$ gives $L_{\z^0}(\z^*) \leq L_{\z^0}(\z^0) \leq L_{\z^0}(\z^*)$. The Identifiability Condition \ass{ass:MLE_ident} gives $\z^0=\z^*$ w.p.1 and we conclude that all subsequential limits of $(\hzn;\, N_0\in \Z_1)$ equal $\z^0$ w.p.1, and hence $\limN\lim_{t\conv\infty} \hzn_t = \z^0$ w.p.1.
%\begin{equation}\label{eqn:scmle4}
%\limN\lim_{t\conv\infty} \hzn_t = \z^0 \quad \text{w.p.1}.
%\end{equation}

%(i)
\noindent (i) \textit{To show $L_{\z^0}(\z^*) \leq \limt L_{\z^0} (\hzn_t) + \e/3$:}\quad
For any $\e>0$, there exists an almost surely finite random integer $N_1(\w,\e)$ so that for all $M$ such that $N_M>N_1(\w,\e)$, the estimate $\hz^{N_M}\teq\lim_{t\conv\infty} \hz_t^{N_M}$ lies in a neighbourhood $\Ne_{\z^*}(\z)$ of $\z^*$ for which $\lvert L_{\z^0}(\z) - L_{\z^0}(\z^*) \rvert < \e/3$, for all $\z\in\Ne_{\z^*}(\z)$, by the continuity of $L_{\z^0}(\cdot)$ on $P$. The uniform continuity of $L_{\z^0}$ on $P$ is shown as follows:
% first method
% pick an arbitrary $\z_1\in P$. Because of the almost sure uniform convergence of $L_N(\z)$ to $L_{\z^0}(\z)$ on $P$ we can find an $N_c$ such that $\lvert L_{N_0}(\z) - L_{\z^0}(\z) \rvert < \e/3$ for all $N_0\geq N_c$ for all $\z\in P$, where $L_{N_0}(\z)\teq L\left(\th^{[1:N_0]};\z\right)$. As $L_{N_0}$ is uniformly continuous over $P$ for each $N_0\in\Z_1$, we can find $\delta_N > 0$ such that $\lvert L_{N_0}(\z_1) - L_{N_0}(\z)\rvert < \e/3$ for $\lvert \z_1 - \z \rvert < \delta_N$. But then we have 
%\begin{equation*}
% \lvert L_{\z^0}(\z_1) - L_{\z^0}(\z)\rvert \leq \lvert L_{\z^0}(\z_1) - L_{N_0}(\z_1) \rvert + \lvert L_{N_0}(\z_1) - L_{N_0}(\z) \rvert + \lvert L_{N_0}(\z) - L_{\z^0}(\z) \rvert \leq \e,
%\end{equation*}
%as long as $\lvert \z_1 - \z \rvert < \delta_N,\, \z \in P$, which means the function $L_{\z^0}(\z)$ is uniformly continuous on $P$.
% second method
pick arbitrary $\z,\z' \in P\subset \tilde P$ such that $\z'\in P$ lies in a coordinate neighborhood $\Ne_\delta(\z)$ of $\z\in P$. We have 
\begin{align*}
\lvert L_{\z^0}(\z)-L_{\z^0}(\z') \rvert & = \lvert -\E_{\z^0}\log f_\z(\th) + \E_{\z^0}\log f_{\z'}(\th) \rvert \\
& \leq \int_\Th \lvert \log f_{\z'}(\th) - \log f_\z(\th) \rvert f_{\z^0}(\th) d\th. 
\end{align*}
Hence for some  $\z''\in \tilde P$ in the line segment $\{\lambda \z + (1-\lambda) \z' ; \, \lambda\in (0,1)\}$, the Mean Value Theorem yields 
\begin{equation}\label{eqn:LPzeta}
\lvert L_{\z^0}(\z)-L_{\z^0}(\z') \rvert \leq \int_\Th  \frac{1}{f_{\z''}(\th)}  \lVert f'_{\z''}(\th) \rVert \lVert \z' - \z \rVert f_{\z^0}(\th) d\th.
\end{equation}
But by \ass{ass:pdf}, $f_{\z''}(\th) > 0$ for all $\z''\in \tilde P$. Then by \eqref{eqn:LPzeta}, for each $\e>0$, there exists $\delta_\e>0$ such that for all $\z,\z' \in P$, $\lVert \z' - \z \rVert<\delta_\e$ implies $\lvert L_{\z^0}(\z)-L_{\z^0}(\z') \rvert < \e$. Hence, $L_{\z^0}(\z)$ is uniformly continuous over $P$.
%  for all $\z\in P$, $\z' \in \Ne_\delta(\z)$ and $\z'' \in \tilde P$. Then for arbitrary $\delta>0$ take an open covering of $P$ by $\delta-$coordinate neighbourhoods and let $\{\Ne_\delta^i(\z);\, 1\leq i \leq M\}$ be a finite subcover. For all $\e>0$, pick the smallest $\delta>0$ among the finite $\delta-$coordinate neighbourhoods so that for all $\e>0$ there exists $\delta>0$ such that for all $\z,\z' \in P$, $\lVert \z - \z' \rVert <\delta$ implies $\lvert L_{\z^0}(\z) - L_{\z^0}(\z') \rvert < \e$. 

%(ii)
\noindent (ii) \textit{To show $\limt L_{\z^0} (\hzn_t) \leq \limt L(\hth_t^{[1:N_0]};\, \hzn_t ) + \e/3$ for all $N\geq N_2(\w,\e)$ for some random $N_2(\w,\e) \in \Z_1$:}\quad
Lemma \ref{lem:scmle1} assures us that we can pick an almost surely finite random integer $N_2(\w,\e)\in \Z_1$ so that for all $N>N_2(\w,\e)$, we have $\lvert L(\th^{[1:N_0]};\, \z)-L_{\z^0}(\z) \rvert < \e/3$ for all $\z\in P$, where $N_0 = \lvert Obs(N) \rvert,\, N_0 \conv \infty$, as $N\conv\infty$, and where $\th^{[1:N_0]} = \limt \hth_t^{[1:N_0]}$. But, from the continuity of $L_{N_0}(\cdot),\,N_0\in\Z_1$, we have $\limt L(\hth_t^{[1:N_0]};\, \z) = L(\th^{[1:N_0]};\, \z)$ for all $\z \in P$, therefore the inequality holds.
%\begin{equation}\label{eqn:scmleLn}
% \lvert L(\th^{[1:N_0]};\, \z)-L_{\z^0}(\z) \rvert < \frac{\e}{3} \quad \text{for all }  \z\in P,
%\end{equation}
%\begin{equation}
% \limt L(\hth_t^{[1:N_0]};\, \z) = L(\th^{[1:N_0]};\, \z) \quad \text{for all } \z \in P.
%\end{equation}

%(iii) 
%\begin{equation}\label{eqn:scmledef}
% L(\th^{[1:N_0]};\, \hzn) \leq L(\th^{[1:N_0]};\, \z)\quad \text{for all }  \z\in P,
%\end{equation}

\noindent (iii) \textit{To show $\limt L(\hth_t^{[1:N_0]};\, \hzn_t ) \leq \limt L(\hth_t^{[1:N_0]};\, \z ),\, \forall \z\in P $:}\quad
This follows from \eqref{eqn:nslmle} since we have $L(\th^{[1:N_0]};\, \hzn) \leq L(\th^{[1:N_0]};\, \z)$ for all $\z\in P$, where $\limt L(\hth_t^{[1:N_0]};\, \hzn_t) = L(\th^{[1:N_0]};\, \hzn)$, and $\limt L(\hth_t^{[1:N_0]};\, \z) = L(\th^{[1:N_0]};\, \z)$.% for all $\z\in P$. 

%(iv)
%\begin{equation}\label{eqn:scmleLn2}
% \lvert L(\th^{[1:N_0]};\, \z) - L_{\z^0}(\z) \rvert < \frac{\e}{3} \quad \text{for all } \z\in P.
%\end{equation}
\noindent (iv) \textit{To show $\limt L(\hth_t^{[1:N_0]};\, \z )  \leq L_{\z^0}(\z) + \e/3,\, \forall \z \in P$:}\quad
Again, we employ Lemma \ref{lem:scmle1}: pick an almost surely finite random integer $N_3(\w,\e)\in \Z_1$ so that for all $N>N_3(\w,\e)$ we have $\lvert L(\th^{[1:N_0]};\, \z) - L_{\z^0}(\z) \rvert < \e/3$ for all $\z\in P$, and let $\hth_t^{[1:N_0]} \conv \th^{[1:N_0]}$ as $t\conv\infty$.
%Combining \eqref{eqn:scmleLp}, \eqref{eqn:scmleLn}, \eqref{eqn:scmledef} and \eqref{eqn:scmleLn2} we obtain

Combining (i)-(iv), yields
\begin{align*}
L_{\z^0}(\z^*) & \leq \limt L_{\z^0} (\hzn_t) + \frac{\e}{3} \leq \limt L(\hth_t^{[1:N_0]};\, \hzn_t ) + \frac{2\e}{3}\\
& \leq \limt L(\hth_t^{[1:N_0]};\, \z ) + \frac{2\e}{3}\\
& \leq L_{\z^0}(\z) + \e \quad \text{w.p.1 for all } \z\in P,
\end{align*}
for all $N_M>\max(N_1,N_2,N_3)(\w,\e)$. 

\paragraph*{Convergence of the Sequence of Estimators $(\hz^{N_0}_t;\, N_0\in\Z_1)$}

Evaluating the relation above at $\z = \z^0$ yields $L_{\z^0}(\z^*) \leq L_{\z^0}(\z^0) + \e$ w.p.1. But this expression is independent of $N_M$, and $\e$ is arbitrary. So, $L_{\z^0}(\z^*)\leq L_{\z^0}(\z^0)$ w.p.1 for all $\z\in P$.
%\begin{equation}\label{eqn:scmleLpzs}
%  L_{\z^0}(\z^*) \leq L_{\z^0}(\z^0) + \e \quad \text{w.p.1}.
%\end{equation}
% Hence, \eqref{eqn:scmle11} and Lemma \ref{lem:scmle2} together give
However, as stated in Lemma \ref{lem:scmle2}, $L_{\z^0}(\z^0) \leq L_{\z^0}(\z)$ w.p.1 for all $\z\in P$, with equality holding if and only if $\z=\z^0$. Therefore, $L_{\z^0}(\z^*) \leq L_{\z^0}(\z^0) \leq L_{\z^0}(\z^*)$ w.p.1, implying $\z^0=\z^*$ w.p.1 by the Identifiability Condition \ass{ass:MLE_ident}, and so all subsequential limits of $(\hzn;\, N_0 \in \Z_1)=(\limt \hzn_t ;\, N_0\in\Z_1)$ equal $\z^0$ w.p.1, or equivalently $\limN \limt \hz^{N_0}_t = \z^0$ w.p.1.
\IEEEQED

\subsection*{Proof of Lemma \ref{lem:scmle1}}\label{lem:scmle1_pro}
%%%%%%%%%%%%%%%%%%%%%%%%%%%%%%%%%%%%%%%%%%%%%%%%%%%%%%%%%%%%%%
% Proposition
%%%%%%%%%%%%%%%%%%%%%%%%%%%%%%%%%%%%%%%%%%%%%%%%%%%%%%%%%%%%%%

By \ass{ass:pdf} the family of densities $f_\z \triangleq \{ f_\z(\th);\th \in \Th, \z\in P\}$ exists for the family of dynamical and cost function parameter distributions $\{ F_\z(.);\, \z \in P\}$. Let $f(\th^{[1:N_0]};\, \z)$, where $N_0 \triangleq \lvert Obs(N) \rvert$, be the likelihood function of $f_\z$ at $\th^{[1:N_0]}$ and let $L_{N_0}(\z) \teq L(\th^{[1:N_0]};\, \z)$ be the continuously differentiable function of $f(\th^{[1:N_0]} ;\, \z)$, given by the scaled log-likelihood function \[ L_{N_0}(\z) \teq L(\th^{[1:N_0]};\, \z) \teq -(1/N_0)\log f(\th^{[1:N_0]};\, \z) \equiv -(1/N_0) \log\left ( \prod_{j\in Obs(N)} f_\z(\th_j) \right), \] $N_0 \triangleq \lvert Obs(N) \rvert$, where $\th^{[1:N_0]} = \{ \th_j;\, j\in Obs(N)\}$. 

The random sequence $L(\z) \teq (L(\th^{[1:N_0]};\z)$; $N_0 \in \Z_1)$ converges w.p.1 for each $\z\in P$ \cite{1988Ca}, where $P$ is a compact set by \ass{ass:pdf}. Then, in order for the almost sure convergence of $L(\z)$ to be uniform over $P$, it is sufficient that the process $((\p L_{N_0}/\p \z)(\z);\,N_0 \in \Z_1)$ exists as a sequence of random variables which is w.p.1 bounded uniformly over $P$, where $L_{N_0}(\z)\teq L(\th^{[1:N_0]};\z)$. This may be shown as follows by the Mean Value Theorem:
\begin{equation*}
\tilde{L}_{K_0,L_0}(\z') \teq L_{K_0}(\z') - L_{L_0}(\z') = \tilde{L}_{K_0,L_0}(\z) + \frac{\p\tilde{L}_{K_0,L_0}}{\p \z} (\z'') (\z'-\z),
\end{equation*}
where $\z'\in P$ lies in an $\e-$coordinate neighbourhood $\Ne_\e(\z)$ of $\z\in P$ and $\z''$ lies on the line segment $\{\lambda \z + (1-\lambda) \z' ; \, \lambda\in (0,1)\}$. Consequently, 
\begin{equation*}
 \lvert \tilde{L}_{K_0,L_0}(\z') \rvert \leq \lvert \tilde{L}_{K_0,L_0}(\z) \rvert + \left\lVert \frac{\p\tilde{L}_{K_0,L_0}}{\p \z}(\z'') \right\rVert \lVert \z' - \z \rVert \quad \text{with} \enspace \lVert \z'-\z \rVert < \e,
\end{equation*}
where the differentiability of $(L_N;\,N\in\Z_1)$ follows from its definition. Let each such $\Ne_\e(\z) \subset \tilde{P}$, choosing a smaller $\e = \e(\z)$, possibly depending upon $\z$, if necessary. Then take an open cover of the compact set $P$ by these $\e-$neighbourhoods and let $\{\Ne_\e^i(\z);\, 1\leq i \leq M\}$ be a finite subcover. By \ass{ass:pdf} for each $j,\, 1\leq j \leq p$, $(\p f_\z/\p \z_j)(\th;\z)$ is bounded uniformly for all $\th\in\Th$. Therefore, $\sup_{\z\in \tilde{P}}( \lVert (\p L_{N_0} / \p \z)(\z) \rVert;\,  N_0 \in \Z_1)<K$. Moreover, by the convergence of the sequences $(L_{N_0}^i(\z);\, 1\leq i \leq M)$ w.p.1 and the boundedness of $( \lVert (\p L_{N_0} / \p \z)(\z) \rVert;\,  N_0 \in \Z_1)$ by $K$ uniformly over $\z\in \tilde{P}$, we obtain $\lvert \tilde{L}_{K_0,L_0}(\z') \rvert \leq \e + 2K\e$ w.p.1 for all $\z' \in P$, for all $K_0,L_0 \geq N(\w,\e)$ for some random $N(\w,\e) \in \Z_1$. But this shows that $(L_{N_0}(\z);\, N_0 \in \Z_1)$ satisfies the Cauchy convergence criterion w.p.1 uniformly over P. Therefore $L_{N_0} \teq L\left(\th^{[1:N_0]};\z\right) \conv L_{\z^0}(\z) \teq -\E_{\z^0}\log f_\z(\th)$ w.p.1, as $N_0\conv\infty$ uniformly over $P$, where $N_0 \teq \lvert Obs(N) \rvert$, and hence as $N\conv\f$ uniformly over $P$.\IEEEQED

%===============================================================================

%===============================================================================
\section{}\label{sec:app_abncee}%{Asymptotic Behaviour of the MF Equations} 
%===============================================================================
%===============================================================================
%%%%%%%%%%%%%%%%%%%%%%%%%%%%%%%%%%%%%%%%%%%%%%%%%%%%%%%%%%%%%%%%%%%%%%%%%%%%%%%%
% Proof of Proposition prop:s_conv
%%%%%%%%%%%%%%%%%%%%%%%%%%%%%%%%%%%%%%%%%%%%%%%%%%%%%%%%%%%%%%%%%%%%%%%%%%%%%%%%
\subsection*{Proof of Proposition \ref{prop:s_conv}}\label{prop:s_conv_pro}
% For the work in this and other sections we shall use the following terminology: let $\{ \mathbf{M}_i,1\leq i \leq k \}$ be a set of real $n \times n$ matrices. Then we shall term the standard relation $\left \lVert \frac{1}{k}\sum_{l=1}^k \mathbf{M}_l   \right \rVert \leq \left ( \frac{1}{k} \sum_{l=1}^k \lVert \mathbf{M}_l \rVert^2 \right )^{1/2}$, the \emph{Jensen's Inequality}.

\noindent 1) \textit{Proof of $\limN\lim_{t\conv\infty} x^*(\tau, \hzn_t) = x^* (\tau, \z^0)$ w.p.1, $t\leq \tau <\infty$:}

Recall that $x^*(\tau,\z^0),\, t\leq \tau < \infty$, is the solution to the MF Equation System \eqref{eqn:ncee}, and $x^*(\tau,\hzn_t),\, t\leq \tau<\infty$, is the solution to the MF-SAC Equation System \eqref{eqn:nceea}. 
Note that the subscript $i$ is suppressed for clarity. 
A contraction mapping argument together with \ass{ass:ind}-\ass{ass:cost_coup} ensure the existence and uniqueness of $x^*(\cdot,\z^0)\in\b{C}_b[0,\infty)$ (see \cite{2007HCM_TAC}). 
\ass{ass:ind}-\ass{ass:cost_coup} also hold for $x^*(\tau,\hzn_t),\, t\leq \tau < \infty,\, t\geq 0$, by Lemma \ref{lem:uni_con}; therefore, the existence and uniqueness properties also hold for $x^*(\tau,\hzn_t)$ for $t\geq 0$. 
Since $x^*(\cdot,\z)$ is a continuous function of $\z$ on $P$, and by Theorem \ref{thm:scmle}, $\limN\limt\hzn_t=\z^0$ w.p.1, 
$\lim_{t\conv\infty}\sup_{t\leq \tau < \infty} \left \lVert x^*(\tau,\hzn_t) - x^*(\tau,\z^0) \right \rVert = O(\e_1(N))$, where $\e_1(N)\conv 0 $ as $N\conv \infty$. 
Therefore, \[ \limN\lim_{t\conv\infty} \left \lVert x^*(\tau,\hzn_t) - x^*(\tau,\z^0) \right \rVert = 0 \, \text{w.p.1}, \, t\leq \tau <\infty.\]

\noindent 2) \textit{Proof of $\limN\lim_{t\conv\infty} s(t;\hth_t^{pr},\hzn_t) = s(t;\th^0,\z^0)$ w.p.1:} 

The solution to the differential equation \eqref{eqn:DE_s} is given by \[s(t;\th^0,\z^0) = - \int_t^\infty \b{\Psi}_{t,\tau}^{-1}(t,\tau,\th^0) \b{Q}(\th^0) x^*(\tau,\z^0) d\tau,\] where $\frac{d}{dt}\b{\Psi}_{t,t_0} = \b{A}_*(\th^0) \b{\Psi}_{t,t_0},\, \b{\Psi}_{t_0,t_0}=\b{I}$, and $x^*$ is generated by the MF equation system \eqref{eqn:ncee}, and $\mathbf{A}_* \triangleq (\mathbf{A}-\mathbf{BR^{-1}B\t\Pi})$. For the certainty equivalence offset function $s(t; \hth_t,\hz_t)$ generated by the MF-SAC Law, we have \[s(t;\hth_t^{pr},\hzn_t)=-\int_t^\infty \mathbf{\Phi}^{-1} (t,\tau,\th^0,\hth^{\tau,t}) \b{Q}(\hth_t) x^*(\tau,\hzn) d\tau,\] where $\frac{d}{dt}\b{\Phi}_{t,t_0} = \b{A}_*(\hth_t) \b{\Phi}_{t,t_0},\, \b{\Phi}_{t_0,t_0}=\b{I}$. We adopt the notation $\b{\Phi}_{t,\tau}\triangleq \b{\Phi}(t,\tau,\th^0,\hth^{\tau,t}),\, \b{\Psi}_{t,\tau} \triangleq \b{\Psi}(t,\tau,\th^0)$ and obtain
% The norm of the difference between these two terms is
\begin{equation*}
\left\lVert s(t;\hth_t,\hzn_t) - s(t;\th^0,\z^0) \right\rVert
= \Big \lVert \int_t^\infty \mathbf{\Phi}_{t,\tau}^{-1} \b{Q} (\hth_t)  x^*(\tau,\hzn_t) d\tau  -  \int_t^\infty \b{\Psi}_{t,\tau}^{-1} \b{Q}(\th^0) x^* (\tau,\z^0) d\tau \Big \rVert.
\end{equation*}

Adding and subtracting $\int_t^\infty \mathbf{\Phi}_{t,\tau}^{-1}  \b{Q}(\hth_t) x^*(\tau,\z^0) d\tau$ and $\int_t^\infty \mathbf{\Phi}_{t,\tau}^{-1} \b{Q}(\th^0) x^*(\tau,\z^0) d\tau$, and using the triangle inequality we get
\begin{equation*}
\begin{aligned}
& \left\lVert s(t;\hth_t,\hzn_t) - s(t;\th^0,\z^0) \right\rVert \leq \Big\lVert \int_t^\infty \mathbf{\Phi}_{t,\tau}^{-1} \b{Q}(\hth_t)\\
& x^*(\tau,\hzn_t) d\tau   - \int_t^\infty \mathbf{\Phi}_{t,\tau}^{-1} \b{Q}(\hth_t) x^*(\tau,\z^0) d\tau \Big\rVert + \Big\lVert \int_t^\infty \mathbf{\Phi}_{t,\tau}^{-1} \\
& \b{Q}(\hth_t) x^*(\tau,\z^0) d\tau  - \int_t^\infty \mathbf{\Phi}_{t,\tau}^{-1} \b{Q}(\th^0) x^*(\tau,\z^0) d\tau \Big\rVert  \\
& +  \Big \lVert \int_t^\infty \b{\Phi}_{t,\tau}^{-1} \b{Q}(\th^0) x^*(\tau,\z^0) d\tau   \\
& - \int_t^\infty \b{\Psi}_{t,\tau}^{-1}  \b{Q}(\th^0) x^* (\tau,\z^0)d\tau  \Big\rVert =:  I_1^{N,t} + I_2^t + I_3^t.
\end{aligned}
\end{equation*}

\noindent (i) \textit{Convergence of $I_1^{N,t}$ and $I_2^t$:}
$\lim_{t\conv \infty} I_1^{N,t} = O(\e_1(N))$, where $\e_1(N)\conv 0$, as $N\conv\infty$, and $\lim_{t\conv \infty} I_2^t = 0$ follows from Lemma \ref{lem:state_transition} and Part 1 of the proof.% Therefore, $\limN\lim_{t\conv \infty} I_1^{N,t} = 0$ w.p.1.

\noindent (ii) \textit{Convergence $I_3^t$:}
From the proof of Lemma \ref{lem:state_transition_conv},
\begin{align*}
I_3^t &= \left\lVert \int_t^\infty \left (\b{\Phi}_{t,\tau}^{-1} - \b{\Psi}_{t,\tau}^{-1} \right ) \b{Q}(\th^0)x^* (\tau,\z^0) d\tau \right\rVert\\
\label{eqn:prop_own_s_conv_1} &= \big \lVert \int_t^\infty \left [ \int_t^\tau \b{\Phi}_{\tau,s} \left ( \b{A}_* (\th^0) - \b{A}_* (\hth_s)  \right ) e^{\b{A}_* (\th^0)\t (s-t)} ds \right ] \\
& \qquad \b{Q}(\th^0) x^* (\tau,\z^0) d\tau \big\rVert,
\end{align*}
$t \leq s \leq \tau < \infty$. Lemma \ref{lem:state_transition} yields the bound $\lVert \b{\Phi}_{\tau,s}\rVert \leq \beta_0 e^{-\rho_{\Phi}(\tau-s)},\,  0 \leq s \leq \tau$, and for the time invariant case, $\lVert e^{A_*^0(s-t)}\rVert \leq \beta_1 e^{-\rho_{A_*^0}(s-t)},\, 0 \leq t \leq s$. For simplicity, set $\rho<\min[\rho_{\Phi},\rho_{A_*^0}]$ and let $T_\w(\e_2)$ be such that $\lVert \b{A}_*(\th^0) - \b{A}_*(\hth_s) \rVert < \e_2,\, s>T_\w(\e_2)$. Then for $T_\w(\e_2) < t\leq s \leq \tau < \infty$ we obtain
\begin{equation}\label{eqn:prop_own_s_conv_2}
\begin{aligned} 
I_3^t & \leq \int_t^\infty \Big[ \int_t^\tau \left\lVert \b{\Phi}_{\tau,s} \right \rVert  \left \lVert \b{A}_*(\th^0) - \b A_*(\hth_s) \right \rVert \\
& \qquad \left\lVert e^{\b{A}_*(\th^0)\t (s-t)} \right\rVert ds  \Big] \lVert \b{Q}^0 \rVert \left \lVert x^* (\tau,\z^0) \right \rVert d\tau\\
& \leq \beta_0 \beta_1 \lambda_0 \e_2 \lVert \b{Q}^0 \rVert \int_t^\infty \left( \int_t^\tau e^{-\rho(\tau-s)}   e^{-\rho(s-t)}  ds \right)  d\tau,
\end{aligned}
\end{equation}
%%
% \nonumber & \qquadsix \qquadfive T_\w(\e_2) < t\leq s \leq \tau < \infty ,\\
%%
where $\b{Q}^0\teq \b{Q}(\th^0)$. The term \eqref{eqn:prop_own_s_conv_2} is satisfied for all arbitrarily small $\e_2>0$ for all sufficiently large $t\geq T_\w(\e_2)$ by use of the bounds $\lVert \b{\Phi}_{t,s} \rVert \leq \beta_0 e^{-\rho(t-s)}$, $\lVert e^{A^0(t-s)} \rVert \leq \beta_1 e^{-\rho(t-s)}$, and $\lambda_0 = \sup_\tau \lVert x^*(\tau,\z^0) \rVert$. Hence, $I_3^t \leq \beta_0 \beta_1 \lambda_0 \e_2 \lVert \b{Q}^0 \rVert \int_t^\infty (\tau-t ) e^{-\rho(\tau-t)} d\tau = \kappa \e_2,\text{ w.p.1, where }\kappa = \beta_0\beta_1 \lambda_0 \lVert \b{Q}^0 \rVert / \rho^2$. By Theorem \ref{thm:str_con_est} $\lVert \b{A}_*(\th^0)-\b{A}_*(\hth_t) \rVert\conv 0$ as $t\conv\infty$; therefore, as $t\conv \infty$, $\e_2\conv 0$. Hence, we obtain $\lim_{t\conv \infty} I_3^t = 0$ w.p.1.

In conclusion we have shown that $\lim_{t\conv\infty} I_1^{N,t} =$ 
$O(\e_1(N))$, and therefore $\limN\lim_{t\conv\infty} I_1^{N,T}=0$ w.p.1. 
In addition, $\lim_{t\conv \infty} I_2^t = 0$ and $\lim_{t\conv \infty} I_3^t = 0$. 
Therefore, $\limN\lim_{t\conv\infty} I^{N,t} = 0$ w.p.1. Hence, 
$\lim_{t\conv\infty} \lVert s(t; \hth_t^{pr},\hzn_t) - s(t;\th^0,\z^0) \rVert = O(\e_1(N))$, and 
$\limN\lim_{t\conv\infty} \lVert s(t; \hth_t^{pr},\hzn_t) - s(t;\th^0,\z^0) \rVert = 0 \text{ w.p.1}$.

\noindent 3) \textit{Proof of ${\hat{u}}^0 (t;\hth_t^{pr},\hzn_t) = -\b{R}^{-1}\hat{\b{B}}_{t}\t (\hat{\b{\Pi}}_{t} x_t +$} \textit{$s(t;\hth_t^{pr},\hzn_t) ) + \xi_k\left [\e(t)-\e(k) \right ]$:}

The solution $\b{\Pi}_\theta \in \R^{n^2},\, \th\in\R^{n(n+m+(n+1)/2)}$, to the algebraic Riccati equation \eqref{eqn:DE_Pi} parametrized by $\th\in \Th$ is a smooth function of $\th$ (see \cite{1984De_TAC}). Hence, $(\b{\Pi}(\hth_t^{pr});t\geq 0)$ satisfies $\b{\Pi}(\hth_t^{pr})<\infty$ w.p.1 for all $t\geq 0$ since $\hth_t^{pr}\in\Th,\, t\geq 0$. It is shown in Part 1 of the proof that the mass signal $x^*(\tau,\hzn_t) \in\b{C}_b[0,\infty),t\leq \tau < \infty$; therefore, $s(t;\hth_t^{pr},\hzn_t)<\infty$ w.p.1 for all $t\geq 0$, evaluated along $\hth_t^{pr},\, t\geq 0$. Hence, $u^0 (t;\hth_t^{pr},\hzn_t)$ is well defined.
%The mass signal $x^*$ is shown in Part 1 of the proof to be bounded, i.e., $x^*\in\b{C}_b[0,\infty)$;
\IEEEQED
%===============================================================================

%===============================================================================
%%%%%%%%%%%%%%%%%%%%%%%%%%%%%%%%%%%%%%%%%%%%%%%%%%%%%%%%%%%%%%%%%%%%%%%%%%%%%%%%
% Proof of Theorem thm:pop_L2_conv
%%%%%%%%%%%%%%%%%%%%%%%%%%%%%%%%%%%%%%%%%%%%%%%%%%%%%%%%%%%%%%%%%%%%%%%%%%%%%%%%
\subsection*{Proof of Theorem \ref{thm:pop_L2_conv}}\label{thm:pop_L2_conv_pro}
We recall the following notation and basic assumptions: $\th_i^0$ denotes the true dynamical parameter of agent $A_i$ in $\Th$ that parametrizes the matrices $[\b{A}_i,\b{B}_i,\b{Q}_i] \in\Th$, which are to be estimated by agent $A_i$, and $\hth_{i,t}=[\hat{\b{A}}_{i,t},\hat{\b{B}}_{i,t},\b{Q}_i]$ is the estimated parameter of agent $A_i$ at time $t$. Note that $\b{Q}_i$ is in the information set of agent $A_i$, therefore does not need to be estimated. We set $\hth_i^{\tau,t}=(\hth_i^s,\tau \leq s \leq t )$, the sample path of the estimated parameter matrices from time $\tau$ to time $t$. The population distribution parameter denoted by $\z^0 \in P$, where $P$ is the parameter set for $F_\z(\cdot)$, parametrizes $F_\z(\cdot)$. Further, the estimated population distribution parameter of agent $A_i$ is denoted as $\hzn_{i,t}$, and $\hz_i^{\tau,t}(N_0) \teq ( \hzn_{i,s},\tau \leq s \leq t)$, is the sample path of the estimated distribution parameter from time $\tau$ to $t$. As shown in Theorem \ref{thm:str_con_est}, under \ass{ass:ind}-\ass{ass:compact}, on the probability space $\ps$, $(\hth_{i,t};\,t\geq 0)$ converges w.p.1 to $\th_i^0$ as $t\conv\infty$, and by Theorem \ref{thm:scmle}, under \ass{ass:pdf} and \ass{ass:MLE_ident}, $(\hzn_{i,t};\,t\geq 0)$ converges w.p.1 to $\z^0$ as $t\conv\infty$ and $N\conv\infty,\, 1\leq i \leq N$. Note that for the optional PCPI, \oass{ass:invB} also needs to be employed. In the sequel, $\b{A}_{\hth_{i,t}},\,\b{B}_{\hth_{i,t}}$ will be used to denote the estimated dynamical parameters whereas $\b{\Pi}_{\hth_{i,t}}$ denotes the solution to \eqref{eqn:nslr}.
%By abuse of notation we use $\th^0$ to denote $\th_i^0$,  the dynamical parameter of agent $A_i$ in $\Th$ that parametrizes the matrices $[\b{A}_i,\b{B}_i] \in\Th$, which are to be estimated by agent $A_i$; this is the self identification procedure and when no confusion can arise $\hth_t$ is used to denote $\hth_{i,t}$, the self estimated parameter of agent $A_i$ at time $t$. $\b{A}^0,\b{B}^0,\b{\Pi}^0$ denote $\b{A}_i(\th_i^0),\b{B}_i(\th_i^0),\b{\Pi}_i(\th_i^0)$ and $\hat{\b{A}}_t,\hat{\b{B}}_t,\hat{\b{\Pi}}_t$ denote $\b{A}_i(\hth_{i,t}),\b{B}_i(\hth_{i,t}),\b{\Pi}_i(\hth_{i,t})$ respectively.
% $\hat{ \th }_\tau^t$ shall denote $\hat{ \b { \th } }_i^{ \tau,t }=( \hth_i^s,\tau \leq s \leq t )$, the sample path of estimated parameter matrices from time $\tau$ to time $t$.
Since the solution $\b{\Pi}_\theta \in \R^{n^2},\, \th\in\R^{n(n+m+(n+1)/2)}$, to the algebraic Riccati equation \eqref{eqn:DE_Pi} parametrized by $\th\in\Th$, is a smooth function of $\th$ (see \cite{1984De_TAC}), $\b{\Pi}_{\hth_{i,t}},t\geq 0$, satisfies $\left \lVert \b{\Pi}_{\hth_{i,t}}-\b{\Pi}_{\th_i^0} \right\rVert \conv 0$ w.p.1 as $t\conv \infty,\, 1\leq i \leq N$. To establish the theorem we first observe $x_i^0(t;\hth_i^{0,t},\hz_i^{0,t}(N_0)),\, t\geq 0,\, 1\leq i \leq N$, is the state of the system subject to the dithered MF Adaptive control law computed from the sum
\begin{equation}\label{eqn:control_def2}
u_i^0(t;\hth_{i,t},\hzn_{i,t})=u_i^{loc}(t;\hth_{i,t})+u_i^{pop}(t;\hth_{i,t},\hzn_{i,t})+u_i^{dit}(t),\quad t\geq 0,\quad 1\leq i \leq N,
\end{equation}
where the control input due to the MF-SAC Law is given by
\[
u_i^0(t;\hth_{i,t},\hzn_{i,t}) = - \b{R}^{-1}\b{B}_{\th_{i,t}}\t  \b{\Pi}_{\th_{i,t}}
x_i(t) -\b{R}^{-1}\b{B}_{\th_{i,t}}\t  s(t;\hth_{i,t},\hzn_{i,t}) + \xi_k\left [\e_i(t)-\e_i(k) \right],\quad t\geq 0.
\]
The term $\left\lVert  x_i^0(t;\hth_{i}^{0,t},\hz_{i}^{0,t}(N_0)) -  x_i^0(t;\th_i^0,\z^0) \right\rVert^2$ will be decomposed into four parts, and convergence properties will be established for each term. We have \[ d\b{\Psi}_i(t,\tau,\th_i^0)=[\b{A}_{\th_i^0} -\b{B}_{\th_i^0} \b{R}^{-1} \b{B}_{\th_i^0}\t \b{\Pi}_{\th_i^0}] \b{\Psi}_i(t,\tau,\th_i^0)dt,\; \b{\Psi}_i(\tau,\tau,\th^0)=\b{I},\] and \[d\b{\Phi}_i(t,\tau,\th_i^0,\hth_i^{\tau,t})=[\b{A}_{\th_i^0} -\b{B}_{\th_i^0} \b{R}^{-1} \b{B}_{\hat\th_{i,t}}\t \b{\Pi}_{\hat\th_{i,t}}] \b{\Phi}_i(t,\tau,\th_i^0,\hth_i^{\tau,t})dt,\; \b{\Phi}_i(\tau,\tau,\th^0,\hth_{i,\tau})=\b{I}.\] Also, in the sequel for clarity we will suppress the subscript $i$ and adopt the notation: $\b{\Phi}_{t,s}\triangleq\b{\Phi}(t,s,\th^0,\hth^{0,t}),\, \b{\Psi}_{t,s}\triangleq\b{\Psi}(t,s,\th^0),\, \b{A}^0 \triangleq \b{A}_{\th^0}$, $\b{B}^0 \triangleq \b{B}_{\th^0},\, \b{\Pi}^0 \triangleq \b{\Pi}_{\th^0},\,s(t)\teq s(t;\th_i^0,\z^0),\, x^0(t)\teq x_i^0(t;\th_i^0,\z^0),\, \hat{\b{A}}_t \triangleq \b{A}_{\hth_t} ,\, \hat{\b{B}}_t \triangleq \b{B}_{\hth_t},\, \hat{\b{\Pi}}_t \triangleq \b{\Pi}_{\hth_t},\, \hat{s}(t)\teq s(t;\hth_{i,t},\hzn_{i,t}),\, \hat{x}_i^0(t)\teq x_i^0(t;\hth_i^{0,t},\hz_i^{0,t}(N_0))$. Displaying the dependence of the fundamental matrix on the parameter estimate trajectory, we use the integral representation and by use of the Cauchy Schwarz Inequality (henceforth termed CS)'s Inequality we obtain
\begin{align}
& \label{eqn:own_stateconv1} \quad \frac{1}{T}\int_0^T \lVert x^0(t;\hth^{0,t},\hz^{0,t}(N_0)) - x^0(t;\th^0,\z^0)\rVert^2 dt  \leq \frac{4}{T}\int_0^T \Big \lVert\b{\Phi}_{t,0} x(0) - \b{\Psi}_{t,0} x(0)\Big \rVert^2 dt\\
& \label{eqn:own_stateconv2}  + \frac{4}{T}\int_0^T \Big \lVert \int_0^t \b{\Phi}_{t,\tau} \b{B}^0 \b{R}^{-1} \hat{\b{B}}\t \hat{s}(t) d\tau - \int_0^t \b{\Psi}_{t,\tau}\b{B}^0 \b{R}^{-1}  {\b{B}^0}\t s(t) d\tau  \Big \rVert^2 dt\\
& \label{eqn:own_stateconv3}  + \frac{4}{T}\int_0^T \Big\lVert \b{\Phi}_{t,0} \int_0^t \b{\Phi}^{-1}_{\tau,0} \b{D} dw(\tau) - \b{\Psi}_{t,0} \int_0^t \b{\Psi}_{\tau,0}^{-1} \b{D} dw(\tau) \Big  \rVert^2 dt\\
& \label{eqn:own_stateconv4}  + \frac{4}{T} \int_0^T \Bigg \lVert \sum _{k=0}^{\lfloor t \rfloor} \Big ( \int_k^{\min[t,k+1]} \b{\Phi}_{t,\tau}\hat{\b{B}}_t  \xi_k [\e(\tau)-\e(k) ] d\tau \Big )   - \sum _{k=0}^{\lfloor t \rfloor} \Big ( \int_k^{\min[t,k+1]}  \b{\Psi}_{t,\tau} \b{B}^0  \xi_k [\e(\tau)-\e(k) ] d\tau \Big ) \Bigg \rVert^2 dt\\
& \nonumber \quad =: I_1^T + I_2^{N,T} + I_3^{T} + I_4^T.
\end{align}
% ACK-LV
% \begin{align}
% %%
% & \nonumber \frac{1}{T}\int_0^T \lVert x_i^0(t;\hth_i^{0,t},\hz_i^{0,t}(N)) - x_i^0(t;\th_i^0,\z^0)\rVert^2 dt \\
% %%
% & \label{eqn:own_stateconv1} \leq \frac{4}{T}\int_0^T \Big \lVert\b{\Phi}_i(t,0,\th_i^0,\hth^{0,t}) x_i(0) - \b{\Psi}_i(t,0,\th_i^0) x_i(0)\Big \rVert^2 dt\\
% %%
% & \qquad \label{eqn:own_stateconv2}+ \frac{4}{T}\int_0^T \Big \lVert \int_0^t \b{\Phi}_i(t,\tau,\th_i^0,\hth_i^{\tau,t}) \b{B}_{\th_i^0} \b{R}^{-1} \b{B}_{\hth_{i,\tau}}\t s(\tau;\hth_{i,\tau},\hzn_{i,\tau}) d\tau - \int_0^t \b{\Psi}_i(t,\tau,\th_i^0)\b{B}_{\th_i^0} \b{R}^{-1}  \b{B}_{\th_i^0}\t s(\tau;\th_i^0,\z^0) d\tau  \Big \rVert^2 dt\\
% %%
% & \qquad \label{eqn:own_stateconv3}+ \frac{4}{T}\int_0^T \Big\lVert \b{\Phi}_i(t,0,\th_i^0,\hth_i^{\tau,t}) \int_0^t \b{\Phi}_i^{-1}(\tau,0,\th_i^0,\hth_i^{\tau,t}) \b{D} dw_i(\tau) - \b{\Psi}_i(t,0,\th_i^0) \int_0^t \b{\Psi}_i^{-1}(\tau,0,\th_i^0) \b{D} dw_i(\tau) \Big  \rVert^2 dt\\
% %%
% & \nonumber \qquad + \frac{4}{T} \int_0^T \Bigg \lVert \sum _{k=0}^{\lfloor t \rfloor} \Big ( \int_k^{\min[t,k+1]} \b{\Phi}_i(t,\tau,\th_i^0,\hth_i^{\tau,t})\b{B}_{\hth_{i,\tau}}  \xi_k [\e_i(\tau)-\e_i(k) ] d\tau \Big )  \\
% %%
% & \qquadtwo \label{eqn:own_stateconv4} - \sum _{k=0}^{\lfloor t \rfloor} \Big ( \int_k^{\min[t,k+1]}  \b{\Psi}(t,\tau,\th_i^0) \b{B}_{\th_i^0}  \xi_k [\e_i(\tau)-\e_i(k) ] d\tau \Big ) \Bigg \rVert^2 dt\\
% %%
% & \nonumber \leq I_1^T + I_2^{N,T} + I_3^{T} + I_4^T.
% \end{align}
% ACK-LV
We will show one by one that the limit supremums $(\lsT)$ of $I_1^T$ \eqref{eqn:own_stateconv1}, $I_3^T$ \eqref{eqn:own_stateconv3} and $I_4^T$ \eqref{eqn:own_stateconv4} are all 0 with probability 1, and the limit supremum $(\lsN\lsT)$ of $I_2^{N,T}$ \eqref{eqn:own_stateconv2} is equal to 0 with probability 1.
\begin{enumerate}[(i)]
%%%%%%%%%%%%%%%%%%%%%%%%%%%%%%%%%%%%%%%%%%%%%%%%%%%%%%%%%%%%%%%%%%
% first submission
% \item \textit{Convergence of $I_1^T$ \eqref{eqn:own_stateconv1}:} We have $I_1^T =  \frac{4}{T}\int_0^T \left \lVert \b{\Phi}(t,0,\th^0,\hth^{0,t}) x(0) - \b{\Psi}(t,0,\th^0) x(0) \right \rVert^2 dt \leq \frac{4 \Vert x(0) \rVert^2}{T}$\\$\times \int_0^T \left \lVert\b{\Phi}_{t,0} - \b{\Psi}_{t,0}  \right \rVert^2 dt$ $=: 4 \Vert x(0)\rVert^2 I_{11}^T$. We have $\b{\Psi}(t,0,\th^0) = \exp(\b{A}_{*}(\th^0)t)$, $\b{A}_{*}(\th^0) \triangleq \b{A}^0-\b{B}^0 \b{R}^{-1}{\b{B}^0}\t\b{\Pi}^0$. As $\b{A}_{*}(\th^0)$ is asymptotically stable due to \ass{ass:conobs}, and $\hth_t \conv \th^0$ w.p.1 as $t\conv \infty$, due to Theorem \ref{thm:str_con_est}, $\lsT I_{11}^T=0$, by Lemma \ref{lem:state_transition_conv}. Therefore, employing \ass{ass:ind}, $\lsT I_1^T = 4 \lVert x(0) \rVert^2 \lsT I_{11}^T = 0$.
% first submission
% second submission
\item \textit{Convergence of $I_1^T$} follows from Lemma \ref{lem:state_transition_conv}.
% second submission
\item \textit{Convergence of $I_2^{N,T}$ \eqref{eqn:own_stateconv2}:} Adding and subtracting $\b{\Phi}(t,\tau,\th^0,\hth^{\tau,t}) \b{B}^0 \b{R}^{-1}  \hat{\b{B}}\t s(\tau;\th^0,\z^0) $ and\\ $\b{\Phi}(t,\tau,\th^0,\hth^{\tau,t}) \b{B}^0 \b{R}^{-1}  {\b{B}^0}\t s(\tau;\th^0,\z^0) $ 
% second submission
using Lemma \ref{lem:state_transition}, Lemma \ref{lem:state_transition_conv}, and \[ \sup_{\tau\geq t}\lVert s(\tau;\hth_\tau,\hzn_\tau) - s(\tau;\th^0,\z^0) \rVert <\e_1(N),\] from Proposition \ref{prop:s_conv} we get $\lsT I_2^{N,T} \leq  O((\e_1(N))^2)$ w.p.1, which implies \[\lsN \lsT I_2^{N,T} = 0 \quad \text{w.p.1}.\]
\item \textit{Convergence of $I_3^T$ \eqref{eqn:own_stateconv3}:}
We have
\begin{equation}\label{eqn:I_3^T}
I_3^T = \frac{4}{T}\int_0^T \left\lVert \b{\Phi}(t,0,\th^0,\hth^{\tau,t}) \int_0^t \b{\Phi}^{-1}(\tau,0,\th^0,\hth^{\tau,t}) \b{D} dw(\tau) - \b{\Psi}(t,0,\th^0) \int_0^t \b{\Psi}^{-1}(\tau,0,\th^0) \b{D} dw(\tau)  \right  \rVert^2 dt,
\end{equation}
% first submission
% where $\frac{d}{dt} \b{\Psi}_{t,t_0}=\mathbf{A}_*(\th^0)  \b{\Psi}_{t,t_0},\; \frac{d}{dt} \b{\Phi}_{t,t_0}=\mathbf{A}_*(\hth_t)  \b{\Phi}_{t,t_0}$ with $\b{\Psi}_{t_0,t_0}=\b{I},\, \b{\Phi}_{t_0,t_0}=\b{I}$. 
% first submission
We use the notation $\b{\Phi}_{t,\tau} \triangleq \b{\Phi}(t,\tau,\th^0,\hth^{\tau,t})$, $\b{\Psi}_{t,\tau} \triangleq \b{\Psi}(t,\tau,\th^0)$. Consider the stochastic differential equations
\begin{align*}
dx_t = \b{A}_*(\hth_t) x_t dt + \b{D} dw(t),\; x(0) = x_0,\\
dy_t = \b{A}_*(\th^0) y_t dt + \b{D} dw(t),\; y(0) = y_0,
\end{align*}
where $x_0 = y_0 < \infty$ by \ass{ass:ind}.
The difference $z_t = x_t - y_t$, satisfies
\[
z_t = \b{\Phi}_{t,0} \int_0^t \b{\Phi}^{-1}_{\tau,0} \b{D} dw(\tau) - \b{\Psi}_{t,0} \int_0^t \b{\Psi}^{-1}_{\tau,0} \b{D} dw(\tau),\quad t\geq 0. 
\]
Alternatively, one can write $dz_t = \b{A}_*(\th^0) z_t dt + [\b{A}_*(\th^0) - \b{A}_*(\hth_{t})]y_t dt$, giving $z_t = \b{\Psi}_{t,t_0} z_0 + \int_0^t \b{\Psi}_{t,s} ( \b{A}_*(\th^0) -  \b{A}_*(\hth_{s})) y_s ds$. Hence we can write \eqref{eqn:I_3^T} as $I_3^T = \frac{4}{T}\int_0^T \Big\lVert \b{\Psi}_{t,t_0} z_0 + \int_0^t \b{\Psi}_{t,s} ( \b{A}_*(\th^0) - \b{A}_*(\hth_{s})) y_s ds  \Big  \rVert^2 dt$, and use the CS Inequality to obtain
\[ 
I_3^T \leq \; \frac{8}{T}\int_0^T  \Big\lVert \b{\Psi}_{t,t_0} z_0 \Big\rVert^2 dt + \frac{8}{T}\int_0^T \Big\lVert \int_0^t \b{\Psi}_{t,s} ( \b{A}_*(\th^0) -  \b{A}_*(\hth_{s})) y_s ds   \Big  \rVert^2  dt.
\]
Now $z_0 = x_0 - y_0 = 0$, since $x_0 = y_0$; therefore, $I_3^T \leq \;  \frac{8}{T}\int_0^T \Big\lVert \int_0^t \b{\Psi}_{t,s} ( \b{A}_*(\th^0) -  \b{A}_*(\hth_{s})) y_s ds   \Big  \rVert^2  dt$.
% Following Lemma \ref{lem:state_transition} we write $\lVert \b{\Psi}_{t,t_0} \rVert \leq \beta_0 e^{-\rho_{A^0}(t-t_0)}$. Also adopting Hypothesis \ass{ass:ind} we get
% \be
% I_{31}^T \leq \frac{8}{T}\int_0^T \big\lVert \beta_0 e^{-\rho_{A^0}(t-t_0)} \big\rVert \big\lVert z_0 \big\rVert dt
% \ee
Let $T_\w(\e_2)$ be such that $t\geq T_\w(\e_2)$ implies $\lVert \hth_t - \th^0 \rVert < \e_2$. Then, $I_{3}^T \leq \frac{8}{T}\int_0^{T_\w} \Big\lVert \cdot  \Big\rVert^2 dt + \frac{8}{T}\int_{T_\w}^T \Big\lVert \cdot  \Big\rVert^2 dt =:   I_{31}^T + I_{32}^T$. 
% first submission
% Showing $\lsT I_{31}^T = 0$ is a special case of the proof of $\lsT I_{32}^T = 0$, and is therefore omitted. 
% \paragraph*{For $I_{32}^T$} $I_{32}^T = \frac{8}{T}\int_{T_\w}^T \Big\lVert \int_0^t \b{\Psi}_{t,s} ( \b{A}_*(\th^0) -  \b{A}_*(\hth_{s})) y_s ds     \Big\rVert^2 dt,\; 0\leq s \leq t,\; T_\w \leq t \leq T$. The inner integral can be separated as $I_{32}^T \leq \frac{16}{T}\int_{T_\w}^T \Big\lVert \int_0^{T_\w} ( \cdot ) ds  \Big\rVert^2 dt +  \frac{16}{T}\int_{T_\w}^T \Big\lVert  \int_{T_\w}^t ( \cdot ) ds  \Big\rVert^2 dt :=  16 I_{321}^T + 16 I_{322}^T$. We will only show $\lsT I_{322}^T = 0$ w.p.1 here, as the solution methodology is the same for both. For $I_{322}^T$, we have $I_{322}^T \; =  \;  \frac{1}{T}\int_{T_\w}^T \Big\lVert  \int_{T_\w}^t \b{\Psi}_{t,s} ( \b{A}_*(\th^0) -  \b{A}_*(\hth_{s})) y_s ds   \Big\rVert^2 dt,\; T_\w \leq \tau \leq t,\; T_\w \leq t \leq T$. Following Lemma \ref{lem:state_transition} we write $\lVert \b{\Psi}_{t,s} \rVert \leq \beta_0 e^{-\rho_{A^0}(t-s)}$ for $t\geq s \geq 0$. We also use the CS Inequality and obtain $I_{322}^T \; =  \;  \frac{1}{T}\int_{T_\w}^T (t-T_\w)   \left ( \int_{T_\w}^t \beta_0^2 e^{-2\rho_{A^0}(t-t_0)} \e_2^2 y_s^2 ds  \right)  dt,\; T_\w \leq \tau \leq t,\; T_\w \leq t \leq T$. $\lsT I_{322}^T = 0$ w.p.1 can be easily shown. This implies $\lsT I_{31}^T = 0$ and $\lsT I_{32}^T = 0$ w.p.1. Therefore, $\lsT I_3^T = 0$ w.p.1. 
% first submission
% second submission
Following Lemma \ref{lem:state_transition} we write $\lVert \b{\Psi}_{t,s} \rVert \leq \beta_0 e^{-\rho_{A^0}(t-s)}$ for $t\geq s \geq 0$. We also use the CS Inequality, and let $\e_2\conv 0$ as $t\conv \f$. We get $\lsT I_{31}^T = 0$ and $\lsT I_{32}^T = 0$ w.p.1. Therefore $\lsT I_3^T = 0$ w.p.1.
% second submission
%%%%%%%%%%%%%%%%%%%%%%%%%%%%%%%% I_4^T %%%%%%%%%%%%%%%%%%%%%%%%%%%%%%%%%%%%%%%%%
\item \textit{Convergence of $I_4^T$ \eqref{eqn:own_stateconv4}:}
We have
%ACK-LV
%\begin{align*}
%& \lsT\frac{1}{T} \int_0^T \Bigg \lVert \sum _{k=0}^{\lfloor t \rfloor} \left ( \int_k^{\min[t,k+1]} \b{\Phi}(t,\tau,\th^0,\hth^{\tau,t})\hat{\b{B}}_\tau \xi_k\left [\e(\tau)-\e(k) \right ] d\tau \right )\\
%& \qquad \qquad -  \sum _{k=0}^{\lfloor t \rfloor} \left ( \int_k^{\min[t,k+1]} \b{\Psi}(t,\tau,\th^0) \b{B}^0  \xi_k\left [\e(\tau)-\e(k) \right ] d\tau \right ) \Bigg \rVert^2 dt,
%\end{align*}
%ACK-LV
\begin{align*}
\lsT\frac{1}{T} \int_0^T \left \lVert  \sum _{k=0}^{\lfloor t \rfloor} \left ( \int_k^{\min[t,k+1]} \left ( \b{\Phi}_{t,\tau}- \b{\Psi}_{t,\tau} \right ) \b{B}^0 \xi_k\left [\e(\tau)-\e(k) \right ] d\tau \right ) \right \rVert^2 dt.
\end{align*}
This term is treated by a direct application of Lemma \ref{lem:kronecker}; therefore, the limit of the time average integral tends to 0.
\end{enumerate}
Overall, we have shown that $\lsT I_1^T = 0$, $\lsT I_2^{N,T} = O((\e_1(N))^2)$, $\lsT I_3^T = 0$ and $\lsT I_4^T = 0$. This implies $\limsup_{T\conv \infty} (1/T)\int_0^T\lVert  x^0(t;\hth^{0,t},\hz^{0,t}(N_0)) - x^0(t;\th^0,\z^0) \rVert^2 = O((\e_1(N))^2)$ w.p.1. Consequently, $\limN\limT (1/T) \int_0^T\lVert  x_i^0(t;\hth_i^{0,t},\hz_i^{0,t}(N_0)) - x_i^0(t;\th_i^0,\z^0) \rVert^2 = 0$ w.p.1, $1\leq i \leq N$.\IEEEQED

%===============================================================================
%%%%%%%%%%%%%%%%%%%%%%%%%%%%%%%%%%%%%%%%%%%%%%%%%%%%%%%%%%%%%%
% Proposition
%%%%%%%%%%%%%%%%%%%%%%%%%%%%%%%%%%%%%%%%%%%%%%%%%%%%%%%%%%%%%%
\begin{prop}\label{prop:pop_u_conv}
For the system \eqref{eqn:dynamics}, let \ass{ass:ind}-\ass{ass:cost_coup}, \ass{ass:pdf}, \ass{ass:MLE_ident} hold. Let $\huio \in \hat{\U}_{MF}^N$ be the MF-SAC input \eqref{eqn:nslu} and $u_i^0 \in \U_{MF}^N $ be the non-adaptive MF-SC input. Then,%, and let $\hxio\conv\xio$ as $t\conv\infty$ and $N\conv\infty,\, 1\leq i \leq N$. Then,
%and $\hth_{i,t} \conv \th_i^0$ w.p.1 as $t\conv\infty$; $\hzn_{i,t}\conv \z_i^0$ w.p.1 as $t\conv \infty,\, N\conv \infty$
\begin{equation*}
\lsN\lsT  \frac{1}{T} \int_0^T \lVert  \huio - u_i^0 \rVert^2 dt = 0\quad \text{w.p.1}, \quad 1\leq i \leq N.
\end{equation*}%\IEEEQED
\end{prop}
%The result is proved in Appendix \ref{prop:pop_u_conv_pro}.
%%%%%%%%%%%%%%%%%%%%%%%%%%%%%%%%%%%%%%%%%%%%%%%%%%%%%%%%%%%%%%%%%%%%%%%%%%%%%%%%
% Proof of Lemma prop:u_conv
%%%%%%%%%%%%%%%%%%%%%%%%%%%%%%%%%%%%%%%%%%%%%%%%%%%%%%%%%%%%%%%%%%%%%%%%%%%%%%%%
\begin{IEEEproof}
We have the term $I^{N,T} =  \frac{1}{T}\int_0^T\lVert u_i^0 (t;\hth_{i,t},\hzn_{i,t}) - u_i^0(t;\th_i^0,\z^0)  \rVert^2 dt$,
which we separate into two parts as $I^{N,T} =  \frac{1}{T}\int_0^{T_\w} \lVert \cdot \rVert^2 dt +  \frac{1}{T}\int_{T_\w}^{T} \lVert \cdot \rVert^2 dt =: I_1^{N,T} + I_2^{N,T}$, where $T_\w$ is a random instant to be determined later. We will only establish $\limN \limT I_2^{N,T} = 0$ w.p.1 here, as $\limN \limT I_1^{N,T} = 0$ w.p.1 is a simpler case of the same argument.
%%%%%%%%%%%%%%%%%%%%% I_2 %%%%%%%%%%%%%%%%%%%
\paragraph*{Convergence of $I_2^{N,T}$}
We have
\begin{multline*}
I_2^{N,T} = \frac{1}{T}\int_{T_\w}^T \lVert u_i^0 (t;\hth_{i,t},\hzn_{i,t}) - u_i^0(t;\th_i^0,\z^0) \rVert^2 dt  = \\
\frac{1}{T}\int_{T_\w}^T \lVert \b{R}^{-1}\hat{\b{B}}\t\hat{\b{\Pi}}_t x_i^0(t;\hth_i^{0,t},\hz_i^{0,t}(N_0)) + \b{R}^{-1}\hat{\b{B}}_t\t s(t;\hth_{i,t},\hzn_{i,t}) \\-  \b{R}^{-1}{\b{B}^0}\t\b{\Pi}^0 x_i^0(t;\th_i^0,\z^0) - \b{R}^{-1}{\b{B}^0}\t s(t;\th_i^0,\z^0) \rVert^2 dt,
\end{multline*}
Dropping the subscript $i$, adopting the notation $\b{B}^0\teq \b{B}(\th^0),\, \hat{\b{B}}_t\teq \b{B}(\hth_t),\, \b{\Pi}^0\teq \b{\Pi}(\th^0),\, \hat{\b{\Pi}}_t\teq \b{\Pi}(\hth_t),\,\uo\triangleq u^0(t;\th^0,\z^0),\, \huo\triangleq \uo(t;\hth_t,\hzn_t),\, \xo\triangleq \xo(t;\th^0,\z^0),\, \hxo \triangleq \xo(t;\th^{0,t},\hz^{0,t}(N_0)),\, s(t) \teq s(t;\th^0,\z^0),\, \hat{s}(t)\teq s(t;\hth_t,\hzn_t)$, and using the CS Inequality, we obtain
\begin{align*}
I_2^{N,T} & \leq  \frac{2}{T}\int_{T_\w}^T  \lVert \b{R}^{-1}\hat{\b{B}}_t\t \hat{\b{\Pi}}_t \hxo(t) - \b{R}^{-1} {\b{B}^0}\t {\b{\Pi}^0} \xo(t) \rVert^2 dt \\ & \qquadtwo + \frac{2}{T}\int_{T_\w}^T \lVert \b{R}^{-1} \hat{\b B}_t\t s(t;\hth_{t},\hzn_{t}) - \b{R}^{-1}{\b{B}^0}\t s(t;\th^0,\z^0)  \rVert^2 dt\\
& =: I_{21}^{N,T} + I_{22}^{N,T}.
\end{align*}

We set $T_\w$ to be the random instant such that $t\geq T_\w$ implies $\lVert \hxo(t)-\xo(t) \rVert < \e_1(N)$ and $\lVert \hat{s}(t)-s(t)\rVert < \e_1(N)$. We obtain $\lsT I_{21}^{N,T} = O((\e_1(N))^2)$ and $\lsT I_{22}^{N,T} = O((\e_1(N))^2)$ from Section \ref{prop:s_conv_pro} 2.(i), which implies 
% remove for now
% $\lsN\lsT I_2^{N,T} \leq \lsN\lsT I_{21}^{N,T} + \lsN\lsT I_{22}^{N,T} = 0$ w.p.1. With $\lsN\lsT I_1^{N,T} = 0$ w.p.1, 
% remove for now
\begin{align*}
\lsN\lsT I^{N,T} & \leq \lsN \lsT I_1^{N,T} + \lsN \lsT I_2^{N,T} \\  & = 0 \text{ w.p.1.}
\end{align*}
\end{IEEEproof}

%===============================================================================
\section{}\label{sec:app_MT}%{Proof of the Main Theorem} 
%===============================================================================

The following five lemmas will be used to prove Proposition \ref{prop:eq_costs} and Proposition \ref{prop:eq_costs_g}. We use the notation $m(\xno(t;\th^{[1:N]},\z^0))\triangleq m((1/N)\sum_{k=1}^N x_k^0(t;\th_i^0,\z^0)+\eta),\, m(\xno(t;\hth^{[1:N]},\hz^{[1:N]}))\triangleq m((1/N)\sum_{k=1}^N x_k^0(t;\hth_{i,t},\hzn_{i,t})+\eta)$, where $m(\cdot)$ is defined in \ass{ass:cost_coup}.
%%%%%%%%%%%%%%%%%%%%%%%%%%%%%%%%%%%%%%%%%%%%%%%%%%%%%%%%%%%%%%
% Lemma 1
%%%%%%%%%%%%%%%%%%%%%%%%%%%%%%%%%%%%%%%%%%%%%%%%%%%%%%%%%%%%%%
\begin{lem}\label{lem:stability}
Let Assumptions \ass{ass:ind}-\ass{ass:cost_coup} hold. For the system \eqref{eqn:dynamics}, the MF control law $u_i(t;\th_i^0,\z^0)$ \eqref{eqn:MF_control} and its corresponding closed-loop solution $x_i^0(t;\th_i^0,\z^0)$ satisfy
\begin{equation}\label{eqn:lem_sta_main}
\sup_{N\geq 1}\max_{1\leq i \leq N}\limsup_{T\conv \infty}\frac{1}{T}\int_0^T \left ( \lVert  x_i^0(t;\th_i^0,\z^0) \rVert^2 + \lVert u_i^0(t;\th_i^0,\z^0) \rVert^2 \right ) dt < \infty.
\end{equation}
\end{lem}
\begin{IEEEproof}

%===============================================================================
%%%%%%%%%%%%%%%%%%%%%%%%%%%%%%%%%%%%%%%%%%%%%%%%%%%%%%%%%%%%%%%%%%%%%%%%%%%%%%%%
% Proof of Lemma lem:stability
%%%%%%%%%%%%%%%%%%%%%%%%%%%%%%%%%%%%%%%%%%%%%%%%%%%%%%%%%%%%%%%%%%%%%%%%%%%%%%%%
%ACK-LV
%\subsection{Proof of Lemma \ref{lem:stability}}\label{lem:stability_pro}

 The same result has been shown to hold in \cite{2008LZ_TAC} (Theorem 4.1) for control action in the form of $u_i^0(t) = u_i^{loc}(t) + u_i^{pop}(t)$ using the notation defined in \eqref{eqn:control_def2}. We are going to repeat this result here for completeness.

 %%%%%%%%%%%%%%%%%% (i) %%%%%%%%%%%%%%%%%
 (i) $\lsT \frac{1}{T} \int_0^T \lVert x_i^0(t;\th_i^0,\z^0) \rVert^2 dt \leq K_1 < \infty$:
 \begin{align}
  \lVert x_i^0(t;\th_i^0,\z^0) \rVert^2 = & \Bigg \lVert e^{\b{A}_{*}(\th_i^0)t}x_i(0) - \int_0^t e^{\b{A}_{*}(\th_i^0)(t-\tau)} \b{B}_i^0 \b{R}^{-1}  {\b{B}_i^0}\t s(\tau;\th_i^0,\z^0) d\tau + \\
 & \qquad \int_0^t e^{\b{A}_{*}(\th_i^0) (t-\tau)} \b{D} dw_i(\tau) \Bigg \rVert^2.
 \end{align}

 Using the CS Inequality, we obtain the inequality:
 \begin{align}
 \lsT \frac{1}{T} \int_0^T \lVert x_i^0(t;\th_i^0,\z^0) \rVert^2 dt & \leq \lsT \frac{3}{T} \int_0^T \Big \lVert e^{\b{A}_{*}(\th_i^0)t}x_i(0) \Big \rVert^2 dt +\\  
 & \lsT \frac{3}{T} \int_0^T \Big \lVert \int_0^t e^{\b{A}_{*}(\th_i^0)(t-\tau)} \b{B}_i^0 \b{R}^{-1}  {\b{B}_i^0}\t s(\tau;\th_i^0,\z^0) d\tau \Big \rVert^2 dt + \\
  & \lsT \frac{3}{T} \int_0^T \Big \lVert \int_0^t e^{\b{A}_{*}(\th_i^0) (t-\tau)} \b{D} dw_i(\tau) \Big \rVert^2 dt,\\
  & \leq \lsT I_1^T + \lsT I_2^T + \lsT I_3^T.
  \end{align}

  For \ass{ass:ind} and \ass{ass:conobs} hold, we obtain $\lsT I_1^T = 0$ w.p.1. 
  
  It is shown in \cite{2008LZ_TAC} (Theorem 4.1) that $\lsT I_2^T \leq \kappa_1 < \infty$ uniformly for all $\th^0\in\Th$.
  
  Using Lemma \ref{lem:ergodicity} we write 
  \be
  \lsT I_3^T = 3 \int_0^\infty \tr \left ( e^{\b{A}_{*}(\b{\th}_i^0) (t-\tau)} \b{D} \b{D}\t e^{\b{A}_{*}(\b{\th}_i^0) (t-\tau)} \right ).
  \ee
  
  We use $\sup_{\th^0 \in \Th}\lVert e^{\b{A}_{*}(\b{\th}_i^0) (t-\tau)} \rVert \leq \beta e^{-\rho (t-s)}$ as shown in Lemma \ref{lem:state_transition} and get 
  \be
  \lsT I_3^T \leq 3 \lVert \b{D} \rVert^2 \beta^2 / 2\rho = \kappa_2 < \infty.
  \ee

  Therefore, 
  \be\label{eqn:sta_1}
  \sup_{N\geq 1}\max_{1\leq i \leq N}\lsT \frac{1}{T} \int_0^T \lVert x_i^0(t;\th_i^0,\z^0) \rVert^2 dt \leq \kappa_1 + \kappa_2 = K_1 < \infty \quad \text{w.p.1}.
  \ee
  
 %%%%%%%%%%%%%%%%%% (ii) %%%%%%%%%%%%%%%%%
 (ii) $\lsT \frac{1}{T} \int_0^T \lVert u_i^0(t;\th_i^0,\z^0)  \rVert^2 dt \leq K_2 < \infty$:
 We have the MF Control Law
  \begin{align}
  u_i^0(t;\th_i^0,\z^0) & =  u_i^{loc}(t;\th_i^0)+u_i^{pop}(t;\th_i^0,\z^0)\\
  & = - \b{R}^{-1} {\b{B}_i^0}\t \left( \b{\Pi}_i^0 x_i(t) +  s(t;\th_i^0,\z^0)\right ),\quad t\geq 0.
  \end{align}
  
 Also, the mass offset function is
 \be s (t;\th_i^0,\z^0) = - e^{-\b{A}_{*}\t (\th_i^0)  t} \int_t^\infty e^{\b{A}_{*}\t (\th_i^0)  \tau } \b{Q}_i x^*(\tau,\z^0) d\tau.\ee
  
 We employ \ass{ass:ind} and obtain $M_{x^*} = \sup_{\tau \geq t} \lVert x^* \rVert$, $M_B = \sup_{\th \in \Th} \lVert \b{B}_{\th}\rVert$, $M_{\b{\Pi}} = \sup_{\th \in \Th} \lVert \b{\Pi}_{\b{\th}}\rVert$, and $M_{\b{Q}} = \sup_{\th \in \Th} \lVert \b{Q}_{\th}\rVert$. Then we obtain
  \be
  \sup_{\th\in\Th}\lVert s_i \rVert  \leq \lVert M_{\b{Q}} \rVert M_{x^*} \beta/\rho \triangleq M_s,\enspace 1\leq i \leq N.
  \ee
  
 Using \eqref{eqn:sta_1}, and the bounds given above, we write
 \be
 \sup_{N\geq 1} \max_{1\leq i \leq N} \lsT{1/T}\int_0^T \lVert u_i^0(t;\th_i^0,\z^0) \rVert^2 dt \leq \lVert \b{R}^{-1} \rVert M_B M_\Pi K_1 +  \lVert \b{R}^{-1} \rVert M_B M_s = K_2 \quad \text{w.p.1}.
  \ee

 Consequently, we get $\sup_{N\geq 1}\max_{1\leq i \leq N}\limsup_{T\conv \infty}\frac{1}{T}\int_0^T \left ( \lVert  x_i^0(t;\th_i^0,\z^0) \rVert^2 + \lVert u_i^0(t,\th_i^0) \rVert^2 \right ) dt \leq K_1+K_2 <\infty$. As both $K_1$ and $K_2$ are independent of $1\leq i \leq N$ and $N\geq 1$, we obtain \eqref{eqn:lem_sta_main}.

\end{IEEEproof}

%%%%%%%%%%%%%%%%%%%%%%%%%%%%%%%%%%%%%%%%%%%%%%%%%%%%%%%%%%%%%%
% Lemma 2
%%%%%%%%%%%%%%%%%%%%%%%%%%%%%%%%%%%%%%%%%%%%%%%%%%%%%%%%%%%%%%
\begin{lem}\label{lem:con_x_m} 
Let Assumptions \ass{ass:ind}-\ass{ass:cost_coup} hold. For the system \eqref{eqn:dynamics}, the closed loop solution $x_i^0(t;\th_i^0,\z^0)$ with the control law $u_i^0(t;\th_i^0,\z^0)$ and the cost-coupling function $m (\xno(t;\th^{[1:N]},\z^0 ) )$ satisfy
\[%\label{eqn:lem_con_x_m}
\sup_{N\geq 1}\max_{1\leq i \leq N} \lsT  \frac{1}{T}\int_0^T \left\lVert    x_i^0(t;\th_i^0,\z^0) - m (\xno(t;\th^{[1:N]},\z^0 ) )    \right\rVert^2 dt < \infty.
\]
\end{lem}
We recall the definition $J_i(u_i,x^*) \triangleq \lsT\frac{1}{T}\int_0^T
\{ \lVert x_i - x^* \rVert_{Q_i}^2 + \lVert u_i \rVert_{R}^2 \}dt$ w.p.1, where $x^*\in\b{C}_b[0,\infty)$ is the solution to the MF Equation System \eqref{eqn:ncee}.
\begin{IEEEproof}
%ACK-LV
%===============================================================================
%%%%%%%%%%%%%%%%%%%%%%%%%%%%%%%%%%%%%%%%%%%%%%%%%%%%%%%%%%%%%%%%%%%%%%%%%%%%%%%%
% Proof of Lemma lem:con_x_m
%%%%%%%%%%%%%%%%%%%%%%%%%%%%%%%%%%%%%%%%%%%%%%%%%%%%%%%%%%%%%%%%%%%%%%%%%%%%%%%%
% ACK-LV
% \subsection{Proof of Lemma \ref{lem:con_x_m}}\label{lem:con_x_m_pro}
 
 Using the CS Inequality we write
 \begin{align}
 \label{eqn:lcxm_1} \frac{1}{T} \int_0^T \left\lVert    x_i^0(t;\th_i^0,\z^0) - m (\xno(t;\th^{[1:N]},\z^0 ) ) \right\rVert^2 dt  & \leq\\ 
 & \frac{2}{T}\int_0^T \left\lVert     x_i^0(t;\th_i^0,\z^0) \right\rVert^2 dt + \frac{2}{T}\int_0^T  \left\lVert m (\xno(t;\th^{[1:N]},\z^0 ) ) \right\rVert^2  dt\\
  & \leq I_1^T + I_2^{N,T}.
 \end{align}
  
 Using Lemma \ref{lem:stability} we get $\lsT I_1^T \leq 2 K_1$, where $K_1$ is given in \eqref{eqn:sta_1}.
  
 For $I_2^{N,T}$ we employ \ass{ass:cost_coup}, and LHS of \eqref{eqn:lcxm_1} can be further bounded by
  \begin{align}
 I_2^{N,T} = & \frac{2}{T}\int_0^T \left\lVert  m\left ( \frac{1}{N}\sum_{k=1}^N x_k^0(t;\th_k^0,\z^0) + \eta \right)\right\rVert^2 dt\\
 & \leq \frac{2\gamma^2}{T}\int_0^T \left\lVert  \frac{1}{N}\sum_{k=1}^N x_k^0(t;\th_k^0,\z^0) + \eta \right\rVert^2 dt.
  \end{align}
 Using the CS Inequality we write
  \begin{align}
 I_2^T & \leq \frac{4\gamma^2}{TN^2}\int_0^T \left\lVert  \sum_{k=1}^N x_k^0(t;\th_k^0,\z^0) \right\rVert^2 dt +  \frac{4\gamma^2}{T}\int_0^T \eta^2 dt\\
  & \leq I_{21}^{N,T} + I_{22}^T.
  \end{align}
 We have $I_{22}^T = 4\gamma^2 \eta^2$. For $I_{21}^{N,T}$ using the CS Inequality again we get
  \be
 I_{21}^{N,T} \leq  \frac{4\gamma^2 N}{TN^2}\int_0^T \lVert   x_i^0(t;\th_i^0,\z^0) \rVert^2dt.
  \ee
 We have shown in Lemma \ref{lem:stability} that $\sup_{N\geq 1}\max_{1\leq i \leq N}\lsT\frac{1}{T}\int_0^T \lVert x_i^0(t;\th_i^0,\z^0) \rVert^2 dt \leq K_1$. Therefore we get the bound
  \be
  I_2^{N,T} \leq \frac{4\gamma^2  K_1}{N} + 4\gamma^2 \eta^2.
  \ee
 We have shown that $\lsT I_1^T \leq 2K_1$. Now we have shown that $\lsT I_2^{N,T} \leq 4\gamma^2  K_1 / N + 4\gamma^2 \eta^2$. Finally we define $K_3 \triangleq 2K_1 + 4\gamma^2  K_1/N + 4\gamma^2 \eta^2$ and finish the proof:
 \be
 \sup_{N\geq 1}\max_{1\leq i \leq N} \lsT  \frac{1}{T}\int_0^T \left\lVert    x_i^0(t;\th_i^0,\z^0) - m (\xno(t;\th^{[1:N]},\z^0 ) )    \right\rVert^2 dt \leq K_3 < \infty.
 \ee%\IEEEQED
\end{IEEEproof}

%%%%%%%%%%%%%%%%%%%%%%%%%%%%%%%%%%%%%%%%%%%%%%%%%%%%%%%%%%%%%%
% Lemma 3
%%%%%%%%%%%%%%%%%%%%%%%%%%%%%%%%%%%%%%%%%%%%%%%%%%%%%%%%%%%%%%
\begin{lem}\label{lem:lim_Ji_opt_N} For the system \eqref{eqn:dynamics} subject to \ass{ass:ind}-\ass{ass:weakdist}, when all agents apply the control generated by \eqref{eqn:MF_control}, the cost function $J_i(u_i^0,u_{-i}^0)$ \eqref{eqn:cost} satisfies 
\[
\limN J_i^N(\uio,\umio) = J_i(\uio,x^*) \quad \text{w.p.1}, \quad 1\leq i \leq N.
\]
\end{lem}
\begin{IEEEproof}
%ACK-LV
%===============================================================================
%%%%%%%%%%%%%%%%%%%%%%%%%%%%%%%%%%%%%%%%%%%%%%%%%%%%%%%%%%%%%%%%%%%%%%%%%%%%%%%%
% Proof of Lemma lem:lim_Ji_opt_N
%%%%%%%%%%%%%%%%%%%%%%%%%%%%%%%%%%%%%%%%%%%%%%%%%%%%%%%%%%%%%%%%%%%%%%%%%%%%%%%%
%ACK-LV
% \subsection{Proof of Lemma \ref{lem:lim_Ji_opt_N}}\label{lem:lim_Ji_opt_N_proof}
 
 From \eqref{eqn:cost} we have the cost function
 \be
 J_i^N(u_i^0,u_{-i}^0) = \lsT\frac{1}{T}\int_0^T \left\{ \lVert x_i^0(t;\th_i^0,\z^0) - m(\xno(t;\th^{[1:N]},\z^0))  \rVert_Q^2 + \lVert u_i^0(t;\th_i^0,\z^0) \rVert_R^2 \right \} dt.
 \ee
 Adding and subtracting $x^*(t,\z^0),\, 0 \leq t \leq T$, to the first integrand on the RHS, we get
 \begin{align}
 J_i^N(u_i^0,u_{-i}^0) \leq  \enspace & J_i^\infty(u_i^0,x^*) \\ 
 &  + \lsT\frac{2}{T}\int_0^T \left \{ \left ( x_i^0(t;\th_i^0,\z^0) -x^*(t,\z^0) \right )\t \b{Q} \left ( x^*(t,\z^0) -  m(\xno(t;\th^{[1:N]},\z^0)) \right ) \right \} dt\\ 
 &   + \lsT\frac{1}{T}\int_0^T \lVert m(\xno(t;\th^{[1:N]},\z^0)) - x^*(t,\z^0) \rVert_Q^2  dt,\\
 =: \enspace & I_1 + I_2^N + I_3^N,
 \end{align}
 where, 
 \be
 J_i^\infty(u_i^0,x^*(t,\z^0))\triangleq \lsT \int_0^T \left \{ \lVert x_i^0(t;\th_i^0,\z^0) - x^*(t,\z^0) \rVert_Q^2 + \lVert u_i^0(t;\th_i^0,\z^0) \rVert_R^2 \right \} dt.
 \ee
 It is shown in \cite[Lemma 6.3]{2008LZ_TAC} that $I_2^N = O(\e_2(N))$ and $I_3^N = o(\e_2(N))$ where $\e_2(N) \conv 0$ as $N\conv \infty$. Therefore, 
 \be
 J_i^N(u_i^0,u_{-i}^0) \leq  J_i(u_i^0,x^*) + O(\e_2(N)).
 \ee
 Adding and subtracting $m(\xno(t;\th^{[1:N]},\z^0))$ to $J_i(u_i^0,x^*)$, and following the same steps above one obtains
 \be
 J_i(u_i^0,x^*) \leq  J_i^N(u_i^0,u_{-i}^0) + O(\e_2(N)).
 \ee
 Hence, one gets
 \be
 \limN J_i^N(u_i^0,u_{-i}^0) = J_i(u_i^0,x^*) \quad \text{w.p.1},\quad 1\leq i \leq N.
 \ee%\IEEEQED
\end{IEEEproof}
%%%%%%%%%%%%%%%%%%%%%%%%%%%%%%%%%%%%%%%%%%%%%%%%%%%%%%%%%%%%%%
% Lemma 4
%%%%%%%%%%%%%%%%%%%%%%%%%%%%%%%%%%%%%%%%%%%%%%%%%%%%%%%%%%%%%%
\begin{lem}\label{lem:lim_Ji_N} Under \ass{ass:ind}-\ass{ass:weakdist}, the set of controls $\U_{MF}^N = \{u_i^0,\, 1\leq i\leq N\}$ is such that when $u_i \in \U_g^N$ is any control adapted to $\F^N$, 
\[ \limN \inf_{u_i \in \U_g^N} J_i^N(u_i,u_{-i}^0) = \limN J_i^N(u_i^0,u_{-i}^0) \quad \text{w.p.1},\quad 1\leq i \leq N.
\]
\end{lem}
\begin{IEEEproof}
%ACK-LV
%===============================================================================
%%%%%%%%%%%%%%%%%%%%%%%%%%%%%%%%%%%%%%%%%%%%%%%%%%%%%%%%%%%%%%%%%%%%%%%%%%%%%%%%
% Proof of Lemma lem:lim_Ji_N
%%%%%%%%%%%%%%%%%%%%%%%%%%%%%%%%%%%%%%%%%%%%%%%%%%%%%%%%%%%%%%%%%%%%%%%%%%%%%%%%
%ACK-LV
% \subsection{Proof of Lemma \ref{lem:lim_Ji_N}}\label{lem:lim_Ji_N_proof}
 
 Let $u_i\triangleq u_i(t;\th_i^0,\z^0)\in\U_g^N$ be a feedback control action and $x_i\triangleq x_i(t;\th_i^0,\z^0)$ be the corresponding closed loop solution. Then, from \eqref{eqn:cost} we have the cost function 
 \be
 J_i^N(u_i,x^*) = \lsT\frac{1}{T}\int_0^T \left\{
 \lVert x_i(t) -x^*(t,\z^0) \rVert_Q^2 + \lVert u_i(t) \rVert_R^2 \right \} dt.
 \ee
 Adding and subtracting $\mnu(t)\triangleq m(x_i(t;\th_i^0,\z^0),\, x_{j\neq i}^0(t;\th^{[1:N]},\z^0))$ we get
 \begin{align}
 J_i^N(u_i,x^*) =  & \lsT\frac{1}{T}\int_0^T \left\{ \lVert x_i(t) - \mnu(t) + \mnu(t) - x^*(t,\z^0) \rVert_Q^2 + \lVert u_i(t) \rVert_R^2 \right \} dt\\
 & \leq \lsT\frac{1}{T}\int_0^T \{ \lVert x_i(t) - \mnu(t)  \rVert_Q^2  +  \lVert  u_i(t)  \rVert_R^2 \} dt \\ 
 & \qquad +\lsT\frac{1}{T}\int_0^T \lVert  x^*(t,\z^0) - \mnu(t) \rVert_Q^2 dt \\ 
 & \qquad   + \lsT\frac{2}{T}\int_0^T \left ( x_i(t) - \mnu(t) \right)\t \b{Q} \left ( \mnu(t) - x^*(t,\z^0)   \right ) dt\\
 & \leq  J_i^N(u_i,u_{-i}^0) + \lsT\frac{1}{T}\int_0^T \lVert  x^*(t,\z^0) - \mnu(t)  \rVert_Q^2 dt\\ 
 & \qquad +  \lsT\frac{2}{T}\int_0^T \left ( x_i(t) - \mnu(t) \right)\t \b{Q} \left ( \mnu(t) - x^*(t,\z^0) \right ) dt\\
 =: \label{eqn:lim_Ji_N_11} & J_i^N(u_i,u_{-i}^0) + I_1^N + I_2^N.
 \end{align}
 
 It is shown in \cite[Lemma 6.3]{2008LZ_TAC} that $I_1^N = o(\e_2(N))$ where $\e_2(N) \conv \infty$ as $N\conv \infty$.
 
 \paragraph*{For $I_2^N$}:
 
 We add and subtract $\mno \triangleq m(\xno(t;\th^{[1:N]},\z^0))$ and obtain
 \begin{align}
 I_2^N = & \lsT\frac{1}{T}\int_0^T \left ( x_i(t) - \mnu(t) \right)\t \b{Q} \left ( \mnu(t) - \mno(t)+ \mno(t) - x^*(t,\z^0) \right ) dt\\
 & \leq \lsT\frac{1}{T}\int_0^T \left ( x_i(t) - \mnu(t) \right)\t \b{Q} \left (\mnu(t) - \mno(t) \right ) dt\\
 & \qquad + \lsT\frac{1}{T}\int_0^T \left ( x_i(t) - \mnu(t) \right)\t \b{Q} \left ( \mno(t) - x^*(t,\z^0) \right ) dt \\
 = & I_{21}^N + I_{22}^N.
 \end{align}
 It is shown in \cite[Lemma 6.4]{2008LZ_TAC} that $\lvert I_{21}^N \rvert = O(\e_2(N))$ and it is shown in \cite[Lemma 6.4]{2008LZ_TAC} that $\lvert I_{22}^N \rvert = O(1/N)$.
 
 As $u_i^0(t;\th_i^0,\z^0)$ is the optimal tracking solution to tracking signal $x^*(t,\z^0)$ \eqref{eqn:opt_tracking}, we obtain
 \be
 J_i(u_i^0,x^*) \leq \inf_{u_i \in \U_g^N} J_i^N(u_i,\umio) + o(\e_2(N)) + O(\e_2(N)) + O(1/N).
 \ee
 
 Adding and subtracting $x^*(t,\z^0)$ to $J_i^N(u_i,u_{-i}^0)$, and following the same steps above one obtains
 \be
 \inf_{u_i \in \U_g^N}J_i^N(u_i,u_{-i}^0)    \leq  J_i(u_i^0,x^*)  + o(\e_2(N)) + O(\e_2(N)) + O(1/N).
 \ee
 
 Hence, $\limN \inf_{u_i \in \U_g^N}  J_i^N(u_i,u_{-i}^0) =  J_i(\uio,x^*)$ w.p.1, $1\leq i \leq N$. It is shown in Lemma \ref{lem:lim_Ji_opt_N} that $\limN J_i^N(u_i^0,u_{-i}^0) = J_i(u_i^0,x^*)$ w.p.1, $1\leq i \leq N$. Therefore, one gets
 \be 
 \limN \inf_{u_i \in \U_g^N}  J_i^N(u_i,u_{-i}^0) = \limN J_i^N(\uio,\umio) \quad \text{w.p.1},\quad 1\leq i \leq N.
 \ee %\IEEEQED
\end{IEEEproof}
%%%%%%%%%%%%%%%%%%%%%%%%%%%%%%%%%%%%%%%%%%%%%%%%%%%%%%%%%%%%%%
% Lemma 5
%%%%%%%%%%%%%%%%%%%%%%%%%%%%%%%%%%%%%%%%%%%%%%%%%%%%%%%%%%%%%%
\begin{lem}\label{lem:mass_conv}
Under the MF-SAC Law and \ass{ass:ind}-\ass{ass:cost_coup}, \ass{ass:pdf}, \ass{ass:MLE_ident}
\begin{equation*}
\limN\lsT \frac{1}{T} \int_0^T \lVert m(\xno(t;\th^{[1:N]},\z^0))  - m(\xno(t;\hth^{1:N},\hz^{[1:N]})) \rVert^2 dt = 0\enspace \text{w.p.1}, \quad 1\leq i \leq N.
\end{equation*}
\end{lem}

\begin{IEEEproof}
%ACK-LV
%===============================================================================
%%%%%%%%%%%%%%%%%%%%%%%%%%%%%%%%%%%%%%%%%%%%%%%%%%%%%%%%%%%%%%%%%%%%%%%%%%%%%%%%
% Proof of Lemma lem:mass_conv
%%%%%%%%%%%%%%%%%%%%%%%%%%%%%%%%%%%%%%%%%%%%%%%%%%%%%%%%%%%%%%%%%%%%%%%%%%%%%%%%
% first submission
% \subsection{Proof of Lemma \ref{lem:mass_conv}}\label{lem:mass_conv_pro}
 We have the equation $I^{N,T} = \frac{1}{T} \int_0^T \lVert m(\xno(t;\th^{[1:N]},\z^0 ) )  - m(\xno( t;\hth^{[1:N]},\hz^{[1:N]} ) )  \rVert^2 dt$,
 where $ \th^{[1:N]} \triangleq \{ \th_i^0,\, 1\leq i \leq N \},\,  \hth^{[1:N]} \triangleq \{ \hth_{i}^{\tau,t},\, 0 \leq \tau \leq T,\, 1\leq i \leq N \}$, and $\hz^{[1:N]}\triangleq \{ \hz^{\tau,t}_i(N),\, 0 \leq \tau \leq T,\, 1\leq i \leq N \}$.
 Employing \ass{ass:cost_coup}, we get the inequality $I^{N,T} \leq \frac{\gamma^2}{T} \int_0^T \left \lVert \frac{1}{N}\sum_{i=1}^N x_i^0(t;\th_i^0,\z^0)  - \frac{1}{N}\sum_{i=1}^N x_i^0(t;\hth_i^{0,t},\hz_i^{0,t}(N)) \right\rVert^2 dt$.
 Using the CS Inequality we get $I^{N,T} \leq \frac{\gamma^2}{T} \int_0^T \frac{1}{N^2} \Big \{ N \sum_{i=1}^N \lVert x_i^0(t) - \hat{x}_i^0(t)\rVert^2 \Big \} dt$, where we use the notation $x_i^0(t) = x_i^0(t;\th_i^0,\z^0)$ and $\hat{x}_i(t) = x_i(t;\hth_i^{0,t},\hz_i^{0,t}(N))$. Applying the supremum limit we get $\lsT I^{N,T} \leq \frac{\gamma^2}{N}\sum_{i=1}^N \Big \{ \lsT \frac{1}{T} \int_0^T \lVert x_i^0(t) - \hat{x}_i^0(t) \rVert^2 dt \Big \}$.
 It is shown in Theorem \ref{thm:pop_L2_conv} that $\lsT \frac{1}{T} \int_0^T \lVert x_i^0(t) - \hat{x}_i^0(t) \rVert^2 dt = O(\e_1(N)^2)$ w.p.1; hence, we get $\lsT I^{N,T} = O(\e_1(N)^2)$ w.p.1, which implies $\lsN$ $\lsT\frac{1}{T}\int_0^T\lVert m (\xno(t;\th^{[1:N]},\z^0 ) ) - m ( \xno(t;\hth^{[1:N]},\hz^{[1:N]} ) ) \rVert^2 dt = 0$ w.p.1.
% \IEEEQED
% first submission
\end{IEEEproof}

%===============================================================================
%%%%%%%%%%%%%%%%%%%%%%%%%%%%%%%%%%%%%%%%%%%%%%%%%%%%%%%%%%%%%%%%%%%%%%%%%%%%%%%%
% Proof of Proposition prop:eq_costs
%%%%%%%%%%%%%%%%%%%%%%%%%%%%%%%%%%%%%%%%%%%%%%%%%%%%%%%%%%%%%%%%%%%%%%%%%%%%%%%%
\subsection*{Proof of Proposition \ref{prop:eq_costs}}\label{prop:eq_costs_pro}
The cost function \eqref{eqn:eq_costs} is repeated here:
\begin{equation*}
J_i^N(\huio,\humio) = \lsT\frac{1}{T}\int_0^T
\big\{ \lVert x_i^0(t;\hth_i^{0,t},\hz_{i}^{0,t}(N_0)) \\ -  m ( \xno(t;\hth^{[1:N]},\hz^{[1:N]} ) ) \rVert_Q^2 + \lVert u_i(t;\hth_{i,t},\hzn_{i,t}) \rVert_R^2 \big\} dt,
\end{equation*}

where $ \hth^{[1:N]} \triangleq \{ \hth_{i}^{0,t},\, 1\leq i \leq N \}$ and $\hz^{[1:N]} \triangleq \{ \hz^{0,t}_i(N_0),\, 1\leq i \leq N \}$.
We expand the term as
\begin{multline*} 
I^{N,T} = \frac{1}{T} \int_0^T \Big\{ \lVert x_i^0(t;\hth_i^{0,t},\hz_{i}^{0,t}(N_0)) - x_i^0(t;\th_i^0,\z^0) \\ + x_i^0(t;\th_i^0,\z^0) - m (\xno(t;\th^{[1:N]},\z^0 ) ) \\
+ m(\xno(t;\th^{[1:N]},\z^0 ) ) - m ( \xno(t;\hth^{[1:N]},\hz^{[1:N]} ) ) \rVert_Q^2 \\ + \lVert u_i^0(t;\hth_{i,t},\hzn_{i,t}) - u_i^0(t;\th_i^0,\z^0) + u_i^0(t;\th_i^0,\z^0)\rVert_R^2 \Big \} dt.
\end{multline*}

In the sequel we adopt the notation $\xio \teq \xio(t;\th_i^0,\z^0)$, $\hxio \teq \xio(t;\hth_i^{0,t},\hz_{i}^{0,t}(N_0))$, $\mno \teq m(\xno(t;\th^{[1:N]},\z^0))$,  $\hmno \triangleq m(\xno(t;\hth^{[1:N]},\hz^{[1:N]}))$,  $\uio \triangleq \uio(t;\th_i^0,\z^0)$,  $\huio \triangleq \uio(t;\th_i^0,\z^0)$, and get the inequality 
\begin{equation}\label{eqn:eq_costs_d}
\begin{aligned}
J_i^N(\huio,\humio) & \leq \lsT\frac{1}{T} \int_0^T \lVert \hxio - \xio\rVert^2_Q dt\\ 
& \nqquad\nqquad + \lsT\frac{1}{T} \int_0^T\lVert \xio - \mno  \rVert^2_Q dt \\ 
& \nqquad\nqquad  + \lsT\frac{1}{T}\int_0^T\lVert \mno  - \hmno \rVert^2_Q dt \\
& \nqquad\nqquad  + \lsT\frac{2\lVert \b{Q}_i \rVert}{T}\int_0^T \left (  \hxio - \xio \right )\t \left ( \xio - \mno \right ) dt\\
& \nqquad\nqquad + \lsT\frac{2\lVert \b{Q}_i \rVert}{T}\int_0^T \left (  \hxio - \xio \right )\t \left (\mno - \hmno \right ) dt\\ 
& \nqquad\nqquad + \lsT\frac{2\lVert \b{Q}_i \rVert}{T}\int_0^T \left (  \xio - \mno \right )\t  \left ( \mno - \hmno \right ) dt\\
& \nqquad\nqquad + \lsT\frac{1}{T}\int_0^T\lVert \huio - \uio  \rVert_R^2 dt+ \lsT\frac{1}{T}\int_0^T\lVert \uio \rVert_R^2dt \\
& \nqquad\nqquad + \lsT\frac{2\lVert \b{R} \rVert}{T}\int_0^T \left ( \huio - \uio  \right )\t \left( \uio \right) dt\\
& \leq \lsT I_1^{N,T} + \lsT I_2^{T} + \lsT I_3^{N,T}\\
& \nqquad\nqquad + \lsT I_4^{N,T} + \lsT I_5^{N,T} + \lsT I_6^{N,T}\\ 
& \nqquad\nqquad + \lsT I_7^{N,T} + \lsT I_8^T + \lsT I_9^{N,T}.
\end{aligned}
\end{equation}

%\begin{enumerate}[(i)]
%%%%%%%%%%%%%%%%%%%%%%%%%%%%%%%%%%%%%%%%%%%%%%%
\noindent (i) \textit{Convergence of $I_1^{N,T}$:} We show in Theorem \ref{thm:str_con_est} that $\hth_i(t) \conv \th_i^0$ w.p.1 as $t\conv \infty$, $1\leq i \leq N$, and in Theorem \ref{thm:scmle} that $\hzn_{i,t}\conv\z^0$ w.p.1, as $t\conv\infty$ and $N\conv\infty,\, 1\leq i \leq N$; therefore the hypotheses for Theorem \ref{thm:pop_L2_conv} are satisfied and $\lsT I_1^{N,T}= O(\e_1(N)^2)$ w.p.1.

%%%%%%%%%%%%%%%%%%%%%%%%%%%%%%%%%%%%%%%%%%%%%%%
\noindent (ii) \textit{$I_2^T$ \& $I_8^T$:} $I_2^T+I_8^T$ equals to the the non-adaptive MF cost function; i.e., $J_i(u_i^0,u_{-i}^0) = \limsup_{T\conv\infty}(I_2^T+I_8^T)$ w.p.1.

%%%%%%%%%%%%%%%%%%%%%%%%%%%%%%%%%%%%%%%%%%%%%%%
\noindent (iii) \textit{Convergence of $I_3^{N,T}$:} We have 
\[ I_3^{N,T}  \leq \frac{\lVert \b{Q}_i \rVert}{T} \int_0^T \lVert \mno(t) - \hmno(t) \rVert^2 dt =: \lVert\b{Q}_i\rVert I_{31}^{N,T}.\] 
From Lemma \ref{lem:mass_conv} we have $\limsup_{T\conv\infty}  I_{31}^{N,T} = O(\e_1(N)^2)$. Therefore, 
\[ \limsup_{T\conv\infty}  I_3^{N,T} = \limsup_{T\conv\infty}  \lVert \b{Q}_i \rVert I_{31}^{N,T} = O(\e_1(N)^2)\quad \text{w.p.1}. \]

%%%%%%%%%%%%%%%%%%%%%%%%%%%%%%%%%%%%%%%%%%%%%%%
\noindent (iv) \textit{Convergence of $I_4^{N,T}$}: We have 
\begin{equation*}
I_4^{N,T} = \frac{2\lVert \b{Q}_i \rVert}{T}\int_0^T \left (  \hxio(t) - \xio(t) \right )\t 
\left ( \xio(t) - \mno(t) \right ) dt .
\end{equation*} 
Applying the CS Inequality we obtain 
\begin{align*} 
I_4^{N,T} & \leq 2 \lVert \b{Q}_i \rVert \left ( \frac{1}{T}\int_0^T \left \lVert \hxio(t) - \xio(t) \right \rVert^2 dt \right )^{1/2} \\
& \qquad \left ( \frac{1}{T}\int_0^T \left\lVert   ( \xio(t) - \mno(t)    \right\rVert^2 dt \right)^{1/2} \\
& =:2\lVert \b{Q}_i \rVert I_{41}^{N,T}\times I_{42}^T.
\end{align*}

We prove in Lemma \ref{lem:con_x_m} that $\limsup_{T\conv\infty} I_{42}^T \leq K_3$\
w.p.1. It is proved in Theorem \ref{thm:str_con_est} that $\hth_i(t) \conv \th_i^0$
w.p.1 as $t\conv\infty$ and $\hzn_{i,t}\conv\z^0$ w.p.1 as $t\conv\infty$
and $N\conv\infty$. Hence, we get $\limsup_{T\conv\infty} I_{41}^{N,T} = O(\e_1(N))$ 
w.p.1. Therefore, 
\begin{align*}
\limsup_{T\conv \infty} I_4^{N,T} &\leq 2\lVert \b{Q}_i \rVert ( \limsup_{T\conv\infty} I_{41}^{N,T} )   ( \limsup_{T\conv\infty} I_{42}^T)\\ 
& = O(\e_1(N)).
\end{align*}
Hence, $\limsup_{T\conv\infty} I_4^{N,T} = O(\e_1(N))$.

%%%%%%%%%%%%%%%%%%%%%%%%%%%%%%%%%%%%%%%%%%%%%%%
\noindent (v) \textit{Convergence of $I_5^{N,T}$:} We have the equation 
\begin{equation*} 
I_5^{N,T} = \frac{2\lVert \b{Q}_i \rVert}{T}\int_0^T\left (  \hxio(t) - \xio(t) \right )\t  \left ( \mno(t) - \hmno(t) \right ) dt .\end{equation*}
Applying the CS Inequality we obtain 
\begin{align*} 
I_5^{N,T} & \leq 2\lVert \b{Q}_i \rVert \left ( \frac{1}{T}\int_0^T \left \lVert \hxio(t) - \xio(t) \right \rVert^2 dt \right )^{1/2} \\ & \qquad \left ( \frac{1}{T}\int_0^T \left\lVert \mno(t) - \hmno(t)    \right\rVert^2 dt \right)^{1/2} \\ 
& = : 2\lVert \b{Q}_i \rVert I_{51}^{N,T}\times I_{52}^{N,T}.
\end{align*} 
We have shown in Theorem \ref{thm:str_con_est} that $\hth_i(t) \conv \th_i^0$ w.p.1 as $t\conv\infty$, and $\hzn_{i,t}\conv\z^0$ as $t\conv\infty$ and $N\conv\infty$ w.p.1. Hence, we get $\limsup_{T\conv\infty} I_{51}^{N,T} = O(\e_1(N))$ w.p.1. The convergence of $I_{52}^{N,T}$ was shown as $\limsup_{T\conv \infty} I_{52}^{N,T} = O(\e_1(N))$ w.p.1 in Lemma \ref{lem:mass_conv}. Therefore, \begin{align*} \limsup_{T\conv \infty} I_5^{N,T} &\leq 2\lVert \b{Q}_i \rVert ( \limsup_{T\conv\infty} I_{51}^{N,T} )   ( \limsup_{T\conv\infty} I_{52}^{N,T}) \\ &= O(\e_1(N)^2). \end{align*} Hence, $\limsup_{T\conv\infty} I_5^{N,T} = O(\e_1(N)^2)$.

%%%%%%%%%%%%%%%%%%%%%%%%%%%%%%%%%%%%%%%%%%%%%%%
\noindent (vi) \textit{Convergence of $I_6^{N,T}$:} We have the equation 
\begin{equation*} 
I_6^{N,T} = \frac{2\lVert \b{Q}_i \rVert}{T}\int_0^T \left (  \xio(t) - \mno(t) \right )\t  \left ( \mno(t) - \hmno(t) \right ) dt.
\end{equation*} 
Applying the CS Inequality we obtain 
\begin{align*} 
I_6^{N,T} & \leq 2\lVert \b{Q}_i \rVert \left ( \frac{1}{T}\int_0^T \left \lVert \xio(t) - \mno(t) \right \rVert^2 dt \right )^{1/2} \\
& \qquad \left ( \frac{1}{T}\int_0^T \left\lVert \mno(t) - \hmno(t) \right\rVert^2 dt \right)^{1/2} \\
& =: 2\lVert \b{Q}_i \rVert I_{61}^T\times I_{62}^{N,T}.
\end{align*} 
Using Lemma \ref{lem:con_x_m}, we get $\limsup_{T\conv\infty} I_{61}^T \leq K_3$ w.p.1. The convergence of $I_{62}^{N,T}$ was shown as 
\[ \limsup_{T\conv \infty} I_{62}^T = O(\e_1(N)) \quad \text{w.p.1} \] 
in Lemma \ref{lem:mass_conv}. Therefore, 
\begin{align*} 
\limsup_{T\conv \infty} I_6^{N,T} &\leq 2\lVert \b{Q}_i \rVert ( \limsup_{T\conv\infty} I_{61}^T ) \times  ( \limsup_{T\conv\infty} I_{62}^{N,T})\\ 
& = O(\e_1(N)) \quad w.p.1.
\end{align*} 
Hence, $\limsup_{T\conv\infty} I_6^{N,T} = O(\e_1(N))$ w.p.1.

%%%%%%%%%%%%%%%%%%%%%%%%%%%%%%%%%%%%%%%%%%%%%%%
\noindent (vii) \textit{Convergence of $I_7^{N,T}$:} We can bound $I_7^{N,T}$ from above as 
\[ I_7^{N,T} \leq \frac{\lVert \b{R} \rVert}{T} \int_0^T \lVert \huio(t) - \uio(t) \rVert^2 dt =: \lVert \b{R} \rVert I_{71}^{N,T}.\] From Proposition \ref{prop:pop_u_conv} we get $\limsup_{T\conv \infty}
I_{71}^{N,T} = O(\e_1(N)^2)$ w.p.1. Therefore, \[\limsup_{T\conv \infty} I_7^{N,T} = \limsup_{T\conv \infty}\lVert \b{R} \rVert I_{71}^{N,T} = O(\e_1(N)^2) \text{ w.p.1}.\]

%%%%%%%%%%%%%%%%%%%%%%%%%%%%%%%%%%%%%%%%%%%%%%%
\noindent (viii) \textit{Convergence of $I_9^{N,T}$:} We have the equation \[ I_9^{N,T} = \frac{2\lVert \b{Q}_i \rVert}{T}\int_0^T \left ( \huio(t)- \uio(t)   \right )\t \left (  \uio(t) \right) dt.\] Applying the CS Inequality we obtain \begin{equation*} 
I_9^{N,T} \leq 2\lVert \b{Q}_i \rVert \left ( \frac{1}{T}\int_0^T \left \lVert \huio(t) - \uio(t) \right \rVert^2 dt \right )^{1/2}  \left ( \frac{1}{T}\int_0^T \left\lVert \uio(t)  \right\rVert^2 dt \right)^{1/2}
 = : 2\lVert \b{Q}_i \rVert I_{91}^{N,T}\times I_{92}^T.
\end{equation*} 
It is shown in Proposition \ref{prop:pop_u_conv} that 
$\limsup_{T\conv\infty} I_{91}^{N,T} = O(\e_1(N))$ w.p.1. We obtain
$\limsup_{T\conv \infty} I_{92}^T \leq K_2$ w.p.1 as shown in Lemma \ref{lem:stability}.
Therefore, 
\begin{align*} \limsup_{T\conv \infty} I_9^{N,T} &\leq 2\lVert \b{Q}_i \rVert ( \limsup_{T\conv\infty} I_{91}^{N,T} ) \times  ( \limsup_{T\conv\infty} I_{92}^T) \\&= O(\e_1(N)) \quad \text{w.p.1}.
\end{align*} 
Hence, $\limsup_{T\conv\infty} I_9^{N,T} = O(\e_1(N))$ w.p.1.
%\end{enumerate}

Overall we have shown that \[\limsup_{T\conv\infty} I^{N,T} \leq \limsup_{T\conv\infty} ( I_2^T + I_8^T ) + O(\e_1(N)) \quad \text{w.p.1}.\] 
Using the same decomposition technique applied in \eqref{eqn:eq_costs_d} we also show that \[ J_i^N(u_i^0,u_{-i}^0) \leq J_i^N (\huio,\humio) +  O(\e_1(N)) \quad \text{w.p.1}.\] 
Consequently, \[\limN J_i^N ( \huio,\humio) = \limN J_i^N(u_i^0,u_{-i}^0 ) \, \text{w.p.1}, \, 1\leq i \leq N.\]
\IEEEQED
%===============================================================================
%%%%%%%%%%%%%%%%%%%%%%%%%%%%%%%%%%%%%%%%%%%%%%%%%%%%%%%%%%%%%%%%%%%%%%%%%%%%%%%%
% Proof of Proposition prop:eq_costs_g
%%%%%%%%%%%%%%%%%%%%%%%%%%%%%%%%%%%%%%%%%%%%%%%%%%%%%%%%%%%%%%%%%%%%%%%%%%%%%%%%
\subsection*{Proof of Proposition \ref{prop:eq_costs_g}}\label{prop:eq_costs_g_pro}
Let $u_i\triangleq u_i(t;\th_i^0,\z^0)\in\U_g^N$ be a feedback control action and $x_i\triangleq x_i(t;\th_i^0,\z^0)$ be the corresponding closed loop solution. LHS of \eqref{eqn:eq_costs_g} is written as 
\begin{equation*} 
J_i^N(u_i,u_{-i}^0)  =  \lsT\frac{1}{T}\int_0^T \Big\{ \left \lVert x_i(t) -  \mnu(t) \right\rVert_{Q_i}^2  + \left\lVert u_i(t) \right\rVert_R^2 \Big\}dt,
\end{equation*} 
where $\mnu(t)\teq m(x_i(t;\th_i^0,\z^0),x_{j\neq i}^0(t;\th^{1:N},\z^0))$. By adding and subtracting\\ $\hmnu\triangleq m(x_i^0(t;\th_i^0,\z^0),\hat{x}_{j\neq i}^0(t;\th^{[1:N]},\hz^{[1:N]}))$ to the integrand, we get
\begin{equation}\label{eqn:eq_costs_g_2}
J_i^N(u_i,u_{-i}^0)  = \lsT\frac{1}{T}\int_0^T
\Big\{  \lVert x_i(t) -  \hmnu(t)  + \hmnu(t)  - \mnu(t) \rVert_{Q_i}^2 +
 \left\lVert u_i(t) \right\rVert_R^2 \Big\}dt.
\end{equation}
Expanding \eqref{eqn:eq_costs_g_2} , we get
\begin{equation}\label{eqn:costs_g_d}
\begin{aligned}
J_i^N(u_i,u_{-i}^0) & = \limsup_{T\conv\infty} \bigg \{ \frac{1}{T} \int_0^T  \left\lVert  x_i(t) -  \hmnu(t) \right \rVert_{Q_i}^2 dt \\
& \nqquad+ \frac{1}{T} \int_0^T  \left\lVert \hmnu(t) - \mnu(t) \right\rVert_{Q_i}^2 dt \\%
& \nqquad+ \frac{2}{T}\int_0^T \left ( x_i(t) -  \hmnu(t)  \right )\t \b{Q}_i \big(   \hmnu(t) \\
& \nqquad - \mnu(t)    \big)dt + \frac{1}{T} \int_0^T  \left\lVert u_i(t) \right\rVert_R^2 dt \bigg\}\\
& =: \lsT \{I_1^{N,T} + I_2^{N,T} + I_3^{N,T} + I_4^T\}\\
J_i^N(u_i,u_{-i}^0) &  \leq \lsT \{ I_1^{N,T} + I_4^T \} + \lsT I_2^{N,T} \\
& \nqquad + \lsT I_3^{N,T}    \quad \text{w.p.1}.
\end{aligned}
\end{equation}
We have $\lsT \{ I_1^{N,T}  +  I_4^T \} =  J_i^N(u_i,\humio)$; therefore, 
\be\label{eqn:eq_costs_g_JiN}
J_i^N(u_i,u_{-i}^0) \leq  J_i^N(u_i,\humio) +\lsT  I_2^{N,T} + \lsT I_3^{N,T}
\ee
w.p.1.
%\begin{enumerate}[(i)]
%%%%%%%%%%%%%%%%%%%%% (i) %%%%%%%%%%%%%%%%%%%%%%%%%%%%%

\noindent (i) \textit{Convergence of $I_2^{N,T}$:} 
% We have $I_2^{N,T} = \frac{1}{T} \int_0^T  \left\lVert \hmnu(t) - \mnu(t) \right\rVert_{Q_i}^2 dt$. 
Lemma \ref{lem:mass_conv} states that \[ \lsT I_2^{N,T} = O((\e_1(N))^2) \quad \text{w.p.1}.\]
%%%%%%%%%%%%%%%%%%%% (ii) %%%%%%%%%%%%%%%%%%%%%%%%%%%%%

\noindent (ii) \textit{Convergence of $I_3^{N,T}$}:
% We have the equation
% \be
% I_3^{N,T} = \frac{2}{T}\int_0^T \left ( x_i(t) -  \hmnu(t) \right )\t \b{Q}_i \left(   \hmnu(t) - \mnu(t)    \right)dt .
% \ee
Applying the CS Inequality we obtain,
\begin{align*}
I_3^{N,T} & \leq 2\lVert \b{Q}_i \rVert \left ( \frac{1}{T}\int_0^T \left \lVert x_i(t) -  \hmnu(t)   \right \rVert^2 dt \right )^{1/2} \\
& \qquad \left ( \frac{1}{T}\int_0^T \left\lVert  \hmnu(t) - \mnu(t)   \right\rVert^2 dt \right)^{1/2}\\
& =: 2\lVert \b{Q}_i \rVert I_{31}^{N,T}\times I_{32}^{N,T}.
\end{align*}
Using Lemma \ref{lem:con_x_m} we obtain ${\limsup_{T\conv\infty} I_{31}^{N,T} \leq K_4}$ w.p.1 and using Lemma \ref{lem:mass_conv}, we get \[\lsT I_{32}^{N,T} = O(\e_1(N)) \quad \text{w.p.1}.\] 

Therefore, 
\begin{align*} 
\limsup_{T\conv \infty} I_3^{N,T} &\leq 2\lVert \b{Q}_i \rVert ( \limsup_{T\conv\infty} I_{31}^{N,T} ) \times  ( \limsup_{T\conv\infty} I_{31}^{N,T})\\ &= O(\e_1(N)).
\end{align*} 
Hence, $\limsup_{T\conv\infty} I_3^{N,T} = O(\e_1(N))$.
% ACK, 28 Mar 2011
% \ack{check K4}
%\end{enumerate}

Repeating \eqref{eqn:eq_costs_g_JiN} here for ease of reference we see that 
\[ 
J_i^N(u_i,u_{-i}^0) \leq  J_i^N(u_i,\humio) + \limsup_{T\conv \infty} \left( I_2^{N,T} + I_3^{N,T} \right)\text{ w.p.1},
\] 
where $\limsup_{T\conv \infty} \left( I_2^{N,T} + I_3^{N,T} \right)=O(\e_1(N))$. Hence, $J_i^N(u_i,u_{-i}^0) \leq  J_i^N(u_i,\humio) + O(\e_1(N))$ w.p.1. Applying the decomposition technique in \eqref{eqn:costs_g_d} for $J_i^N(u_i,\humio)$, one can also get $J_i^N(u_i,\humio) \leq J_i^N(u_i,u_{-i}^0) + O(e_1(N))$ w.p.1, which implies the claim that $\limN J_i^N(u_i,u_{-i}^0) = \limN J_i^N(u_i,\humio )$ w.p.1, $1\leq i \leq N$. Therefore, 
\begin{equation} 
\limN \inf_{u_i\in\U_g^N}J_i^N(u_i,\humio) = \limN \inf_{u_i \in \U_g^N} J_i^N(u_i,u_{-i}^0) \quad \text{ w.p.1, }1\leq i \leq N.
\end{equation}
\IEEEQED

%===============================================================================
%%%%%%%%%%%%%%%%%%%%%%%%%%%%%%%%%%%%%%%%%%%%%%%%%%%%%%%%%%%%%%%%%%%%%%%%%%%%%%%%
% Proof of Theorem thm:CESACMF
%%%%%%%%%%%%%%%%%%%%%%%%%%%%%%%%%%%%%%%%%%%%%%%%%%%%%%%%%%%%%%%%%%%%%%%%%%%%%%%%
\subsection*{Proof of Theorem \ref{thm:CESACMF}}\label{thm:CESACMF_pro}
%v3
%The proof consists of the unification of the principal theorems proved earlier, the Propositions \ref{prop:eq_costs} and \ref{prop:eq_costs_g} above and the key lemmas established in Appendix \ref{sec:app_MT}.

First, it is evident that Theorem \ref{thm:str_con_est} gives (a), Theorem \ref{thm:scmle} gives (b), and Theorem \ref{thm:L2_stable} gives (c). Second, using a technique similar to that used in \cite[Theorem 6.2]{2008LZ_TAC}, it is shown in Proposition \ref{prop:eq_costs} that
\begin{equation}\label{eqn:CES_J1}
J_i^N(\huio,\humio) \leq J_i^N(u_i^0,u_{-i}^0 ) + O( \e_1 (N) ) \quad \text{w.p.1},
\end{equation}
where $\e_1(N) \conv 0$ as $N \conv \infty$. Then, Lemma \ref{lem:lim_Ji_N} gives
\begin{equation}\label{eqn:CES_J2}
J_i^N(u_i^0,u_{-i}^0) \leq \inf_{u_i \in \U_g^N} J_i^N(u_i,u_{-i}^0) + o(\e_2(N)) + O(\e_2(N))  + O(1/N) \quad \text{w.p.1},\quad 1\leq i \leq N,
\end{equation}
where $\e_2(N) \conv 0$ as $N \conv \infty$. Finally, Proposition \ref{prop:eq_costs_g} states that
\be\label{eqn:CES_J3}
\inf_{u_i \in \U_g^N} J_i^N(u_i,u_{-i}^0) \leq \inf_{u_i \in \U_g^N} J_i^N(u_i,\humio) + O(\e_1(N)),
\ee
w.p.1, $1\leq i \leq N$, where $\e_1(N) \conv 0$ as $N \conv \infty$.

Equations \eqref{eqn:CES_J1}, \eqref{eqn:CES_J2}, \eqref{eqn:CES_J3} together then give the first inequality in
\begin{equation*}%\label{eqn:LZ}
J_i^N( \huio,\humio ) - \e(N) \leq \inf_{u_i \in \U_g^N} J_i^N(u_i, \humio) \leq \\ J_i^N(\huio,\humio),
\end{equation*}
w.p.1, $1\leq i \leq N$, while the second is immediate, where $\e(N)= O (\e_1(N)) + O(\e_2(N)) + o(\e_2(N)) + O(1/N) $. This concludes the proof for (d). 

Claim (e) restates Proposition \ref{prop:eq_costs}, and claim (f) is a consequence of the $\e$-Nash property (d), with the existence of the limits given by \eqref{eqn:eq_costs}.%  It is proved in (d) that $\e\conv\ 0$ as $N\conv\infty$; therefore, (f) follows. 
\IEEEQED

\end{document}